\documentclass[11pt,a4paper,reqno]{amsart}
\usepackage[margin=25mm]{geometry}
\usepackage[T1]{fontenc}
\usepackage{lmodern}
\usepackage{microtype}
\usepackage{amsmath,amssymb,mathtools}
\usepackage{mathrsfs}
\usepackage{enumitem}
\usepackage{booktabs}
\usepackage{xcolor}
\usepackage{tikz}
\usetikzlibrary{arrows.meta,shapes.geometric}
\usepackage[colorlinks=true,allcolors=blue]{hyperref}
\usepackage{bookmark}

\hypersetup{
  pdftitle={Spectral reciprocity for reductive groups},
  pdfsubject={Graded Fourier reciprocity for finite subgroups of SU(2) and reductive targets}
}

\theoremstyle{plain}
\newtheorem{theorem}{Theorem}[section]
\newtheorem{proposition}[theorem]{Proposition}
\newtheorem{lemma}[theorem]{Lemma}
\newtheorem{corollary}[theorem]{Corollary}
\newtheorem{conjecture}[theorem]{Conjecture}
\theoremstyle{definition}
\newtheorem{definition}[theorem]{Definition}
\newtheorem{example}[theorem]{Example}

\theoremstyle{remark}
\newtheorem{remark}[theorem]{Remark}

\numberwithin{equation}{section}

\DeclareMathOperator{\Hom}{Hom}
\DeclareMathOperator{\Lie}{Lie}
\DeclareMathOperator{\Sym}{Sym}
\DeclareMathOperator{\im}{im}
\DeclareMathOperator{\Stab}{Stab}
\DeclareMathOperator{\Ad}{Ad}
\DeclareMathOperator{\ad}{ad}
\DeclareMathOperator{\rootht}{ht}
\DeclareMathOperator{\SO}{SO}

\newcommand{\SU}{\operatorname{SU}}
\newcommand{\UU}{\operatorname{U}}
\newcommand{\CC}{\mathbb C}
\newcommand{\ZZ}{\mathbb Z}
\newcommand{\RR}{\mathbb R}
\newcommand{\FF}{\mathbb F}
\newcommand{\PP}{\mathbb P}
\newcommand{\cO}{\mathcal O}
\newcommand{\fg}{\mathfrak g}
\newcommand{\ft}{\mathfrak t}
\newcommand{\fl}{\mathfrak l}
\newcommand{\fh}{\mathfrak h}
\newcommand{\fa}{\mathfrak a}
\newcommand{\cE}{\mathcal E}
\newcommand{\cP}{\mathcal P}
\newcommand{\cR}{\mathcal R}
\newcommand{\cH}{\mathcal H}
\newcommand{\cM}{\mathcal M}
\newcommand{\wtE}{\widetilde{\mathcal E}}
\newcommand{\wtP}{\widetilde{\mathcal P}}
\newcommand{\wtR}{\widetilde{\mathcal R}}
\newcommand{\pair}[2]{\langle #1,#2\rangle}

\DeclareMathOperator{\coker}{coker}
\DeclareMathOperator{\Tor}{Tor}
\DeclareMathOperator{\Fun}{Fun}
\DeclareMathOperator{\id}{id}
\newcommand{\Dic}{\operatorname{Dic}}
\newcommand{\rh}{\varrho^\vee}
\newcommand{\eexp}[1]{\exp(2\pi i #1)}

\title[Spectral reciprocity for reductive groups]{Spectral reciprocity for reductive groups}

\author{Baizhou Jin}
\address{Baizhou Jin, Qiuzhen College, Tsinghua University, Beijing, China}
\email{jinbz21@mails.tsinghua.edu.cn}

\author{Zichang Wang}
\address{Zichang Wang, Qiuzhen College, Tsinghua University, Beijing, China}
\email{w-zc21@mails.tsinghua.edu.cn}

\date{November 2026}

\begin{document}

\begin{abstract}
For every finite subgroup $\Gamma\subset\SU(2)$, every connected compact or connected complex reductive group $G$, and every finite central subgroup $B\subset Z(G)$, we construct a degree-preserving linear Fourier isomorphism between the polynomial class modules of $(G,B)$ and $((G/B)^\vee,\widehat B)$. The transform exchanges central twists with characters of the actual grades and preserves the invariant polynomial contributions of the full centralizers. In degree zero it yields central-coefficient reciprocity under the Poincar\'e pairing on $S^3/\Gamma$. 

In particular, $\Gamma$ admits equally many conjugacy classes of homomorphisms into $G$ and its Langlands dual, which verifies a conjecture of Kojima and Tachikawa.
\end{abstract}
\maketitle
\tableofcontents

\clearpage
\section{Introduction}
\subsection{Main result}
Let $\Gamma\subset\SU(2)$ be a finite subgroup and let $G$ be a connected compact or connected complex reductive group. Write
\[
N_\Gamma(G)=N(\Gamma,G)=|\Hom(\Gamma,G)/G|
\]
for the number of $G$-conjugacy classes of homomorphisms from $\Gamma$ into $G$.

Throughout this paper, the full root datum $(X,\Phi,Y,\Phi^\vee)$ of $G$ is equipped with its perfect integral pairing $\pair{\cdot}{\cdot}\colon X\times Y\to\ZZ$, and the dual group $G^\vee$ is defined by the dual root datum $(Y,\Phi^\vee,X,\Phi)$. Kojima and Tachikawa formulated the following conjecture for compact simple Lie groups~\cite[Conjecture~1.2]{kojima-tachikawa}.

\begin{conjecture}[Kojima--Tachikawa]\label{Conj.num}
Let $\Gamma\subset\SU(2)$ be finite, and let $G$ be a connected compact simple Lie group with Langlands dual $G^\vee$. Then
\[
N_\Gamma(G)=N_\Gamma(G^\vee).
\]
\end{conjecture}

To state the refinement, fix a finite central subgroup $B\subset Z(G)$ and put
\[
G^\dagger=(G/B)^\vee,\qquad \widehat B=\ker(G^\dagger\to G^\vee)\simeq\Hom(B,\UU(1)).
\]
All group cohomology below has trivial coefficient action. For each $\omega\in H^2(\Gamma,B)$, choose a normalized cocycle $c_\omega$ representing it and define
\[
\begin{split} \mathscr R_{B,\omega}(\Gamma,G)= \bigl\{f\colon\Gamma\to G:\ &f(1)=1,\\ &f(\gamma)f(\delta)=c_\omega(\gamma,\delta)f(\gamma\delta) \text{ for all }\gamma,\delta\in\Gamma\bigr\}/G. \end{split}
\]
Here $G$ acts by simultaneous conjugation. Set
\[
V_B(\Gamma,G)=\bigoplus_{\omega\in H^2(\Gamma,B)} \CC[\mathscr R_{B,\omega}(\Gamma,G)], \qquad F(\Gamma;B)=H^1(\Gamma,B)\times\widehat{H^2(\Gamma,B)}.
\]
The action on a basis vector in the $\omega$-summand is
\[
(\eta,\xi)\cdot[f]=\xi(\omega)[\eta f], \qquad (\eta f)(\gamma)=\eta(\gamma)f(\gamma).
\]
Changing the chosen cocycles gives isomorphic $F(\Gamma;B)$-modules. These are the twisted class modules of~\cite[Definitions~1.5 and~1.7]{kojima-tachikawa}.

Using $H^i(\Gamma,C)\simeq H^i(S^3/\Gamma,C)$ for $i=1,2$, the cup product, coefficient evaluation, and the oriented fundamental class of $S^3/\Gamma$ give a perfect pairing and hence an isomorphism
\[
\iota_B\colon H^1(\Gamma,B)\xrightarrow{\sim} \widehat{H^2(\Gamma,\widehat B)}.
\]
The resulting swap
\[
\begin{split} s_B\colon F(\Gamma;B)&\xrightarrow{\sim}F(\Gamma;\widehat B),\\ (\eta,\xi)&\longmapsto \bigl(\iota_{\widehat B}^{-1}(\xi),\iota_B(\eta)\bigr) \end{split}
\]
is the spherical Poincar\'e swap of~\cite[Definition~1.8]{kojima-tachikawa}. The refined conjecture~\cite[Conjecture~1.9]{kojima-tachikawa} asserts the following.

\begin{conjecture}[Kojima--Tachikawa]\label{Conj.fou}
Let $\Gamma\subset\SU(2)$ be finite, let $G$ be a connected compact simple Lie group, and let $B\subset Z(G)$ be finite. With $G^\dagger=(G/B)^\vee$, there is a complex-linear isomorphism
\[
\mathcal U_B\colon V_B(\Gamma,G)\xrightarrow{\sim} V_{\widehat B}(\Gamma,G^\dagger)
\]
satisfying $\mathcal U_B(a\cdot v)=s_B(a)\cdot\mathcal U_B(v)$ for every $a\in F(\Gamma;B)$ and $v\in V_B(\Gamma,G)$.
\end{conjecture}

Djokovi\'c~\cite[Theorem~1]{djokovic} proved the cyclic case $\Gamma=C_n$, which counts conjugacy classes of elements whose orders divide $n$ in complex semisimple groups. Kojima and Tachikawa~\cite{kojima-tachikawa} proved Conjecture~\ref{Conj.num} for $(\SU(n),\operatorname{PU}(n))$ and $(\operatorname{Sp}(n),\SO(2n+1))$ for every finite subgroup of $\SU(2)$, and for $(\operatorname{PSp}(n),\operatorname{Spin}(2n+1))$ for $\Gamma\in\{2T,2O,2I\}$. They also proved several cases of Conjecture~\ref{Conj.fou}, including cyclic groups with arbitrary central coefficients and the full-center $\SU(n)$ case for every such $\Gamma$.

We now consider the projective representation class, which is a graded refinement for connected compact or connected complex reductive $G$, with the same notation $B$, $G^\dagger$, and $\widehat B$. The detailed definitions are in \eqref{eq:classes}, Definition~\ref{def:projective-modules}(1) and \eqref{eq:cyclic-classes}. We write $\mathscr P_{\Gamma,B}(G)$ for the polynomial class module, whose degree-zero part is the space of complex-valued functions on the corresponding actual projective class set:
\begin{center}
\begin{tabular}{ccc}
\toprule
$\Gamma$ & $\mathscr P_{\Gamma,B}(G)$ & $r_\Gamma$\\
\midrule
$\Dic_m\ (m\ge2)$ & $\cP_B(G)$ & $2$\\
$2T,2O,2I$ & $\wtP_{\Gamma,B}(G)$ & $1$\\
$C_n\ (n\ge1)$ & $\cP_B(C_n,G)$ & $1$\\
\bottomrule
\end{tabular}
\end{center}
The three modules are defined in \eqref{eq:gradeddef}, Definition~\ref{def:projective-modules}(3), and \eqref{eq:cyclic-classes-mods}, respectively. In each case a class contributes $\Sym(\mathfrak s^*)^{W_{\mathrm{full}}}$, where $\mathfrak s$ is the complex Lie algebra of a maximal torus of the connected centralizer (complexified when $G$ is compact), and $W_{\mathrm{full}}$ is the action of its normalizer in the full centralizer. Linear forms have degree one. The actual grade belongs to $B^{r_\Gamma}$: for $\Dic_m$, the group $H^2(\Gamma,\ZZ)$ is $(\ZZ/2\ZZ)^2$ when $m$ is even and $\ZZ/4\ZZ$ when $m$ is odd, so keeping the pair of relator grades in $B^2$ treats both parities uniformly, whereas for $2T,2O,2I$ and $C_n$ the group $H^2(\Gamma,\ZZ)$ is cyclic and a single grade in $B$ suffices. The central translations $T_k$ and grade multipliers $M_\chi$ are given by \eqref{eq:operators}, \eqref{eq:heisenberg-ops}, and \eqref{eq:cyclic-operators}. Our main result is the following.

\begin{theorem}\label{thm:unified}
Let $\Gamma\subset\SU(2)$ be finite, let $G$ be connected compact or connected complex reductive, and let $B\subset Z(G)$ be finite.
\begin{enumerate}[label=\textup{(\arabic*)},leftmargin=2em]
\item There is a complex-linear, degree-preserving Fourier isomorphism
\[
\mathcal F_{\Gamma,B}\colon\mathscr P_{\Gamma,B}(G) \xrightarrow{\sim}\mathscr P_{\Gamma,\widehat B}(G^\dagger)
\]
such that, for $k\in B^{r_\Gamma}$ and $\chi\in\widehat{B^{r_\Gamma}}$,
\[
\mathcal F_{\Gamma,B}T_k=M_k^\dagger\mathcal F_{\Gamma,B}, \qquad \mathcal F_{\Gamma,B}M_\chi=T_{\chi^{-1}}^\dagger\mathcal F_{\Gamma,B}.
\]
Here $k$ is evaluated as a character of $\widehat B^{\,r_\Gamma}$; in additive notation the second translation is $T_{-\chi}^\dagger$. The transform preserves the polynomial invariants of the full centralizers, and its degree-zero restriction compares the function spaces on the actual projective class sets.
\item In particular, the ordinary representation counts satisfy
\[
|\Hom(\Gamma,G)/G|=|\Hom(\Gamma,G^\vee)/G^\vee|.
\]
\item Choosing one actual grade for each class in $H^2(\Gamma,B)$ and taking polynomial degree zero gives an $F(\Gamma;B)$-module isomorphic to $V_B(\Gamma,G)$. There is an $s_B$-equivariant complex-linear isomorphism
\[
V_B(\Gamma,G)\xrightarrow{\sim}V_{\widehat B}(\Gamma,G^\dagger),
\]
where $s_B$ is the spherical Poincar\'e swap defined above.
\end{enumerate}
\end{theorem}

Assertions \textup{(2)} and \textup{(3)} prove Conjectures~\ref{Conj.num} and~\ref{Conj.fou} and extend them to connected reductive groups. The actual-grade Fourier transform in \textup{(1)} is expressed in the chosen carry or relator coordinates. To obtain \textup{(3)}, we pass to cohomology classes and compare those coordinates with the spherical Poincar\'e pairing. The dihedral comparison is Corollary~\ref{cor:poincare-packets}, the polyhedral comparison is Proposition~\ref{prop:packet-reciprocity}, and the cyclic comparison is Corollary~\ref{cor:cyclic-packets}.

The proof of \textup{(1)} for the dihedral groups is completed in Section~\ref{sec:graded}, and the polyhedral proof is completed in Section~\ref{sec:polyhedral-conclusion}; Proposition~\ref{prop:cyclic-fourier} treats the cyclic case. These exhaust the finite subgroups of $\SU(2)$~\cite{coxeter}. For \textup{(2)}, take $B=1$ and polynomial degree zero in \textup{(1)}. Then the two spaces have bases $\Hom(\Gamma,G)/G$ and $\Hom(\Gamma,G^\vee)/G^\vee$, whose finiteness is recorded in Lemma~\ref{lem:compact-reduction}. Section~\ref{sec:count} also gives an independent dihedral counting proof by Burnside's lemma.

\subsection{Compact reduction and finiteness of conjugacy classes}
We first recall why representation counts and full-centralizer polynomial gradings for a connected complex reductive group $G$ may be computed inside a maximal compact subgroup $K\subset G$.

\begin{lemma}\label{lem:compact-reduction}
Let $G$ be a connected complex reductive group with maximal compact subgroup $K$, and let $\Gamma$ be a finite group.
\begin{enumerate}[label=\textup{(\arabic*)},leftmargin=2em]
\item Inclusion induces a bijection $\Hom(\Gamma,K)/K\xrightarrow{\sim}\Hom(\Gamma,G)/G$, and for every $\rho\in\Hom(\Gamma,K)$, the inclusion $C_K(\rho)\hookrightarrow C_G(\rho)$ is a maximal compact subgroup inducing an isomorphism on component groups and Weyl groups.
\item The quotient $\Hom(\Gamma,G)/G$ is finite.
\end{enumerate}
\end{lemma}
\begin{proof}
For any finite subgroup $F\subset G$, the function $x\mapsto\sum_{f\in F}d(x,fo)^2$ on the nonpositively curved symmetric space $G/K$ is proper and strictly convex. Its unique minimum is fixed by $F$, which conjugates $F$ into $K$. If $g\in G$ conjugates two images in $K$, unitarity implies that $g^*g$ centralizes the first image. In the polar decomposition $g=ku$ with $k\in K$ and $u$ positive, the factor $u=(g^*g)^{1/2}$ therefore centralizes the first image as well, so $k\in K$ gives the same conjugation. The same polar-decomposition argument applies to $C_G(\rho)$.

To see finiteness, note that at any $\rho\in\Hom(\Gamma,G)$, the Zariski tangent space to the representation variety $\Hom(\Gamma,G)$ is $Z^1(\Gamma,\fg_{\Ad\rho})$, whereas the tangent space to the $G$-orbit $G\cdot\rho$ is $B^1(\Gamma,\fg_{\Ad\rho})$. Since $\Gamma$ is finite, averaging over $\Gamma$ gives $Z^1(\Gamma,\fg_{\Ad\rho})=B^1(\Gamma,\fg_{\Ad\rho})$. Consequently, the inequalities
\[
\dim(G\cdot\rho)\le\dim_\rho\Hom(\Gamma,G)\le\dim T_\rho\Hom(\Gamma,G)
\]
are equalities. Every point $\rho$ is therefore smooth on $\Hom(\Gamma,G)$, and its locally closed orbit contains an open neighborhood of $\rho$; translating by $G$ shows that each orbit is open. Quasi-compactness of the affine variety $\Hom(\Gamma,G)$ implies that there are finitely many orbits.
\end{proof}

The same reduction applies to the projective classes used below: the lifts of their finitely many elements, together with $B$, generate a finite subgroup, since its quotient by $B$ is an image of $\Gamma$. Applying Lemma~\ref{lem:compact-reduction} to this finite group shows that the full normalizer action on a maximal torus agrees for a centralizer and its maximal compact subgroup, so the full-centralizer polynomial modules may be computed in either setting.

Part~\ref{part:dihedral} treats $\Dic_m$ for $m\ge2$, Part~\ref{part:polyhedral} treats $2T,2O,2I$, and Appendix~\ref{sec:cyclic} treats $C_n$.

\subsection{Physical context}

The physical motivation comes from electric--magnetic duality in four-dimensional $\mathcal N=4$ supersymmetric gauge theory; see Vafa and Witten~\cite{vafa-witten}. Gauge holonomy on $S^3/\Gamma$ relates this setting to homomorphisms from $\Gamma$.

Earlier physical studies relate ADE representation counts to Chern--Simons states and Ehrhart polynomials, and analyze supersymmetric Casimir-energy degeneracy on $S^3/\Gamma$~\cite{JuEhrhart,JuCasimir}. Choi and Tachikawa~\cite[\S\S4.3--4.4]{CT} match electric and magnetic symmetry actions with Chern--Simons line operators and use the modular $S$-matrix to recover the $\SU(n)$--$\operatorname{PU}(n)$ counting equality. A related correspondence for ground states of type $A_1$ class $\mathcal S$ theories is studied in~\cite{AVGJ}, with the Chern--Simons gauge group determined by $\Gamma$ through the McKay correspondence. Homomorphisms $\Gamma\to E_8$ and their physical induction have also been studied through six-dimensional superconformal theories~\cite{FR,BZ}.

\subsection{Acknowledgements}
The authors would like to thank Shihao Wang and Zhiwei Zheng for many helpful discussions. This work is supported by the National Key R\&D Program of China (No. 2024YFA1014700).

\subsection{AI disclosure}
Most key arguments were originally developed through multiple rounds of interaction with ChatGPT 6.0 Astra, and ChatGPT 6.1 Sol and Gemini 4 Neon were used for writing and editorial refinement. The authors independently verified and finalized all mathematical content and take full responsibility for the paper.

\clearpage
\part{Binary dihedral groups}\label{part:dihedral}

\section{Main results and outline}

In this part $\Gamma=\Dic_m$, $m\ge2$. We use the two relator grades in $B^2$ to prove the dihedral cases of Theorem~\ref{thm:unified}, obtaining a finite abelian model uniformly in $m$ and in the full root datum of $G$.

\subsection{Main results}

Let
\begin{equation}\label{eq:presentation}
\Gamma=\Dic_m=\langle a,b\mid a^m b^{-2}=1,\ (bab^{-1}a)^{-1}=1\rangle, \qquad m\ge2.
\end{equation}
The relations give $a^{2m}=1$ and reduce every word to $a^j b^\epsilon$ with $0\le j<2m$ and $\epsilon\in\{0,1\}$. The quaternion realization, $a=\exp(\pi i/m)$ and $b=j$, has these $4m$ distinct elements, so the group has order $4m$. Write
\[
N_m(G)=|\Hom(\Gamma,G)/G|.
\]
By Lemma~\ref{lem:compact-reduction}, throughout this part reductive groups are assumed to be connected compact.

The Langlands dual $G^\vee$ is defined by the full dual root datum. In particular, we work with the full character and cocharacter lattices and allow nontrivial central quotients and connected centers.

To formulate the Fourier comparison, let $B\subset Z(G)$ be a finite central subgroup and define
\begin{equation}\label{eq:classes}
\cR_B(G)=\{(t,u)\in G^2:t^m u^{-2}\in B,\ (utu^{-1}t)^{-1}\in B\}/G.
\end{equation}
Its actual relator grade is
\[
\gamma(t,u)=(t^m u^{-2},(utu^{-1}t)^{-1})\in B^2.
\]
We use additive notation for tori and finite coefficient groups. Put
\begin{equation}\label{eq:D}
D=\begin{pmatrix}m&-2\\-2&0\end{pmatrix},\qquad G^\dagger=(G/B)^\vee,\qquad B^\dagger=\widehat B=\Hom(B,\UU(1))\subset Z(G^\dagger),
\end{equation}
where $B^\dagger=\ker(G^\dagger\to G^\vee)$. For a representative $p=(t,u)$ of a class in $\cR_B(G)$, let $D_p=C_G(t,u)$, choose a maximal torus $S_p\subset D_p^0$, and let $W_p$ be the image of the full normalizer $N_{D_p}(S_p)$ on $\mathfrak s_p=\Lie(S_p)_{\CC}$. Define the polynomial class module
\begin{equation}\label{eq:gradeddef}
\cP_B(G)=\bigoplus_{[p]\in\cR_B(G)} \Sym(\mathfrak s_p^*)^{W_p}.
\end{equation}
These invariant polynomial spaces are canonically identified under changes of maximal torus and conjugation of $p$, since the ambiguity lies in the full normalizer. We grade $\cP_B(G)$ by polynomial degree, with linear forms of degree one; its degree-zero part is canonically $\Fun(\cR_B(G),\CC)$.

For $k=(k_1,k_2)\in B^2$, let $\tau_k(t,u)=(k_1t,k_2u)$. Central twisting leaves $D_p$ unchanged and hence identifies the corresponding polynomial spaces. Using these identifications, define operators on $\cP_B(G)$ by
\begin{equation}\label{eq:operators}
(T_k f)(p)=f(\tau_{-k}p),\qquad (M_\chi f)(p)=\chi(\gamma(p))f(p),\qquad \chi\in\widehat{B^2}.
\end{equation}
Here $f(p)$ denotes the invariant polynomial attached to the class of $p$; in degree zero these are the usual operators on scalar-valued functions.

In these coordinates, Theorem~\ref{thm:unified}\textup{(1)} asks for a degree-preserving map $\mathcal F_B\colon\cP_B(G)\xrightarrow{\sim}\cP_{B^\dagger}(G^\dagger)$ satisfying
\begin{equation}\label{eq:fourier-rules}
\mathcal F_B T_k=M_k^\dagger\mathcal F_B,\qquad \mathcal F_B M_\chi=T_{-\chi}^\dagger\mathcal F_B.
\end{equation}
Its degree-zero restriction is the map $\mathcal F_{B,0}$ between function spaces on the actual class sets. For $B=1$, these sets are $\Hom(\Gamma,G)/G$ and $\Hom(\Gamma,G^\vee)/G^\vee$, so their dimension equality gives Theorem~\ref{thm:unified}\textup{(2)}; an independent proof by Burnside's lemma is given in Section~\ref{sec:count}.

\subsection{Proof outline}

In Section~\ref{sec:normal}, we construct an adapted maximal torus $T=C_{C_G(t)^0}(S)$ for any twisted dihedral pair $(t,u)$ by choosing a maximal torus $S\subset C_G(t,u)^0$. This puts the pair in torus normal form relative to an involution $w\in W$, characterizes adaptation by root regularity in $\Psi=\{\alpha\in\Phi:w\alpha=-\alpha\}$, and identifies the full centralizer action on $\Lie(S)$ with the Weyl stabilizer in $C_W(w)$.

In Section~\ref{sec:contraction}, we introduce an integral torus map $F_w$ whose kernel $K_{B,w}$ may have positive-dimensional components. An affine contraction on the fundamental alcoves of $\Psi$ linearizes the Tits shift $\rh$ and induces a $C_W(w)/W(\Psi)$-equivariant strong deformation retraction onto the adapted orbits, identifying $\cR_B(G)$ with $\coprod_{[w]\,\mathrm{admissible}} E_{B,w}/C_W(w)$, where $E_{B,w}=\pi_0 K_{B,w}$.

In Section~\ref{sec:lattice}, we establish equivariant Pontryagin duality for the component groups of possibly singular lattice maps. Passing from $(G,B)$ to $(G^\dagger,B^\dagger)$ transposes $F_w$, identifying $E_{B^\dagger,w}(G^\dagger)$ with $\widehat{E_{B,w}(G)}$ as $C_W(w)$-modules.

In Section~\ref{sec:count}, Burnside's lemma applied to the dual $C_W(w)$-modules gives an independent proof of the ordinary dihedral counting reciprocity.

In Sections~\ref{sec:fourier} and~\ref{sec:graded}, we construct the finite Fourier transform on each admissible stratum and tensor it with $C_W(w)$-equivariant polynomial isomorphisms on the $w$-fixed Cartan subspaces, establishing the degree-zero and graded dihedral cases of Theorem~\ref{thm:unified}\textup{(1)}.

In Section~\ref{sec:packets}, we pass from relator coordinates to cohomological packets and deduce packet reciprocity in evaluation coordinates; the comparison with the spherical Poincar\'e pairing on $S^3/\Dic_m$ is carried out in Appendix~\ref{sec:poincare}. Section~\ref{sec:examples} illustrates the construction on tori, $Q_8$, and $\SU(3)$.

\section{Adapted tori and Weyl conjugacy}\label{sec:normal}

We use the compact Lie group and maximal torus structure theory in~\cite[Ch.~IV--V]{brocker-tomdieck}, and the root-system and affine-Weyl-group background in~\cite[Ch.~IV--VI]{Bourbaki}.

Fix a maximal torus $T=V/\Lambda$, with exponential $x\mapsto\eexp{x}$, full cocharacter lattice $\Lambda=X_*(T)$, and dual character lattice $\Lambda^*=X^*(T)$. Let $\Phi\subset V^*$ be the root system and $W=N_G(T)/T$ the Weyl group. Every root is integral on $\Lambda$, and $W$ fixes $Z(G)$ pointwise.

For the coefficient subgroup $B$, define
\begin{equation}\label{eq:coefficientlattice}
\Lambda_B=\{v\in V:\eexp{v}\in B\},\qquad \Lambda_B/\Lambda=B,\qquad T/B=V/\Lambda_B.
\end{equation}
Every root is integral on $\Lambda_B$, since $B$ is central. Torus coordinates $x,y$ are additive; a tilde denotes a real lift when one is needed. Expressions such as $\alpha(x)\notin\ZZ$ are independent of the lift.

\subsection{Adapted maximal tori}

We construct normal forms by choosing a maximal torus containing $t$ and normalized by $u$. Choosing this torus through a maximal torus of the connected pair centralizer also controls conjugacy, since two such choices can be aligned by an element that fixes the pair.

\begin{lemma}\label{lem:adapted}
Let $J$ be a connected compact Lie group, let $\alpha$ be an automorphism with $\alpha^2=1$, and let $S$ be a maximal torus of $(J^\alpha)^0$. Then $C_J(S)$ is a maximal torus of $J$, and all tori obtained this way are conjugate under $(J^\alpha)^0$.
\end{lemma}
\begin{proof}
The group $H=C_J(S)$ is connected and $\alpha$-stable, and $S\subset Z(H)^0$. Write $\fh=\Lie H$ and $\mathfrak s=\Lie S$. Since $S$ is a maximal torus of $(J^\alpha)^0$, we have
\[
\fh^\alpha =\Lie C_{(J^\alpha)^0}(S) =\mathfrak s \subseteq\mathfrak z(\fh).
\]

The compact Lie algebra $\fh$ is reductive, so its derived algebra $\mathfrak d=[\fh,\fh]$ is semisimple and has zero intersection with $\mathfrak z(\fh)$. It is $\alpha$-stable because
\[
\alpha[\fh,\fh]=[\alpha\fh,\alpha\fh]=[\fh,\fh],
\]
so $\mathfrak d^\alpha=\fh^\alpha\cap\mathfrak d=0$. Since $\alpha$ is an involution, it acts as $-1$ on $\mathfrak d$. For $X,Y\in\mathfrak d$, the automorphism identity gives
\[
-[X,Y]=\alpha[X,Y]=[\alpha X,\alpha Y]=[X,Y],
\]
so $\mathfrak d$ is abelian and hence zero. Thus $H$ is a torus. It contains any maximal torus of $J$ containing $S$, so it is maximal. Finally, any two choices of $S$ are conjugate under the connected group $(J^\alpha)^0$, and the same conjugation identifies their centralizers in $J$.
\end{proof}

For a pair $(t,u)$, a maximal torus $T_0\subset G$ containing $t$ and normalized by $u$ is called \emph{adapted to $(t,u)$} if $t$ is regular in the Levi subgroup
\[
L=C_G((T_0^{\Ad u})^0),
\]
or equivalently, if $C_L(t)^0=T_0$. Relative to the fixed torus $T$, we call a pair in torus normal form \emph{adapted} if $T$ is adapted to that pair.

Let $(t,u)$ be a pair with $\gamma(t,u)=(b_1,b_2)\in B^2$, and put $J=C_G(t)^0$. The element $t$ lies in $J$, since any maximal torus of $G$ containing $t$ is contained in $C_G(t)^0$; in particular, $J$ has the same rank as $G$ and $t\in Z(J)$. Also $u$ normalizes $J$, because
\[
utu^{-1}=b_2^{-1}t^{-1}
\]
has the same centralizer as $t$, and the relation $u^2=b_1^{-1}t^m$ makes $\Ad u$ an involution on $J$.

Choose a maximal torus $S$ of $(J^{\Ad u})^0=C_G(t,u)^0$, and put $L=C_G(S)$. Then $L$ is a Levi subgroup of $G$ containing $t$ and $u$, with
\[
C_L(t,u)^0=S=Z(L)^0.
\]
Since $C_L(t)^0=C_J(S)$, Lemma~\ref{lem:adapted} gives
\begin{equation}\label{eq:adaptedtorus}
T=C_J(S)=C_L(t)^0
\end{equation}
as a maximal torus of $J$, hence of $G$ and $L$. It contains $t$, is normalized by $u$, and has $t$ regular in $L$, so $T$ is adapted to $(t,u)$.

\subsection{Root subsystems of involutions}

Fix a $W$-invariant positive inner product on $V$. For an involution $w\in W$, put
\[
V^\pm=\ker(w\mp1),\qquad S=(T^w)^0,\qquad L=C_G(S).
\]
Then $L$ is a Levi subgroup of $G$ with maximal torus $T$. Since $\Lie S=V^+$, the roots of $L$ are precisely those roots of $G$ that vanish on $V^+$. We therefore set
\begin{equation}\label{eq:rootdata}
\begin{gathered}
\Psi=\Phi(L,T) =\{\alpha\in\Phi:\alpha|_{V^+}=0\} =\{\alpha\in\Phi:w\alpha=-\alpha\},\\
W_0=W(L,T)=W(\Psi),\qquad C=C_W(w),\qquad A=1+w.
\end{gathered}
\end{equation}

A generic vector of $V^+$ is annihilated precisely by the roots in $\Psi$, so its Weyl stabilizer is $W_0$. Since $w$ fixes this vector, $w\in W_0$. The coroot span of $\Psi$ is contained in $V^-$. Conversely, $W_0$ fixes the orthogonal complement of that span pointwise. As $w\in W_0$ and $w=-1$ on $V^-$, the coroot span must be all of $V^-$. Thus $w$ acts as $-1$ on the reflection space of $\Psi$ and is its longest Weyl element for any choice of positive roots.

The Weyl centralizer of $w$ is
\begin{equation}\label{eq:normalizer}
C_W(w)=N_W(\Psi).
\end{equation}
Indeed, an element commuting with $w$ preserves $\Psi$; conversely, an element normalizing $\Psi$ preserves its coroot span $V^-$ and its orthogonal complement $V^+$, on which $w$ acts as $-1$ and $1$, respectively, and therefore commutes with $w$. In particular, $W_0$ is a normal subgroup of $C$.

Write $\Psi=\coprod_i\Psi_i$, choose positive roots, and put
\[
\rh_i=\frac12\sum_{\alpha\in\Psi_i^+}\alpha^\vee, \qquad \rh=\sum_i\rh_i.
\]
The pinned compact Tits representative $n=n_w$ can be chosen in the derived group of $L$, where it satisfies~\cite[Lemma~5.1]{adams}
\begin{equation}\label{eq:tits}
n^2=\eexp{\rh}.
\end{equation}
This identity follows by multiplying a reduced expression for the longest element by its reverse: the rank-one identity $n_\alpha^2=\alpha^\vee(-1)$ and the braid relations show that each positive coroot contributes once, giving the exponential of half their sum. It also descends from the simply connected cover of $[L,L]$ to every isogeny form.

Set $T_-=(1-w)T$ and $R=T/T_-$. The torus $T_-$ has Lie algebra $V^-$. Writing $t=\eexp{x}$ and $u=\eexp{y}n$, conjugation by $T$ changes $y$ by exactly $T_-$. The equations become
\begin{equation}\label{eq:beta}
\beta(x,[y])=(mx-Ay-\rh,-Ax)\in B^2.
\end{equation}
This is an equality in $T^2$, with $\rh$ read modulo $\Lambda$. The term $Ay$ is well defined on $R$ because $A(1-w)=0$, so quotienting by torus conjugation removes precisely the $T_-$ direction of the second parameter.

\begin{lemma}\label{lem:criterion} A pair in the above torus normal form is adapted if and only if
\begin{equation}\label{eq:regular}
\alpha(x)\notin\ZZ\qquad(\alpha\in\Psi).
\end{equation}
These conditions are also equivalent to $S=(T^w)^0$ being a maximal torus of $C_G(t,u)^0$. When they hold,
\[
T=C_{C_G(t)^0}(S).
\]
\end{lemma}
\begin{proof}
Put $S=(T^w)^0$ and $L=C_G(S)$. Since $\Phi(L,T)=\Psi$, the roots of $C_L(t)^0$ relative to $T$ are precisely those $\alpha\in\Psi$ for which $\alpha(x)\in\ZZ$. Thus \eqref{eq:regular} is equivalent to $C_L(t)^0=T$, which is the definition of adaptation.

Set $J=C_G(t)^0$ and $D=C_G(t,u)$. Connectedness of the torus centralizer $C_J(S)$ gives
\[
C_J(S)=C_L(t)^0.
\]
Suppose the pair is adapted, so this group is $T$. Any torus $S'\subset D^0$ containing $S$ centralizes both $t$ and $S$, and hence lies in $C_J(S)=T$. Since $u$ centralizes $S'$, we obtain
\[
S'\subseteq(T^{\Ad u})^0=(T^w)^0=S,
\]
so $S$ is a maximal torus of $D^0$.

Conversely, suppose $S$ is a maximal torus of $D^0=(J^{\Ad u})^0$. Lemma~\ref{lem:adapted} shows that $C_J(S)$ is a maximal torus of $J$. It contains $T$, which is itself a maximal torus of $J$, so $C_J(S)=T$. Consequently $C_L(t)^0=T$, and the pair is adapted. The formula for $T$ follows as well.
\end{proof}

Define $X_{B,w}\subset T\times R$ by \eqref{eq:beta} and \eqref{eq:regular}. It is finite: applying $1-w$ to the first equation uses $(1-w)B=0$ and $(1-w)\rh=2\rh\in\Lambda$ to give $m(1-w)x=0$, which together with $Ax=-b_2$ gives
\begin{equation}\label{eq:finite}
2mx=-mb_2.
\end{equation}
There are finitely many possibilities for $x$, and for each one, $A\colon R\to(T^w)^0$ is an isogeny, so the first equation allows finitely many $[y]$.

\subsection{Normalizer action}

Conjugation by a Weyl lift may change the chosen representative $n$ by a translation in the second coordinate. We show that this correction lies in $T_-$ and hence vanishes in $R$, so the induced Weyl action on the parameter space is linear.

\begin{lemma}\label{lem:linearaction} For every $v\in C=C_W(w)$ and every lift $\dot v\in N_G(T)$,
\begin{equation}\label{eq:commutator}
\dot v n\dot v^{-1}n^{-1}\in T_-.
\end{equation}
Consequently the action on $X_{B,w}$ induced by conjugation is
\[
v\cdot(x,[y])=(vx,[vy]).
\]
\end{lemma}
\begin{proof}
First suppose $v\in W_0$, and choose its lift in $[L,L]$. The commutator belongs to $T$, because $v$ commutes with $w$, and also to $[L,L]$. Since $T\cap[L,L]=T_-$ (the maximal torus of $[L,L]$, with Lie algebra $V^-$), this proves \eqref{eq:commutator} for these lifts.

Next suppose $v$ preserves $\Psi^+$. A lift normalizes $L$ and permutes the simple root spaces of its derived group, with scalar factors of absolute value one relative to a compact-compatible pinning. The map
\[
T\longrightarrow\prod_{\alpha\text{ simple in }\Psi}\UU(1), \qquad t\longmapsto(\alpha(t))_\alpha
\]
is surjective, since its differential has full rank by linear independence of the simple roots. Multiplying the lift by a suitable element of $T$ makes it permute the pinned simple root vectors exactly. It then permutes the rank-one Tits representatives, and by the braid relations it fixes their canonical longest-word product $n$, so its commutator with $n$ is the identity.

Every $v\in N_W(\Psi)$ is a product of an element of $W_0$ and one preserving $\Psi^+$. The commutator identity and invariance of $T_-$ prove the claim for a lift of every such $v$. Changing a lift by $t\in T$ changes the commutator by $t\,w(t)^{-1}\in T_-$. Finally, conjugating $\eexp{y}n$ adds precisely this commutator to $vy$, which vanishes in $R$. Since the action is induced by conjugation, it preserves the actual relator grade.
\end{proof}

\subsection{Normal forms and conjugacy classes}

We now show that conjugation by $T$ and $N_G(T)$ accounts for all conjugacies in $G$.

\begin{proposition}\label{prop:normalform} There is a grade-preserving bijection
\begin{equation}\label{eq:normalform}
\cR_B(G)\simeq\coprod_{[w],\,w^2=1} X_{B,w}/C_W(w).
\end{equation}
\end{proposition}
\begin{proof}
Choose one representative $w$ of each conjugacy class of involutions in $W$, together with the representative $n_w$ fixed above. For $(x,[y])\in X_{B,w}$, choose a representative $y\in T$ and form
\[
t=\eexp{x},\qquad u=\eexp{y}n_w.
\]
Equation~\eqref{eq:beta} says precisely that both relator values of this pair lie in $B$. Changing $y$ by an element of $T_-=(1-w)T$ changes the pair by conjugation in $T$, so its $G$-conjugacy class depends only on $(x,[y])$. By Lemma~\ref{lem:linearaction}, the action of $C_W(w)$ on these parameters is induced by conjugation in $N_G(T)$. We therefore obtain a well-defined map
\[
\varphi\colon \coprod_{[w],\,w^2=1}X_{B,w}/C_W(w) \longrightarrow \cR_B(G).
\]

To prove surjectivity, let $(t,u)$ represent a class in $\cR_B(G)$, and choose a maximal torus $S$ of $C_G(t,u)^0$. The adapted-torus construction gives a maximal torus
\[
T_0=C_{C_G(t)^0}(S)
\]
of $G$ adapted to $(t,u)$. Conjugate $T_0$ to the fixed torus $T$. The image of $u$ in $W$ is an involution; conjugating further by an element of $N_G(T)$, we may arrange that this image is our chosen representative $w$. We can then write
\[
t=\eexp{x},\qquad u=\eexp{y}n_w.
\]
The relator equations give \eqref{eq:beta}, and adaptation gives \eqref{eq:regular} by Lemma~\ref{lem:criterion}. Thus $(x,[y])\in X_{B,w}$.

To prove injectivity, suppose that $(x,[y])\in X_{B,w}$ and $(x',[y'])\in X_{B,w'}$ determine pairs $(t,u)$ and $(t',u')$ conjugate by some $g\in G$. Put
\[
D=C_G(t,u),\qquad D'=C_G(t',u'),\qquad S=(T^w)^0,\qquad S'=(T^{w'})^0.
\]
By Lemma~\ref{lem:criterion}, $S$ and $S'$ are maximal tori of $D^0$ and $(D')^0$, respectively, and
\[
T=C_{C_G(t)^0}(S) =C_{C_G(t')^0}(S').
\]

Since $gD^0g^{-1}=(D')^0$, the tori $gSg^{-1}$ and $S'$ are maximal tori of the same connected compact group $(D')^0$. Choose $d\in(D')^0$ such that
\[
d(gSg^{-1})d^{-1}=S',
\]
and set $h=dg$. Because $d$ centralizes the second pair, $h$ still conjugates $(t,u)$ to $(t',u')$. The identities for $T$ above give
\[
hTh^{-1} =C_{C_G(t')^0}(hSh^{-1}) =C_{C_G(t')^0}(S') =T,
\]
so $h\in N_G(T)$.

Let $v=hT\in W$. Passing $huh^{-1}=u'$ to $W$ gives $vwv^{-1}=w'$. Since $w$ and $w'$ were chosen as representatives of conjugacy classes, $w=w'$ and $v\in C_W(w)$. Lemma~\ref{lem:linearaction} then yields
\[
(x',[y'])=(vx,[vy]),
\]
so the two parameters represent the same point of $X_{B,w}/C_W(w)$.

Finally, the grade of the pair associated to $(x,[y])$ is $\beta(x,[y])$, and conjugation fixes both relator values because they lie in the central subgroup $B$. Hence $\varphi$ preserves the actual relator grade.
\end{proof}

The classification is therefore
\[
\begin{aligned}
\text{pairs modulo }G &\longleftrightarrow\coprod_{[w]}X_{B,w}/C\\
&\longleftrightarrow\coprod_{[w]}(X_{B,w}/W_0)/(C/W_0),
\end{aligned}
\]
with $C$ and $W_0$ taken in each $w$-summand. The last expression separates the root-reflection quotient, which will be described by alcoves, from the residual Weyl action.

Central translations also preserve this normal form: multiplying $t$ or $u$ by a central element leaves their centralizers and the involution on $J$ unchanged. Thus the same adapted tori may be used, the involution stratum is fixed, and $k=(k_1,k_2)\in B^2$ acts by $(x,[y])\mapsto(x+k_1,[y+k_2])$.

\subsection{Full centralizers}

For the graded comparison in Section~\ref{sec:graded}, we also need the action of the full centralizer on a maximal torus of its identity component.

\begin{lemma}\label{lem:centralizer} For the pair corresponding to $p=(x,[y])\in X_{B,w}$, the action of the full group $N_{C_G(t,u)}(S)$ on $\Lie S=V^+$ is the restriction of $\Stab_C(p)$ to $V^+$.
\end{lemma}
\begin{proof}
Put
\[
D=C_G(t,u),\qquad J=C_G(t)^0,\qquad C_p=\Stab_C(p).
\]
By Lemma~\ref{lem:criterion}, we have $S=(T^w)^0$ and $T=C_J(S)$. If $h\in N_D(S)$, then $h$ normalizes $J$ because it centralizes $t$, and
\[
hTh^{-1}=C_{hJh^{-1}}(hSh^{-1})=C_J(S)=T.
\]
Thus $h$ has a Weyl image $v=hT$. Projecting $huh^{-1}=u$ to $W$ gives $vwv^{-1}=w$, and Lemma~\ref{lem:linearaction} shows that $v$ fixes $p=(x,[y])$. We obtain a homomorphism
\[
\pi\colon N_D(S)\longrightarrow C_p,\qquad h\longmapsto hT.
\]

To prove surjectivity, take $v\in C_p$ and choose a lift $\dot v\in N_G(T)$. Since $v$ fixes $p$, Lemma~\ref{lem:linearaction} gives
\[
\dot v t\dot v^{-1}=t,\qquad \dot v u\dot v^{-1}=\delta u,\qquad \delta\in T_-.
\]
As $T_-=(1-w)T$, choose $a\in T$ such that $a\,w(a)^{-1}=\delta^{-1}$ and put $h=a\dot v$. Since $(\Ad u)|_T=w$, we have
\[
hth^{-1}=t,\qquad huh^{-1}=\delta\,a\,w(a)^{-1}u=u.
\]
Moreover, $\dot v$ normalizes $S$ because $v$ commutes with $w$, while $a$ centralizes $S$. Hence $h\in N_D(S)$ and $\pi(h)=v$.

The kernel consists of the elements of $T$ that centralize $u$, namely $T^w$. We thus obtain an exact sequence
\[
1\longrightarrow T^w\longrightarrow N_D(S) \xrightarrow{\ \pi\ }C_p\longrightarrow1.
\]
The kernel acts trivially on $\Lie S=V^+$, and the action of $h$ on $V^+$ is the restriction of its Weyl image $\pi(h)$. Thus $N_D(S)$ and $C_p$ have the same image on $V^+$.
\end{proof}

\section{Alcove contraction and component groups}\label{sec:contraction}

We now identify $X_{B,w}/W_0$ with the $W_0$-orbits of a finite abelian group. The construction uses the kernel of an integral torus map, which may have positive-dimensional components; after quotienting by $W_0$, an alcove contraction identifies its connected components with the adapted parameter orbits.

For each irreducible factor $\Psi_i$, let $V_i^-$ be its coroot space, let $h_i$ be its Coxeter number, and let $p_i\colon V\to V_i^-$ be the orthogonal projection. Define
\begin{equation}\label{eq:F}
H=\sum_i h_i p_i,\qquad F_w=\begin{pmatrix}m-H&-A\\-A&0\end{pmatrix}.
\end{equation}
The operator $H$ is chosen so that contraction toward $\rh_i/h_i$ removes the affine term $\rh_i$ in each root-system factor.

\subsection{Integral torus maps}

\begin{lemma}\label{lem:integrality} Let $Q_i^\vee\subset P_i^\vee$ be the coroot and coweight lattices of $\Psi_i$, and put $Q_\Psi^\vee=\bigoplus_i Q_i^\vee$. Then
\begin{equation}\label{eq:integrality}
Q_i^\vee\subset\Lambda,\qquad p_i\Lambda_B\subset P_i^\vee,\qquad h_iP_i^\vee\subset Q_i^\vee.
\end{equation}
In particular, $H\Lambda_B\subset Q_\Psi^\vee\subset\Lambda$, and $F_w$ preserves both $\Lambda^2$ and $\Lambda_B^2$.
\end{lemma}
\begin{proof}
Every coroot of $\Psi_i$ is a cocharacter of $T$. For $\kappa\in\Lambda_B$ and $\alpha\in\Psi_i$, centrality of $B$ gives
\[
\alpha(p_i\kappa)=\alpha(\kappa)\in\ZZ,
\]
which proves the second inclusion.

For the third, choose a Coxeter element $c_i\in W(\Psi_i)$. It has order $h_i$ and no fixed vectors on $V_i^-$, so its norm $1+c_i+\cdots+c_i^{h_i-1}$ is zero there. Every Weyl reflection acts trivially on $P_i^\vee/Q_i^\vee$, so applying this norm to a coweight gives $h_i\kappa\equiv0\pmod{Q_i^\vee}$. This proves the third inclusion and thus $H\Lambda_B\subset Q_\Psi^\vee\subset\Lambda$.

Since $A=1+w$ and multiplication by $m$ preserve both $\Lambda$ and $\Lambda_B$, the last assertion follows.
\end{proof}

Define the compact abelian kernel and its finite component group by
\begin{equation}\label{eq:kernel}
K_{B,w}=\ker\bigl(F_w\colon V^2/\Lambda^2\longrightarrow V^2/\Lambda_B^2\bigr), \qquad E_{B,w}=\pi_0K_{B,w}.
\end{equation}
For a point of this kernel, $F_w(x,y)$, read in $T^2$, belongs to $B^2$. We call this its \emph{linear grade}:
\begin{equation}\label{eq:lineargrade}
\gamma^{\mathrm{lin}}(x,y)=F_w(x,y)\in B^2.
\end{equation}
It is constant on connected components, since $B^2$ is finite. Moreover, $F_w$ commutes with $C$, because $C$ permutes the factors $\Psi_i$ and their projections while preserving their Coxeter numbers. Because $C$ fixes $B$ pointwise, the linear grade is also constant on $C$-orbits.

The projection $T^2\to T\times R$ has connected kernel $\{0\}\times T_-$, which lies in $K_{B,w}$ since $A$ vanishes on $T_-$. The image of $K_{B,w}$ under this projection is
\[
\bar K_{B,w}=\{(x,[y])\in T\times R:((m-H)x-Ay,-Ax)\in B^2\}.
\]
Set
\[
\mathcal K=\bar K_{B,w}/W_0.
\]
Quotienting a compact group by a connected subgroup preserves its component group, and quotienting next by a finite group identifies the orbits of components; hence
\begin{equation}\label{eq:components}
\pi_0(\mathcal K)\simeq E_{B,w}/W_0.
\end{equation}

\subsection{The full-lattice alcove quotient}

Choose a closed fundamental alcove in each $V_i^-$:
\[
\fa_i=\{z\in V_i^-:\alpha_{ij}(z)\ge0\ (1\le j\le\ell_i),\ \theta_i(z)\le1\},
\]
where $\theta_i$ is the highest root. Write
\[
\Delta=\Lambda/Q_\Psi^\vee,\qquad \mathcal A=V^+\times\prod_i\fa_i\times R.
\]
The real eigenspace decomposition need not split the cocharacter lattice, so a single lattice element can act on several coordinates simultaneously.

\begin{lemma}\label{lem:residualquotient} The group $W_0$ acts trivially on $R$, and
\begin{equation}\label{eq:affinenormal}
Q_\Psi^\vee\rtimes W_0\triangleleft\Lambda\rtimes W_0, \qquad (\Lambda\rtimes W_0)/(Q_\Psi^\vee\rtimes W_0)\simeq\Delta.
\end{equation}
There is a homeomorphism
\begin{equation}\label{eq:residualquotient}
(T\times R)/W_0\simeq\mathcal A/\Delta.
\end{equation}
For $\lambda\in\Lambda$, write $\lambda=\lambda_++\sum_i\lambda_i$ in the real decomposition, with $\lambda_i=p_i\lambda$. Let $\omega_{\lambda_i}$ be the unique alcove-preserving element of $P_i^\vee\rtimes W(\Psi_i)$ whose translation class is $\lambda_i\bmod Q_i^\vee$. Then the action on the right is
\begin{equation}\label{eq:deltaaction}
[\lambda]\cdot(v,(z_i),[y])= \bigl(v+\lambda_+,(\omega_{\lambda_i}z_i),[y]\bigr).
\end{equation}
\end{lemma}
\begin{proof}
For $\alpha\in\Psi$, the difference $s_\alpha z-z$ lies in $V^-$, so $W_0$ acts trivially on $T/T_-=R$. For a lattice vector,
\[
s_\alpha\lambda-\lambda =-\alpha(\lambda)\alpha^\vee\in Q_\Psi^\vee.
\]
Since these reflections generate $W_0$, every $v\in W_0$ satisfies $v\lambda-\lambda\in Q_\Psi^\vee$. The map
\[
\Lambda\rtimes W_0\longrightarrow\Delta,\qquad (\lambda,v)\longmapsto[\lambda]
\]
is therefore a homomorphism with kernel $Q_\Psi^\vee\rtimes W_0$, and \eqref{eq:affinenormal} follows.

We may now take the quotient of $V\times R$ in two stages. The affine Weyl group acts on each $V_i^-$ with fundamental domain $\fa_i$ and acts trivially on $V^+$ and $R$, so its quotient is $\mathcal A$. Taking the residual quotient by $\Delta$ proves \eqref{eq:residualquotient}.

Finally, $\lambda_i\in P_i^\vee$ by Lemma~\ref{lem:integrality}. In the extended affine Weyl group, the alcove stabilizer maps isomorphically to $P_i^\vee/Q_i^\vee$, so folding translation by $\lambda_i$ back into $\fa_i$ gives $\omega_{\lambda_i}$. This proves \eqref{eq:deltaaction}.
\end{proof}

In \eqref{eq:deltaaction}, the same $\lambda$ determines the translation in $V^+$ and the symmetry of every alcove factor.

\begin{example}\label{ex:latticecoupling} For $G=\SU(3)$, let $w=s_{\alpha_1}$ and use the simple-coroot basis of its cocharacter lattice. Then $\Psi=A_1$, and we may write
\[
x=s v_0+t\alpha_1^\vee,\qquad v_0=\alpha_2^\vee+\tfrac12\alpha_1^\vee\in V^+, \qquad 0\le t\le\tfrac12.
\]
The generator $[\alpha_2^\vee]$ of $\Delta=\Lambda/\ZZ\alpha_1^\vee$ acts by
\[
(s,t,[y])\longmapsto(s+1,\tfrac12-t,[y]),
\]
simultaneously translating the $V^+$ coordinate and reflecting the alcove.
\end{example}

\subsection{Image of the alcove contraction}

Write
\begin{equation}\label{eq:affinerootsum}
\theta_i=\sum_{j=1}^{\ell_i}a_{ij}\alpha_{ij}, \qquad h_i=1+\sum_{j=1}^{\ell_i}a_{ij}.
\end{equation}
Set $d_i=\rh_i/h_i$. Then
\[
\alpha_{ij}(d_i)=\frac1{h_i}\quad(1\le j\le\ell_i), \qquad \theta_i(d_i)=1-\frac1{h_i}.
\]
The function $1-\theta_i$ is the simple root corresponding to the affine vertex. Every extended alcove symmetry permutes $\alpha_{i1},\ldots,\alpha_{i\ell_i},1-\theta_i$, and since $d_i$ is the unique point where all these functions take the value $1/h_i$, every such symmetry fixes $d_i$.

For $m\ge h_i$, define
\begin{equation}\label{eq:q}
q_i(z)=\frac{\rh_i+(m-h_i)z}{m} =d_i+\frac{m-h_i}{m}(z-d_i).
\end{equation}

\begin{lemma}\label{lem:alcoveimage} Define
\[
U_{m,i}=
\left\{\xi\in\fa_i:
\alpha_{ij}(\xi)\ge 1/m\quad(1\le j\le\ell_i),\quad
\theta_i(\xi)\le 1-1/m
\right\}.
\]
If $m<h_i$, then $U_{m,i}$ is empty. If $m\ge h_i$, then
\begin{equation}\label{eq:qimage}
q_i(\fa_i)=U_{m,i}\subset\operatorname{int}(\fa_i).
\end{equation}
For $m>h_i$, this map has the unique inverse
\begin{equation}\label{eq:inverse}
z_i=\frac{m\xi_i-\rh_i}{m-h_i}.
\end{equation}
For $m=h_i$, its image is the single point $d_i$ and its fiber is all of $\fa_i$.

Every point of $X_{B,w}$, after folding into the open alcoves, has its $i$th coordinate in $U_{m,i}$. In particular, $X_{B,w}$ is empty if $m<h_i$ for some factor.
\end{lemma}
\begin{proof}
Using $\alpha_{ij}(\rh_i)=1$ and $\theta_i(\rh_i)=h_i-1$, we obtain
\begin{equation}\label{eq:qcoordinates}
\begin{aligned}
\alpha_{ij}(q_i(z)) &=\frac1m+\frac{m-h_i}{m}\alpha_{ij}(z), &&1\le j\le\ell_i,\\
\theta_i(q_i(z)) &=\frac{h_i-1}{m}+\frac{m-h_i}{m}\theta_i(z).
\end{aligned}
\end{equation}
When $m>h_i$, the map $q_i$ is invertible on $V_i^-$. These identities match the conditions $\alpha_{ij}(z)\ge0$ and $\theta_i(z)\le1$ with the defining inequalities of $U_{m,i}$. Consequently, \eqref{eq:qimage} holds and the inverse \eqref{eq:inverse} lies in $\fa_i$. The defining bounds of $U_{m,i}$ also place it in $\operatorname{int}(\fa_i)$.

For any $\xi\in U_{m,i}$, equation~\eqref{eq:affinerootsum} gives
\[
\frac{h_i-1}{m} =\sum_{j=1}^{\ell_i}\frac{a_{ij}}m \le \sum_{j=1}^{\ell_i}a_{ij}\alpha_{ij}(\xi) =\theta_i(\xi) \le\frac{m-1}{m}.
\]
Thus $U_{m,i}$ is empty when $m<h_i$. When $m=h_i$, both inequalities must be equalities. Since every $a_{ij}$ is positive, this forces $\alpha_{ij}(\xi)=1/h_i$ for every $j$, and hence $\xi=d_i$. In this case $q_i$ is constant at $d_i$, so its image is $U_{h_i,i}=\{d_i\}$ and its fiber is all of $\fa_i$.

Now take an actual parameter and fold its first coordinate into the alcoves. The adaptation criterion places each coordinate $\xi_i$ in the interior, so
\[
\alpha_{ij}(\xi_i)>0\quad(1\le j\le\ell_i), \qquad \theta_i(\xi_i)<1.
\]
Apply $\alpha_{ij}$ to the first relator equation. These roots annihilate $B$ and $Ay$, and $\alpha_{ij}(\rh)=1$, so $m\alpha_{ij}(\xi_i)$ is an integer. The highest-root expansion then shows that $m\theta_i(\xi_i)$ is also an integer. Consequently,
\begin{equation}\label{eq:inequalities}
m\alpha_{ij}(\xi_i)\in\ZZ_{\ge1} \quad(1\le j\le\ell_i), \qquad m\theta_i(\xi_i)\le m-1,
\end{equation}
so $\xi_i\in U_{m,i}$.
\end{proof}

\subsection{Linearization of the relator equations}

Assume henceforth that $m\ge h_i$ for every factor. On $\mathcal A$, apply $q_i$ in each alcove and keep $v\in V^+$ and $[y]\in R$ fixed. Each $q_i$ commutes with $\omega_{\lambda_i}$, since it is an affine contraction toward the fixed point $d_i$. The construction therefore descends through the $\Delta$-quotient to a map on the folded orbit space.

For a folded real representative $x=v+\sum_i z_i$, write $qx=v+\sum_i q_i(z_i)$. The defining formulas give
\[
mqx-\rh=(m-H)x,\qquad Aqx=Ax.
\]
Consequently, for any real lift $\widetilde y$ of $[y]$, the equality
\[
(mqx-A\widetilde y-\rh,-Aqx) =((m-H)x-A\widetilde y,-Ax)
\]
holds in $V^2$. In torus coordinates this is the linearization identity
\begin{equation}\label{eq:linearization}
\beta(qx,[y])=F_w(x,y).
\end{equation}
Thus $q$ sends the linear kernel into the actual relator locus, and \eqref{eq:qimage} ensures adaptation. We obtain
\begin{equation}\label{eq:qonspaces}
q\colon\mathcal K=\bar K_{B,w}/W_0\longrightarrow \mathcal X:=X_{B,w}/W_0,
\end{equation}
which carries the linear grade to the actual grade.

Conversely, fold any actual parameter into the open alcoves. Lemma~\ref{lem:alcoveimage} gives an inverse coordinate from \eqref{eq:inverse} in every factor with $m>h_i$, while in a factor with $m=h_i$, the actual coordinate is $d_i$ and any $z_i\in\fa_i$ is an inverse. Keeping the plus and $R$ coordinates unchanged, identity \eqref{eq:linearization} then shows that every such inverse lies in $\mathcal K$, proving surjectivity of \eqref{eq:qonspaces}.

\subsection{Sections and equivariant deformation retractions}

The inverse is unique in every noncritical factor, whereas in a critical factor the fiber is the entire alcove. We choose a section and contract these critical factors to its image.

For a folded actual coordinate $\xi_i$, set
\[
s_i(\xi_i)=
\begin{cases}
\dfrac{m\xi_i-\rh_i}{m-h_i},&m>h_i,\\
d_i,&m=h_i.
\end{cases}
\]
Keeping the other coordinates fixed gives a section
\[
s\colon\mathcal X\longrightarrow\mathcal K,\qquad q\circ s=\id_{\mathcal X}.
\]
The noncritical inverse commutes with every extended alcove symmetry, and the chosen critical point $d_i$ is fixed by every such symmetry, so $s$ descends modulo $\Delta$. It is continuous, since $\mathcal X$ is finite and discrete.

For a folded kernel point, define $H_t$ by keeping the plus and $R$ coordinates fixed and replacing $z_i$ by
\begin{equation}\label{eq:homotopy}
z_i(t)=
\begin{cases}
z_i,&m>h_i,\\
(1-t)z_i+t d_i,&m=h_i,
\end{cases}
\qquad 0\le t\le1.
\end{equation}

On a critical factor $V_i^-$, both $m-H$ and $A$ vanish. Thus, setting
\[
\delta x(t)=t\sum_{i:m=h_i}(d_i-z_i),
\]
we have
\begin{equation}\label{eq:homotopykernel}
F_w(x+\delta x(t),\widetilde y)=F_w(x,\widetilde y)\qquad\text{in }V^2.
\end{equation}
Thus the entire path satisfies the same kernel equations and has the same linear grade.

Moreover, for any extended alcove symmetry $\omega$ in a critical factor,
\[
\omega((1-t)z_i+t d_i) =(1-t)\omega(z_i)+t\omega(d_i) =(1-t)\omega(z_i)+t d_i.
\]
The homotopy therefore commutes with \eqref{eq:deltaaction} and descends to a continuous homotopy in $\mathcal K$ satisfying
\begin{equation}\label{eq:retractionidentities}
H_0=\id_{\mathcal K},\qquad H_1=sq,\qquad H_t\circ s=s\quad(0\le t\le1).
\end{equation}
It is a strong deformation retraction onto the finite section $s(\mathcal X)$, so each component of $\mathcal K$ contains a unique section point. Combining this with \eqref{eq:components} gives the following correspondence.

\begin{proposition}\label{prop:contraction} If $m<h_i$ for some $i$, then $X_{B,w}=\varnothing$. Otherwise $\bar K_{B,w}/W_0$ strongly deformation retracts onto $s(X_{B,w}/W_0)$, and
\begin{equation}\label{eq:44}
X_{B,w}/W_0\simeq\pi_0(\mathcal K)\simeq E_{B,w}/W_0.
\end{equation}
This bijection is $C/W_0$-equivariant, identifies actual grades with the linear grades on components, and intertwines central translations by $B^2$.
\end{proposition}
\begin{proof}
It remains to verify equivariance and compatibility with the grades. Represent each coset in $C/W_0$ by an element preserving $\Psi^+$. Such an element permutes the factors and their alcoves, preserves their Coxeter numbers, and sends each $\rh_i$ and $d_i$ to the corresponding vectors in the target factor, while acting linearly on the plus coordinate and on $R$. Hence $q$, $s$, and $H_t$ commute with this action and are $C/W_0$-equivariant.

For a central translation $k=(k_1,k_2)\in B^2$, choose a lift $\kappa_1\in\Lambda_B$ of $k_1$. Since $p_i\kappa_1\in P_i^\vee$, its folded action is \eqref{eq:deltaaction} with $\kappa_1$ in place of $\lambda$. The same fixed-center argument shows that $q$, $s$, and $H_t$ commute with it. Translation of the second parameter acts by $[y]\mapsto[y+k_2]$ and commutes with these maps as well.

These translations preserve the kernel locus: $H\kappa_1\in Q_\Psi^\vee\subset\Lambda$ by Lemma~\ref{lem:integrality}, and $w$ fixes $B$, so for $k=(k_1,k_2)\in B^2$, as an equality in $T^2$,
\begin{equation}\label{eq:FD}
F_wk=(mk_1-2k_2,-2k_1)=Dk\in B^2.
\end{equation}
This identifies the change of linear grade with the change of actual relator grade. Finally, \eqref{eq:linearization} identifies the grades pointwise, and finiteness of $B^2$ makes them constant on each component.
\end{proof}

\begin{figure}[htbp]
\centering
\begin{tikzpicture}[
  x=10cm,y=10cm,
  font=\small,
  zero point/.style={
    circle,fill=blue!70!black,draw=white,
    line width=.4pt,inner sep=0pt,minimum size=5.5pt
  },
  central point/.style={
    diamond,fill=orange!85!black,draw=white,
    line width=.4pt,inner sep=0pt,minimum size=7pt
  },
  boundary point/.style={
    circle,fill=black!45,
    line width=.4pt,inner sep=0pt,minimum size=4pt
  },
  q arrow/.style={
    -{Stealth[length=2mm]},dashed,black!55,
    line width=.7pt,shorten <=4pt,shorten >=4pt
  }
]


\coordinate (O)  at (0,0);
\coordinate (R)  at (1,0);
\coordinate (T)  at (.5,.5);
\coordinate (qO) at (.25,.125);
\coordinate (qR) at (.75,.125);
\coordinate (qT) at (.5,.375);
\coordinate (d)  at (.5,.25);

\fill[black!2] (O)--(R)--(T)--cycle;
\fill[blue!9] (qO)--(qR)--(qT)--cycle;

\begin{scope}
  \clip (O)--(R)--(T)--cycle;
  \draw[step=.125,black!13,very thin]
    (0,0) grid (1,.5);
\end{scope}

\draw[-{Stealth[length=2mm]},black!50]
  (0,0)--(1.08,0) node[right] {$x_1$};
\draw[-{Stealth[length=2mm]},black!50]
  (0,0)--(0,.56) node[above] {$x_2$};

\draw[black,line width=.9pt]
  (O)--(R)--(T)--cycle;
\draw[blue!65!black,line width=1pt]
  (qO)--(qR)--(qT)--cycle;

\draw[q arrow] (O)--(qO)
  node[midway,above,sloped] {$q$};
\draw[q arrow] (R)--(qR);
\draw[q arrow] (T)--(qT);

\foreach \v in {0,...,4} {
  \pgfmathtruncatemacro{\lastu}{8-\v}
  \foreach \u in {\v,...,\lastu} {
    \pgfmathtruncatemacro{\onwall}
      {(\v==0)||(\u==\v)||(\u+\v==8)}
    \ifnum\onwall=1\relax
      \node[boundary point]
        at ({\u/8},{\v/8}) {};
    \else
      \pgfmathtruncatemacro{\grade}{mod(\u+\v-3,2)}
      \ifnum\grade=0\relax
        \node[zero point]
          at ({\u/8},{\v/8}) {};
      \else
        \node[central point]
          at ({\u/8},{\v/8}) {};
      \fi
    \fi
  }
}

\path (O)--(T)
  node[midway,above,sloped] {$\alpha_1=0$};
\path (T)--(R)
  node[midway,above,sloped] {$\theta=1$};
\node[below] at (.5,-.015) {$\alpha_2=0$};

\node[below left]  at (O) {$(0,0)$};
\node[below right] at (R) {$(1,0)$};
\node[above] at (T) {$(\tfrac12,\tfrac12)$};
\node at (.79,.06) {$\mathfrak a$};
\node[text=blue!65!black]
  at (.5,.18) {$q(\mathfrak a)$};
\node[above right,inner sep=2pt] at (d) {$d$};

\node at (.5,.67)
  {$\Psi=B_2,\qquad m=8,\qquad B=\{0,b\}$};
\node at (.5,.605)
  {$q(z)=d+\tfrac12(z-d),\qquad d=(\tfrac12,\tfrac14)$};

\node[zero point] at (.02,-.105) {};
\node[anchor=west] at (.05,-.105)
  {$\beta=(0,0)$: 6 interior points};

\node[central point] at (.56,-.105) {};
\node[anchor=west] at (.59,-.105)
  {$\beta=(b,0)$: 3 interior points};

\node[boundary point] at (.02,-.165) {};
\node[anchor=west] at (.05,-.165)
  {16 boundary solutions in $B^2$; not adapted};

\node[anchor=west] at (.005,-.235)
  {$q(\mathfrak a)=
    \{x:x_1-x_2\ge\tfrac18,\ x_2\ge\tfrac18,\
         x_1+x_2\le\tfrac78\}.$};

\end{tikzpicture}
\caption{Alcove contraction in type $B_2$ for $m=8$. The six blue circles have actual grade $(0,0)$, and the three orange diamonds have grade $(b,0)$. The sixteen gray boundary points satisfy the relator grade condition but fail adaptation for the $w=-1$ stratum, so they do not belong to $X_{B,w}$. The shaded triangle is $q(\fa)$.}
\label{fig:b2-contraction}
\end{figure}
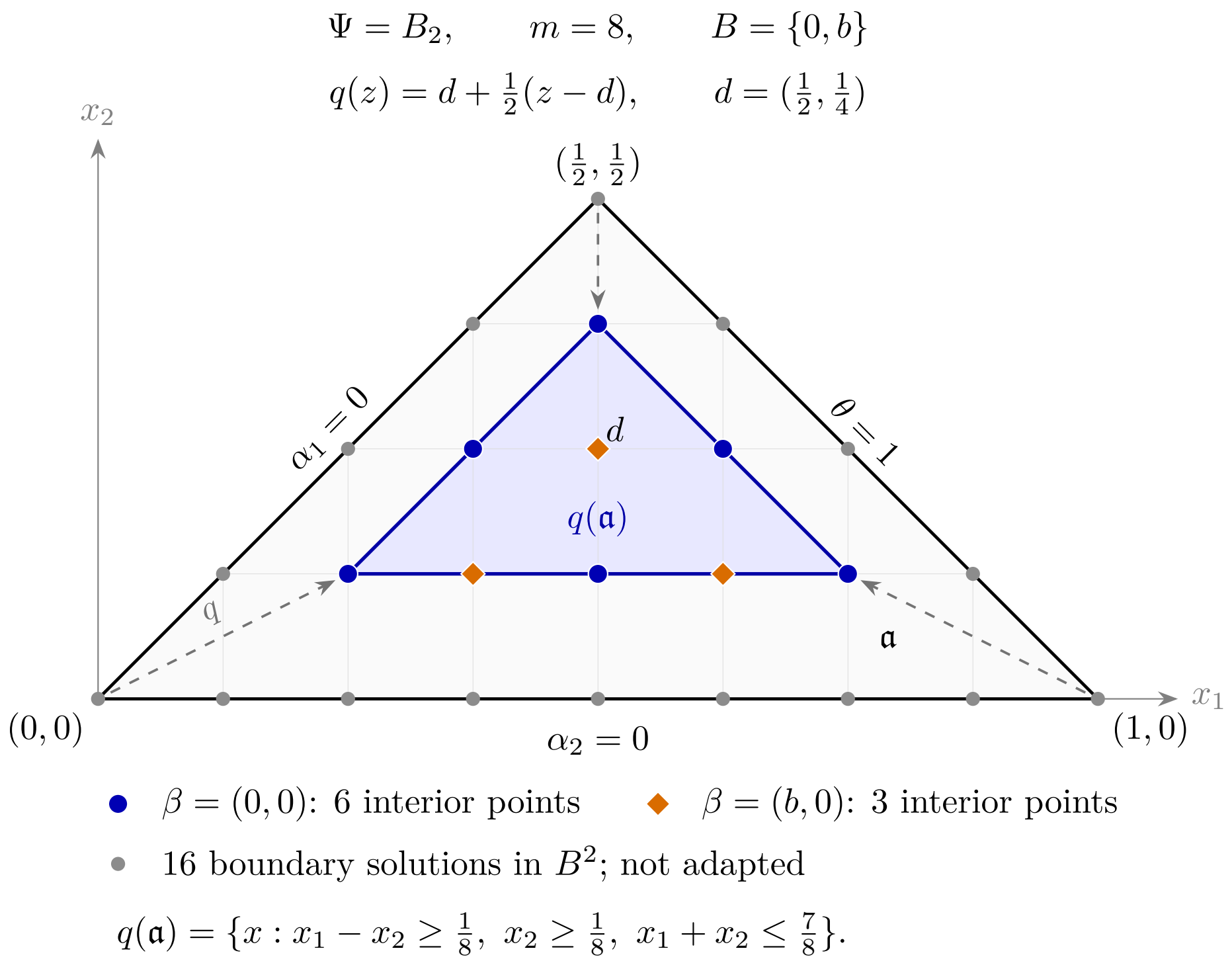

\begin{example}\label{ex:b2-contraction}
Take $G=\operatorname{Spin}(5)$, $w=-1$, $m=8$, and $B=Z(G)=\{0,b\}$ (illustrated in Figure~\ref{fig:b2-contraction}). In the coordinates $\alpha_1(x)=x_1-x_2$, $\alpha_2(x)=x_2$, one has
\[
\Lambda=\{(r,s)\in\ZZ^2:r+s\in2\ZZ\},\qquad \Lambda_B=\ZZ^2,\qquad b=[(1,0)]\in T=\RR^2/\Lambda.
\]
Here $\Psi=B_2$, $h=4$, and $\rh=(2,1)$, so the contraction center is $d=\rh/h=(1/2,1/4)$ and $q(z)=d+\tfrac12(z-d)$. Since $A=0$ and $R=0$,
\[
K_{B,w}=\{z\in T:4z\in B\}\times T,\qquad E_{B,w}=\pi_0K_{B,w}\simeq\{z\in T:4z\in B\}.
\]
Folding the first factor into $\fa$ gives the nine points of $\tfrac14\ZZ^2\cap\fa$. The identity $\beta(qz)=([4z],0)$, with brackets denoting reduction modulo $\Lambda$, carries these representatives of $E_{B,w}/W_0$ to the nine adapted points in the smaller triangle.
\end{example}

\begin{example}\label{ex:criticalcontraction} Take $G=\SU(2)$, $w=-1$, $m=2$, and $B=1$. In the coordinate $T=\RR/\ZZ$, one has $\alpha(x)=2x$, $\rh=1/2$, $h=2$, and $\fa=[0,1/2]$. Here $A=0$ and $F_w=0$, so $K_{1,w}=T^2$. Since $R=0$, its folded quotient is $\mathcal K=[0,1/2]$. The actual adapted equation is
\[
2\xi-\tfrac12\in\ZZ,\qquad 0<\xi<\tfrac12,
\]
whose unique solution is $\xi=1/4=d$. Thus $\mathcal X=\{1/4\}$, and
\[
H_t(z)=(1-t)z+t/4
\]
contracts $\mathcal K$ to its section point. The singular linear kernel contributes one component, hence one class in this stratum. Section~\ref{sec:examples} includes the other stratum and the total count.
\end{example}

\begin{remark} Fix a representative $\widetilde\xi\in\mathcal A$ of a point of $\mathcal X$. Before taking the residual quotient, its fiber is the product of the critical alcoves, since all other coordinates have unique inverses. Hence the fiber after quotienting is
\[
q^{-1}([\widetilde\xi])\simeq \left(\prod_{i:m=h_i}\fa_i\right)/\Delta_{\widetilde\xi}.
\]
Indeed, any preimage can be moved by $\Delta$ so that its image equals $\widetilde\xi$, after which the remaining identifications are given by its stabilizer $\Delta_{\widetilde\xi}$. That stabilizer is finite (its plus translation is zero and it embeds in $\prod_i P_i^\vee/Q_i^\vee$) and fixes the product of the points $d_i$, so this fiber is contractible.
\end{remark}

\begin{corollary}\label{cor:model} Call $w$ \emph{admissible} if $m\ge h_i$ for every factor of $\Psi(w)$ (vacuously when $\Psi$ is empty). Then
\begin{equation}\label{eq:classmodel}
\cR_B(G)\simeq\coprod_{[w]\ \mathrm{admissible}} E_{B,w}/C_W(w),
\end{equation}
compatibly with actual grades and central translations.
\end{corollary}
\begin{proof}
Proposition~\ref{prop:normalform} gives the actual conjugacy classes as $X_{B,w}/C$, and nonadmissible strata are empty by Proposition~\ref{prop:contraction}. In each admissible stratum, normality of $W_0$ in $C$ and the equivariance in Proposition~\ref{prop:contraction} give
\[
X_{B,w}/C \simeq (X_{B,w}/W_0)/(C/W_0) \simeq (E_{B,w}/W_0)/(C/W_0) \simeq E_{B,w}/C,
\]
compatibly with grades and central translations.
\end{proof}

We now compare the resulting finite $C_W(w)$-modules for the two dual root data.

\section{Equivariant lattice duality}\label{sec:lattice}

We first establish equivariant duality for the component groups of integral torus maps (which may be singular), and then identify the transposed map with the map associated to the dual root datum and compare the induced Weyl actions.

\subsection{Component groups of lattice maps}

\begin{lemma}\label{lem:lattice} Let $L_0\subset V_0$ and $L_1\subset V_1$ be full lattices in finite-dimensional real vector spaces, and let $F\colon V_0\to V_1$ be linear with $FL_0\subset L_1$. Put
\[
E(F)=\pi_0\ker(V_0/L_0\xrightarrow{F}V_1/L_1).
\]
Then
\begin{equation}\label{eq:torsion}
E(F)\simeq\frac{L_1\cap\im_{\RR} F}{FL_0} =\Tor\coker(F\colon L_0\to L_1).
\end{equation}
The transpose gives a torus map in the reverse direction,
\[
V_1^*/L_1^*\xrightarrow{F^t}V_0^*/L_0^*.
\]
For real lifts $x\in V_0$ and $y\in V_1^*$ of kernel points on these two sides, the formula
\begin{equation}\label{eq:pairing}
P([x],[y])=\eexp{\pair{y}{Fx}},\qquad Fx\in L_1,\quad F^ty\in L_0^*,
\end{equation}
defines a perfect pairing $E(F)\times E(F^t)\to\UU(1)$, where $[x]$ and $[y]$ denote connected components.

Suppose a finite group $C$ acts on $L_0$ and $L_1$ with $F c_0=c_1 F$, where $c_j$ denotes the action of $c\in C$ on $V_j$. Give the dual spaces the contragredient actions $c_j^\vee=(c_j^{-1})^t$. Then $F^t c_1^\vee=c_0^\vee F^t$, and $P$ is $C$-invariant and induces an isomorphism of $C$-modules
\[
E(F^t)\xrightarrow{\ \sim\ }\widehat{E(F)},\qquad [y]\longmapsto \bigl([x]\longmapsto\eexp{\pair{y}{Fx}}\bigr),
\]
where the character action is
\begin{equation}\label{eq:dualaction}
(c\chi)(e)=\chi(c^{-1}e).
\end{equation}
\end{lemma}
\begin{proof}
Send a kernel point $x\bmod L_0$ to $Fx\bmod FL_0$. This is a surjection onto the middle group in \eqref{eq:torsion}, with kernel the image of $\ker_{\RR} F$, namely the connected torus $\ker_{\RR} F/(\ker_{\RR} F\cap L_0)$. The quotient $(L_1\cap\im_{\RR} F)/FL_0$ is finite, since $\im_{\RR} F$ is a rational subspace and $FL_0$ is a full lattice in it; this gives the first equality in \eqref{eq:torsion}. A class in $L_1/FL_0$ is torsion precisely when its representative lies in the rational span of $FL_0$, namely $\im_{\RR} F$. This gives the second equality.

Changing $x$ by $\lambda\in L_0$ changes the exponent of \eqref{eq:pairing} by $\pair{F^ty}{\lambda}\in\ZZ$, and changing $y$ by $\mu\in L_1^*$ changes it by $\pair{\mu}{Fx}\in\ZZ$. Moving $x$ in $\ker_{\RR} F$ or $y$ in $\ker_{\RR} F^t$ changes the exponent by zero. Thus $P$ is well defined on components.

In Smith normal form, $F$ is diagonal with positive entries $d_1,\ldots,d_r$ and zero entries. For a nonzero entry, the kernel coordinates are $a/d_j$ and $b/d_j$ modulo $\ZZ$, and the induced pairing is
\[
(a,b)\longmapsto\eexp{ab/d_j},
\]
which is perfect, whereas the zero directions contribute only connected kernel tori or free cokernel summands and hence no component groups.

Finally, taking transposes gives $F^t c_1^\vee=c_0^\vee F^t$, and
\[
\pair{c_1^\vee y}{F c_0x} =\pair{(c_1^{-1})^t y}{c_1Fx} =\pair{y}{Fx},
\]
so $P$ is $C$-invariant. Writing $\chi_{[y]}([x])=P([x],[y])$, invariance gives
\[
\chi_{c[y]}([x])=P(c^{-1}[x],[y])=(c\chi_{[y]})([x]),
\]
which is the equivariance for \eqref{eq:dualaction}.
\end{proof}

\subsection{Langlands-dual root data}

The coefficient inclusion and its dual are
\begin{equation}\label{eq:duallattices}
\Lambda\subset\Lambda_B,\qquad \Lambda_B^*\subset\Lambda^*.
\end{equation}
The group $G/B$ has cocharacter lattice $\Lambda_B$ and character lattice $\Lambda_B^*$. Hence $G^\dagger=(G/B)^\vee$ has cocharacter lattice $\Lambda_B^*$, while the quotient $G^\dagger/B^\dagger=G^\vee$ has cocharacter lattice $\Lambda^*$. In these coordinates,
\[
B^\dagger=\Lambda^*/\Lambda_B^*\simeq\widehat{\Lambda_B/\Lambda}=\widehat B,
\]
with evaluation given by $\eexp{\pair{\eta}{\kappa}}$ for $\kappa\in\Lambda_B$ and $\eta\in\Lambda^*$.

The domain and target lattices of the kernel construction therefore exchange as follows:
\begin{center}
\begin{tabular}{lll}
\toprule
Data & $G$ with $B$ & $G^\dagger$ with $B^\dagger$\\
\midrule
Real torus space & $V$ & $V^*$\\
Cocharacter lattice & $\Lambda$ & $\Lambda_B^*$\\
Lattice after coefficient quotient & $\Lambda_B$ & $\Lambda^*$\\
Domain lattice for the pair map & $L_0=\Lambda^2$ & $L_1^*=(\Lambda_B^*)^2$\\
Target lattice for the pair map & $L_1=\Lambda_B^2$ & $L_0^*=(\Lambda^*)^2$\\
\bottomrule
\end{tabular}
\end{center}
These lattices match the setup of Lemma~\ref{lem:lattice}.

\begin{proposition}\label{prop:rootduality} Identify the Weyl groups by root--coroot exchange. The admissible involutions agree on the two sides, and for each such involution, the pairing of Lemma~\ref{lem:lattice} gives
\begin{equation}\label{eq:perfect}
E_{B^\dagger,w}(G^\dagger)\simeq\widehat{E_{B,w}(G)}
\end{equation}
as $C_W(w)$-modules, with the Weyl actions used in the finite abelian models of Corollary~\ref{cor:model}.
\end{proposition}
\begin{proof}
On the dual real space, the Weyl action corresponding to $c\in W$ is
\begin{equation}\label{eq:weylcontragredient}
c^\vee=(c^{-1})^t,
\end{equation}
which preserves both dual lattices in \eqref{eq:duallattices}. For the involution $w$ itself, $w^{-1}=w$, so $w^\vee=w^t$. The roots negated by $w$ become their coroots under duality, so $\Psi(w)$ becomes $\Psi(w)^\vee$, whose irreducible factors have the same Coxeter numbers. In particular, the inequalities defining admissibility are unchanged.

Use the dual of the chosen $W$-invariant inner product on $V^*$. The projection onto the coroot space of the corresponding dual factor is $p_i^t$. Thus $A^\vee=A^t$, $H^\vee=H^t$, and the dual kernel map is
\begin{equation}\label{eq:transpose}
F_w^\dagger=F_w^\vee =\begin{pmatrix}m-H^t&-A^t\\-A^t&0\end{pmatrix} =F_w^t.
\end{equation}
Here the underlying real linear map is the same for $G^\dagger$ and $G^\vee$, with the domain and target lattices given in the table above.

Applying Lemma~\ref{lem:lattice} with $L_0=\Lambda^2$ and $L_1=\Lambda_B^2$ gives a perfect pairing $P_w$ and the group isomorphism \eqref{eq:perfect}. To compare the Weyl actions, let $C=C_W(w)$ and let $c$ act diagonally on $V^2$. Since $F_wc=cF_w$, for real kernel lifts $\mathbf{x}\in V^2$ and $\mathbf{y}\in(V^*)^2$ we have
\begin{align}
P_w(ce,c^\vee e^\dagger) &=\eexp{\pair{c^\vee\mathbf{y}}{F_wc\mathbf{x}}}\notag\\
&=\eexp{\pair{(c^{-1})^t\mathbf{y}}{cF_w\mathbf{x}}} =P_w(e,e^\dagger).\label{eq:rootpairinginvariance}
\end{align}
Therefore $e^\dagger\mapsto P_w(-,e^\dagger)$ intertwines the dual Weyl action with the character action \eqref{eq:dualaction}. Finally, Lemma~\ref{lem:linearaction} identifies the normal-form actions on both sides with these linear Weyl actions, and Proposition~\ref{prop:contraction} carries them to the component models equivariantly.
\end{proof}

\section{Counting reciprocity}\label{sec:count}

Before constructing the Fourier transform, we recall an elementary counting consequence of Proposition~\ref{prop:rootduality}.

\begin{lemma}\label{lem:burnside} If a finite group $C$ acts by automorphisms on a finite abelian group $E$, then
\[
|E/C|=|\widehat E/C|.
\]
\end{lemma}
\begin{proof}
For $c\in C$, a fixed character factors through $E/(c-1)E$. Thus
\[
|\widehat E^{\,c}|=|E/(c-1)E|=|\ker(c-1\colon E\to E)|=|E^c|,
\]
where the middle equality follows from $|\im(c-1)|=|E|/|\ker(c-1)|$. Burnside's lemma now gives
\[
|E/C|=\frac1{|C|}\sum_{c\in C}|E^c| =\frac1{|C|}\sum_{c\in C}|\widehat E^{\,c}|=|\widehat E/C|.\qedhere
\]
\end{proof}

\begin{proof}[Proof of Theorem~\ref{thm:unified}\textup{(2)} for binary dihedral groups]
Set $B=1$. Corollary~\ref{cor:model} gives
\[
N_m(G)=\sum_{[w]\ \mathrm{admissible}}|E_{1,w}(G)/C_W(w)|.
\]
The dual root datum has the same involution classes, centralizers, and admissibility conditions. Proposition~\ref{prop:rootduality} identifies its finite group with $\widehat{E_{1,w}(G)}$, equivariantly for the $C_W(w)$-actions. Lemma~\ref{lem:burnside} equates each summand, and summing over the admissible involution classes proves the dihedral counting assertion.
\end{proof}

\section{Fourier reciprocity in degree zero}\label{sec:fourier}

The finite transforms constructed here give the degree-zero specialization of Theorem~\ref{thm:unified}\textup{(1)} and will be tensored with polynomial factors in Section~\ref{sec:graded}.

\subsection{Central operators and relator grades}

Recall that for $k=(k_1,k_2)\in B^2$, $\tau_k(t,u)=(k_1t,k_2u)$, and $D$ is defined in \eqref{eq:D}. A direct calculation with the two relators gives
\begin{equation}\label{eq:gradechange}
\gamma(\tau_k p)=\gamma(p)+Dk.
\end{equation}
The operators \eqref{eq:operators} therefore satisfy
\begin{equation}\label{eq:heisenberg}
M_\chi T_k=\chi(Dk)\,T_kM_\chi.
\end{equation}

Fix an admissible involution $w$ and write
\[
F=F_w,\qquad E=E_{B,w}(G),\qquad E^\dagger=E_{B^\dagger,w}(G^\dagger).
\]
For $e\in E$ and $e^\dagger\in E^\dagger$, choose real lifts $\mathbf x\in V^2$ and $\mathbf y\in(V^*)^2$ of kernel points in these components. By Corollary~\ref{cor:model} and Proposition~\ref{prop:rootduality}, the actual relator grades in the component models are
\[
\begin{aligned}
\gamma\colon E&\longrightarrow B^2, & \gamma(e)&=F\mathbf x\bmod\Lambda^2,\\
\gamma^\dagger\colon E^\dagger &\longrightarrow(B^\dagger)^2, & \gamma^\dagger(e^\dagger) &=F^t\mathbf y\bmod(\Lambda_B^*)^2.
\end{aligned}
\]

For $k\in B^2$ and $\chi\in\widehat{B^2}=(B^\dagger)^2$, choose representatives $\kappa\in\Lambda_B^2$ and $\eta\in(\Lambda^*)^2$. Since $F$ preserves both $\Lambda_B^2$ and $\Lambda^2$, these representatives determine kernel points on the two sides. Write $[\kappa]\in E$ and $[\eta]\in E^\dagger$ for their component classes. Under the finite abelian model,
\[
(T_kf)(e)=f(e-[\kappa]),\qquad (M_\chi f)(e)=\chi(\gamma(e))f(e).
\]
The coefficient pairing is $\chi(k)=\eexp{\pair{\eta}{\kappa}}$; in particular, $\gamma^\dagger(e^\dagger)$ is a character of $B^2$, and $\gamma^\dagger(e^\dagger)(k)$ denotes its value at $k$.

Let $P_w\colon E\times E^\dagger\to\UU(1)$ be the perfect pairing of Lemma~\ref{lem:lattice}. By definition of $P_w$,
\begin{align}
P_w([\kappa],e^\dagger) &=\eexp{\pair{F^t\mathbf y}{\kappa}} =\gamma^\dagger(e^\dagger)(k), \label{eq:ratio1}\\
P_w(e,[\eta]) &=\eexp{\pair{\eta}{F\mathbf x}} =\chi(\gamma(e)). \label{eq:ratio2}
\end{align}

\subsection{The finite Fourier transform}

Define
\begin{equation}\label{eq:transform}
(\mathcal F_w f)(e^\dagger)=\frac1{\sqrt{|E|}} \sum_{e\in E}P_w(e,e^\dagger)f(e).
\end{equation}
Since $P_w$ is perfect, $|E|=|E^\dagger|$, and character orthogonality~\cite{serre} gives the inverse
\begin{equation}\label{eq:inversefourier}
(\mathcal F_w^{-1}g)(e)=\frac1{\sqrt{|E|}} \sum_{e^\dagger\in E^\dagger}\overline{P_w(e,e^\dagger)}g(e^\dagger),
\end{equation}
which is unitary for counting measure on $E$ and $E^\dagger$.

Changing the summation variable in \eqref{eq:transform} and using \eqref{eq:ratio1} gives
\[
(\mathcal F_wT_kf)(e^\dagger) =P_w([\kappa],e^\dagger)(\mathcal F_wf)(e^\dagger) =(M_k^\dagger\mathcal F_wf)(e^\dagger),
\]
while \eqref{eq:ratio2} gives
\[
(\mathcal F_wM_\chi f)(e^\dagger) =(\mathcal F_wf)(e^\dagger+[\eta]) =(T_{-\chi}^\dagger\mathcal F_wf)(e^\dagger).
\]
The minus sign in the second identity follows from the convention $T_kf(e)=f(e-k)$ and the positive exponent in the Fourier kernel. The pairing is $C=C_W(w)$-invariant, so $\mathcal F_w$ is $C$-equivariant and restricts to an isomorphism
\[
\Fun(E,\CC)^C\simeq\Fun(E^\dagger,\CC)^C.
\]
Corollary~\ref{cor:model} identifies these invariant functions with functions on the actual class sets. Taking the direct sum over admissible $[w]$ gives the degree-zero map $\mathcal F_{B,0}$ of Theorem~\ref{thm:unified}\textup{(1)}.

\section{Graded Fourier reciprocity}\label{sec:graded}

Recall the polynomial class module \eqref{eq:gradeddef}, where each pair $p=(t,u)\in G^2$ contributes $\Sym((\Lie S)_{\CC}^*)^{W_p}$ for a maximal torus $S\subset C_G(t,u)^0$ and $W_p$ the image of $N_{C_G(t,u)}(S)$.

\begin{proposition}\label{prop:gradedmodel} There is a degree-preserving isomorphism
\begin{equation}\label{eq:gradedmodel}
\cP_B(G)\simeq\bigoplus_{[w]\ \mathrm{admissible}} \left(\Fun(E_{B,w},\CC)\otimes\Sym((V_w^+)_{\CC}^*)\right)^{C_W(w)},
\end{equation}
intertwining the operators $T_k$ and $M_\chi$ of \eqref{eq:operators} on the left with the restrictions of $T_k\otimes\id$ and $M_\chi\otimes\id$ on the right. Here $T_k$ and $M_\chi$ act on $\Fun(E_{B,w},\CC)$ as in Section~\ref{sec:fourier}, and they commute with $C_W(w)$ since $C_W(w)$ fixes $B$ pointwise.
\end{proposition}
\begin{proof}
Fix an admissible $w$ and write
\[
\begin{gathered}
X=X_{B,w},\qquad E=E_{B,w},\qquad C=C_W(w),\\
Q=C/W_0,\qquad P=\Sym((V_w^+)_{\CC}^*).
\end{gathered}
\]
Let $\cP_{B,w}(G)$ denote the contribution of this stratum to $\cP_B(G)$. For $p\in X$, put $C_p=\Stab_C(p)$. Lemma~\ref{lem:centralizer} identifies the invariant polynomial space associated to its pair with $P^{C_p}$.

For a finite $C$-set $Y$, a $C$-equivariant function $f\colon Y\to P$ is determined by its values at orbit representatives, and its value at $y$ must lie in $P^{C_y}$. Consequently,
\[
\cP_{B,w}(G) \simeq\bigoplus_{[p]\in X/C}P^{C_p} \simeq\left(\Fun(X,\CC)\otimes P\right)^C.
\]
Since $W_0$ acts trivially on $V_w^+$, it acts trivially on $P$ as well. Every $C$-equivariant $P$-valued function is therefore constant on $W_0$-orbits and descends to a $Q$-equivariant function. Proposition~\ref{prop:contraction} gives a $Q$-equivariant bijection $X/W_0\simeq E/W_0$, so
\[
\begin{aligned}
\left(\Fun(X,\CC)\otimes P\right)^C &\simeq\left(\Fun(X/W_0,\CC)\otimes P\right)^Q\\
&\simeq\left(\Fun(E/W_0,\CC)\otimes P\right)^Q\\
&\simeq\left(\Fun(E,\CC)\otimes P\right)^C.
\end{aligned}
\]
These isomorphisms preserve polynomial degree, and taking the direct sum over admissible $[w]$ gives \eqref{eq:gradedmodel}.

Central twisting leaves the pair centralizer unchanged, so the same maximal torus and polynomial variables may be used. Before taking invariants, central translation moves the parameter and leaves the polynomial value unchanged. The correspondence of Proposition~\ref{prop:contraction} intertwines central translations, so the resulting isomorphism is compatible with the translation operators. If a translation fixes a class only up to conjugacy, write $\tau_{-k}e=ce$ with $c\in C$; equivariance then gives
\[
(T_kf)(e)=f(\tau_{-k}e)=f(ce)=c\cdot f(e).
\]
Thus the finite model also reproduces the induced action on the polynomial factor when a translation fixes a class only up to conjugacy. Finally, the correspondence preserves actual grades, so it is also compatible with the grade multipliers.
\end{proof}

\begin{proof}[Proof of Theorem~\ref{thm:unified}\textup{(1)} for binary dihedral groups]
Fix an admissible $w$ and put $C=C_W(w)$, $U=V_w^+$, and $U^\dagger=(V^*)_w^+$. Restriction identifies $U^\dagger$ with $U^*$ as real $C$-representations. Averaging a positive definite inner product gives a $C$-equivariant real isomorphism $j_w\colon U\to U^\dagger$, whose inverse induces the degree-preserving polynomial isomorphism
\[
J_w=\Sym((j_w^{-1})^*)\colon \Sym(U_{\CC}^*)\xrightarrow{\ \sim\ } \Sym((U^\dagger)_{\CC}^*),\qquad (J_wp)(y)=p(j_w^{-1}y).
\]

Proposition~\ref{prop:rootduality} identifies the dual component group with $\widehat{E_{B,w}}$ as a $C$-module. Thus the finite Fourier transform \eqref{eq:transform} and $J_w$ give a degree-preserving isomorphism $(\mathcal F_w\otimes J_w)^C$ between the corresponding summands of \eqref{eq:gradedmodel} for $(G,B)$ and $(G^\dagger,B^\dagger)$. Taking the direct sum over admissible $[w]$ and using Proposition~\ref{prop:gradedmodel} on both sides defines $\mathcal F_B$.

Before taking $C$-invariants, central translations and grade multipliers act on the finite parameter factor. The operator identities established in Section~\ref{sec:fourier} therefore remain valid after tensoring with $J_w$, taking $C$-invariants, and transporting through Proposition~\ref{prop:gradedmodel}, proving \eqref{eq:fourier-rules}. In degree zero $J_w$ is the identity on constants, so the restriction of $\mathcal F_B$ is $\mathcal F_{B,0}$ constructed in Section~\ref{sec:fourier}.
\end{proof}

For $B=1$, the ordinary proof also gives the Hilbert-series equality without constructing the Fourier transform. Molien's formula reads
\begin{equation}\label{eq:molien}
\operatorname{Hilb}\cP_1(G;u)= \sum_{[w]\ \mathrm{admissible}}\frac1{|C|} \sum_{c\in C}\frac{|E_{1,w}^{\,c}|}{\det\bigl(1-u\,c|_{(V_w^+)^*}\bigr)}.
\end{equation}
The numerator is unchanged under finite abelian duality by Lemma~\ref{lem:burnside}, and the denominator is unchanged under duality of real finite-group representations. Setting $u=0$ recovers the ordinary count.

\section{Cohomological packets}\label{sec:packets}

All cohomology in this section has trivial coefficient action. Let $X$ be the presentation $2$-complex of \eqref{eq:presentation}, with one vertex, edges $a,b$, and two oriented $2$-cells $e_1,e_2$ attached along
\[
r_1=a^mb^{-2},\qquad r_2=(bab^{-1}a)^{-1}.
\]
In the integral cellular chain groups,
\[
\partial_2(e_1)=ma-2b,\qquad \partial_2(e_2)=-2a,
\]
so the matrix of $\partial_2$ in these bases is $D$. We may form $B\Gamma$ by attaching cells of dimension at least three to $X$. Since $\det D=-4$, the map $\partial_2$ is injective over $\ZZ$, so the chain identity $\partial_2\partial_3=0$ forces the integral cellular boundary $\partial_3$ to vanish.

For any abelian coefficient group $C$, the low-degree cellular cochain complex is consequently
\[
C\xrightarrow{\ 0\ } \underbrace{C^2}_{\text{coordinates on }a,b} \xrightarrow{\ D_C^t\ } \underbrace{C^2}_{\text{coordinates on }e_1,e_2} \xrightarrow{\ 0\ }C^3(B\Gamma,C).
\]
Using $D^t=D$, we obtain
\begin{equation}\label{eq:cohomology}
\begin{aligned}
H^1(\Gamma,C)&=\ker D_C^t=\ker D_C,\\
H^2(\Gamma,C)&=\coker D_C^t=\coker D_C.
\end{aligned}
\end{equation}

For $C=\ZZ$, let $e_i^*$ be the cochain taking value one on $e_i$ and zero on the other $2$-cell, and put $\rho_1=[e_1^*]$, $\rho_2=[e_2^*]$ in $H^2(\Gamma,\ZZ)$. These are the images of the two standard coordinate vectors in the cochain cokernel. They generate $H^2(\Gamma,\ZZ)$ with relations $m\rho_1=2\rho_2$ and $2\rho_1=0$, just as the classes of $a,b$ generate $H_1(\Gamma,\ZZ)=\Gamma_{\mathrm{ab}}$ with relations $ma=2b$ and $2a=0$.

Write $M=S^3/\Gamma$ for the orbit space of the left action $\gamma\cdot x=\gamma x$. The classifying map $M\to B\Gamma$ identifies cohomology in degrees one and two, since its homotopy fiber is $S^3$. Hence the same coordinates give
\[
H_1(M,\ZZ)=\coker D,\qquad H^2(M,\ZZ)=\coker D^t=\coker D,
\]
using the generator classes $a,b$ for $H_1(M,\ZZ)$ and $\rho_1,\rho_2$ for $H^2(M,\ZZ)$. These are the coordinates in which Appendix~\ref{sec:poincare} computes Poincar\'e duality.

Now let $B$ be finite. In the central-extension interpretation of $H^2(\Gamma,B)$, the value on $e_i$ is the central error obtained by evaluating its attaching word $r_i$ on chosen generator lifts. Thus an actual relator grade $b=(b_1,b_2)\in B^2$ represents the cohomology class $[b]$, and changing the lifts by $k\in B^2$ changes this vector to $b+D_Bk$. Put
\[
A_B=\ker D_B=H^1(\Gamma,B),\qquad Q_B=\coker D_B=H^2(\Gamma,B).
\]
The group $A_B$ consists of the central twists that preserve the actual grade, while $Q_B$ records the extension class after changes of lifts.

Symmetry of $D$ gives perfect coordinate evaluation pairings
\[
Q_B\times A_{B^\dagger}\longrightarrow\UU(1), \qquad A_B\times Q_{B^\dagger}\longrightarrow\UU(1),
\]
whose comparison with the spherical Poincar\'e pairing is given in Proposition~\ref{prop:poincare}.

For a grade $b\in B^2$, let $\cP_b(G)$ denote its graded summand. The subgroup $A_B$ preserves it. For $e\in\widehat{A_B}=Q_{B^\dagger}$ and $\omega=[b]\in Q_B$, define $n_G(e,\omega;u)$ to be the Hilbert series of the $e$-eigenspace in $\cP_b(G)$. It is independent of the chosen representative $b$, since translations identify all grades in the same coset of $\im D_B$ and commute with $A_B$.

Choose a representative $b_\omega\in B^2$ for each $\omega\in Q_B$, and put
\[
\cH_B(G)=\bigoplus_{\omega\in Q_B}\cP_{b_\omega}(G).
\]
The operators $T_k$ ($k\in A_B$) and $M_\chi$ ($\chi\in A_{B^\dagger}$) define commuting actions of $A_B\times A_{B^\dagger}$ on this graded space, and the joint character $(k,\chi)\mapsto\pair{k}{e}\pair{\omega}{\chi}$ occurs with Hilbert series $n_G(e,\omega;u)$.

\begin{corollary}\label{cor:packets} With the evaluation coordinates above,
\begin{equation}\label{eq:packet}
n_G(e,\omega;u)=n_{G^\dagger}(-\omega,e;u).
\end{equation}
\end{corollary}
\begin{proof}
Apply Theorem~\ref{thm:unified}\textup{(1)} in each polynomial degree. For $k\in A_B$, the identity $\mathcal F_B T_k=M_k^\dagger\mathcal F_B$ carries the electric character $e$ to the dual magnetic class $e$, while for $\chi\in\ker D_{B^\dagger}$, the identity $\mathcal F_B M_\chi=T_{-\chi}^\dagger\mathcal F_B$ carries the magnetic class $\omega$ to the dual electric character $-\omega$. This compares the joint eigenspaces summed over the entire magnetic coset. Such a coset contains $|\im D_B|$ grades, all with the same eigenspace dimensions, and $|\im D_{B^\dagger}|=|\im D_B|$ because $D_B$ and $D_{B^\dagger}$ are Pontryagin dual. Dividing by $|\im D_B|$ proves \eqref{eq:packet}.
\end{proof}

Each fixed magnetic-grade space has a real form, since its permutation action and its real invariant polynomials are defined over $\RR$. Complex conjugation therefore gives
\begin{equation}\label{eq:realform}
n_G(e,\omega;u)=n_G(-e,\omega;u).
\end{equation}

\section{Examples}\label{sec:examples}

\subsection{Tori}

If $G$ is a torus of rank $r$, the only involution is $w=1$, $\Psi$ is empty, and $F_w=D\otimes\id_V$. Since $|\det D|=4$, its ordinary kernel has $4^r$ points. Thus $N_m(G)=4^r$, as also follows from $|\Gamma_{\mathrm{ab}}|=4$, and the graded Hilbert series is $4^r/(1-u)^r$.

\subsection{The quaternion group \texorpdfstring{$Q_8$}{Q8}}

Let $m=2$, so $\Gamma=Q_8$. For $w=1$, the kernel of $D$ on $T^2$ consists of four points, all fixed by the Weyl involution, giving the four central characters. For $w=-1$, one has $\Psi=A_1$, $h=2$, $A=0$, and $F_w=0$. Its kernel is $T^2$, which has a single connected component and contributes the one noncommutative class (see Example~\ref{ex:criticalcontraction} for the deformation retraction in this stratum). Hence
\[
N_2(\SU(2))=5,\qquad \operatorname{Hilb}\cP_1(\SU(2);u)=\frac4{1-u^2}+1.
\]

\subsection{A nonsplit lattice in \texorpdfstring{$\SU(3)$}{SU(3)}}

The nontrivial involution class in $W=S_3$ is represented by a simple reflection. In the simple-coroot basis,
\[
w=\begin{pmatrix}-1&1\\0&1\end{pmatrix},\qquad H=1-w,
\]
and the integral matrix is
\[
F_w=\begin{pmatrix}
m-2&1&0&-1\\
0&m&0&-2\\
0&-1&0&0\\
0&-2&0&0
\end{pmatrix}.
\]
If $r=m-2>0$, its cokernel has generators $e_2,e_4$ and the single torsion relation $2r e_2=0$, since the other relations eliminate $e_1=-2e_2$ and $e_3=r e_2-2e_4$. Thus $E_{1,w}\simeq\ZZ/(2r)$, and $w$ acts on its generator by inversion, since $we_2=e_1+e_2=-e_2$. It contributes $r+1=m-1$ orbits. For $m=2$, the same relations give a torsion-free cokernel, so $E_{1,w}=0$ and the contribution is one. The $w=1$ stratum has five classes, so $N_m(\SU(3))=m+4$.

The same count follows from the irreducible representations of $\Dic_m$: there are four one-dimensional complex irreducibles and $m-1$ two-dimensional irreducibles. A determinant-one three-dimensional representation is either a triple of characters with product one (five multisets, for either abelianization $\ZZ/4$ or $(\ZZ/2)^2$), or a two-dimensional irreducible plus the inverse of its determinant character.

Example~\ref{ex:latticecoupling} describes the residual lattice action before taking components.

\clearpage
\part{Polyhedral groups}\label{part:polyhedral}

\section{Main results and outline}\label{sec:polyhedral-setup}

In this part we fix $\Gamma\in\{2T,2O,2I\}$. By Lemma~\ref{lem:compact-reduction}, throughout this part reductive groups are assumed to be connected compact, with complexification denoted $G_{\CC}$; we write $G^\dagger=(G/B)^\vee$ and $B^\dagger=\widehat B$.

\subsection{Main results}
Recall the spherical presentation
\[
\Gamma=\langle x,y,z\mid x^2=y^3=z^k=xyz\rangle,\qquad k\in\{3,4,5\},
\]
corresponding to $2T,2O,2I$ for $k=3,4,5$, respectively. Abelianization yields $x=3z$, $y=2z$, and $(6-k)z=0$, so that $\Gamma_{\mathrm{ab}}\simeq C_\nu$ with $\nu=6-k\in\{3,2,1\}$.

Let $X$ be the presentation $2$-complex of $\Gamma$, with one vertex, three $1$-cells corresponding to $x,y,z$, and three $2$-cells attached along the defining relations. Its cellular boundary map $C_2(X,\ZZ)\to C_1(X,\ZZ)$ is represented by
\[
\begin{pmatrix}
1 & -1 & -1\\
-1 & 2 & -1\\
-1 & -1 & k-1
\end{pmatrix},
\]
whose determinant is $k-6=-\nu\ne 0$, so $H_2(X,\ZZ)=0$. Since $B\Gamma$ is obtained from $X$ by attaching cells of dimension at least three, $H_2(\Gamma,\ZZ)=0$. Together with $H_1(\Gamma,\ZZ)\simeq\Gamma_{\mathrm{ab}}\simeq C_\nu$, the universal coefficient theorem gives, for every finite trivial $\Gamma$-module $B$,
\begin{equation}\label{eq:polyhedral-cohomology}
H^1(\Gamma,B)=B[\nu],\qquad H^2(\Gamma,B)=B/B^\nu.
\end{equation}

Fix the abelianization map
\[
q\colon\Gamma\longrightarrow C_\nu=\ZZ/\nu\ZZ,
\qquad q(z)=1,
\]
and let $\ell(g)\in\{0,\ldots,\nu-1\}$ be the standard integer representative of $q(g)$. Writing $B$ multiplicatively, define the $2$-cocycle
\[
c_b(g,h)=b^{(\ell(g)+\ell(h)-\ell(gh))/\nu}\qquad(b\in B).
\]
This is the pullback of the standard cyclic cocycle on $C_\nu$. By naturality of the universal coefficient theorem, $q^*$ is an isomorphism on $H^2$, and the second identification in \eqref{eq:polyhedral-cohomology} is explicitly
\[
B/B^\nu\xrightarrow{\sim}H^2(\Gamma,B),\qquad bB^\nu\longmapsto[c_b].
\]
Thus every $B$-valued $2$-cocycle is cohomologous to some $c_b$. We retain all actual labels $b\in B$, rather than choosing a single representative per cohomology class, to record transport between these cocycle representatives.

\begin{definition}\label{def:projective-modules}
Let $P$ be a connected compact Lie group and $B\subset Z(P)$ a finite central subgroup.
\begin{enumerate}[label=\textup{(\arabic*)},leftmargin=2em]
\item For $b\in B$, let $\cR_b(P)$ denote the finite set of $P$-conjugacy classes of maps $f\colon\Gamma\to P$ satisfying $f(g)f(h)=c_b(g,h)f(gh)$ for all $g,h\in\Gamma$.
\item A class $[f]\in\cR_b(P)$ is \emph{elliptic} if $C_P(f)^0=Z(P)^0$. Write $\cE_b(P)\subset\CC[\cR_b(P)]$ for the permutation subspace spanned by the elliptic classes at grade $b$, and write $E_\Gamma(P)$ for the set of ordinary ($b=1$) elliptic classes.
\item Define the unreduced projective class module, the unreduced elliptic module, and the graded projective class module by
\[
\begin{gathered}
\wtR_{\Gamma,B}(P)=\bigoplus_{b\in B}\CC[\cR_b(P)],\qquad
\wtE_{\Gamma,B}(P)=\bigoplus_{b\in B}\cE_b(P),\\
\wtP_{\Gamma,B}(P)=\bigoplus_{b\in B}\ \bigoplus_{[f]\in\cR_b(P)}\Sym(\mathfrak s_f^*)^{W_f},
\end{gathered}
\]
where $S_f$ is a maximal torus of $C_P(f)^0$, $\mathfrak s_f=\Lie(S_f)_{\CC}$, and $W_f=\im\bigl(N_{C_P(f)}(S_f)\to\operatorname{GL}(\mathfrak s_f)\bigr)$ is the finite group induced by the full normalizer on $\mathfrak s_f$. Linear forms in $\mathfrak s_f^*$ are assigned degree one, so that the degree-zero component of $\wtP_{\Gamma,B}(P)$ is $\wtR_{\Gamma,B}(P)$.
\end{enumerate}
\end{definition}

Let $\pi\colon P\to P/B$ be the quotient map. For fixed $b$, a map with multiplication error $c_b$ projects to an ordinary homomorphism $\Gamma\to P/B$. By Lemma~\ref{lem:compact-reduction}, these homomorphisms have finitely many conjugacy classes. For each class admitting a projective lift of grade $b$, choose a representative $\bar\rho$ and one such lift $f_0$. A conjugating element in $P/B$ lifts to $P$, so every projective class over $[\bar\rho]$ has a representative $f$ with $\pi\circ f=\bar\rho$. Write $f(g)=\lambda(g)f_0(g)$ with $\lambda(g)\in B$; since $f$ and $f_0$ have the same multiplication error, $\delta\lambda=1$, so $\lambda\in\Hom(\Gamma,B)=B[\nu]$. Each quotient class therefore contributes at most $|B[\nu]|$ projective conjugacy classes; in particular, $\cR_b(P)$ is finite.

When $\Gamma$ is clear from the context, we abbreviate $\wtE_{\Gamma,B}(P)$ to $\wtE_B(P)$. For $k\in B$, we have
\[
c_{k^\nu b}=c_b\,\delta(k^\ell),\qquad
(k^\ell)(g)=k^{\ell(g)}.
\]
Thus each $k\in B$ acts by translation between grades in the same cohomology class, and each $\chi\in B^\dagger$ acts by a grade multiplier; on $\wtR_{\Gamma,B}(P)$, $\wtE_B(P)$, and $\wtP_{\Gamma,B}(P)$, their actions are
\begin{equation}\label{eq:heisenberg-ops}
T_k(b,f)=(k^\nu b,k^\ell f),\qquad M_\chi|_b=\chi(b).
\end{equation}
They satisfy the Weyl--Heisenberg commutation relation
\[
M_\chi T_k=\chi(k^\nu)T_kM_\chi.
\]
For $k\in B[\nu]$, $T_k$ preserves each grade and is the usual $H^1(\Gamma,B)$-action by central twisting. Since $C_P(k^\ell f)=C_P(f)$, these operators preserve $\wtE_B(P)$ and act on $\wtP_{\Gamma,B}(P)$ by transport of the polynomial summands. Such transport uses a conjugator aligning the chosen representatives and maximal tori; its ambiguity lies in the full normalizer defining $W_f$, so its action on invariant polynomials is independent of this choice.

\begin{definition}\label{def:packets}
For $e\in\widehat{B[\nu]}$ and $\omega\in B/B^\nu$, the \emph{elliptic electric--magnetic packet multiplicity} is
\[
n_P^{\mathrm{ell}}(e,\omega)=\dim\{v\in\cE_b(P):T_\eta v=e(\eta)v\text{ for all }\eta\in B[\nu]\},
\]
where $b\in B$ is any lift of $\omega$. Because the transport operators $T_k$ commute with $T_\eta$ ($\eta\in B[\nu]$) and identify the grades $b$ and $k^\nu b$, this dimension is independent of the chosen representative $b$.
\end{definition}

\begin{remark}\label{rem:kt-module} To compare $\wtR_{\Gamma,B}(P)$ directly with the cohomological module introduced in~\cite[Definition~1.7]{kojima-tachikawa}, choose one representative $b_\omega\in B$ for each cohomology class $\omega\in B/B^\nu=H^2(\Gamma,B)$ and set
\[
V_B(\Gamma,P)=\bigoplus_{\omega\in B/B^\nu}\CC[\cR_{b_\omega}(P)].
\]
This module preserves a single actual grade per cohomology class, equipped with the full action of $H^1(\Gamma,B)=B[\nu]$ via $T_\eta$ on each summand and the evaluation action of $\widehat{H^2(\Gamma,B)}=\widehat{B/B^\nu}$ on the index $\omega$. Choosing transports $T_k$ gives isomorphisms between the modules arising from different representatives as representations of $H^1(\Gamma,B)\times\widehat{H^2(\Gamma,B)}$. Such an isomorphism need not be canonical, since two choices of $k$ differ by an element of $B[\nu]$, but the representation isomorphism class and its joint multiplicities are independent of these choices. Thus the degree-zero packet multiplicities of Proposition~\ref{prop:packet-reciprocity} are the joint multiplicities in $V_B(\Gamma,P)$.
\end{remark}

\begin{theorem}\label{thm:polyhedral-main}
Let $\Gamma\in\{2T,2O,2I\}$, let $G$ be a connected compact reductive group, and let $B\subset Z(G)$ be a finite central subgroup. There is a complex-linear, degree-preserving Fourier isomorphism
\begin{equation}\label{eq:final-graded-fourier}
\mathcal F\colon\wtP_{\Gamma,B}(G)\xrightarrow{\sim}\wtP_{\Gamma,B^\dagger}(G^\dagger)
\end{equation}
satisfying $\mathcal F T_k=M_k^\dagger\mathcal F$ and $\mathcal F M_\chi=T_{\chi^{-1}}^\dagger\mathcal F$ for all $k\in B$ and $\chi\in B^\dagger$. Its degree-zero restriction gives a linear isomorphism $\wtR_{\Gamma,B}(G)\xrightarrow{\sim}\wtR_{\Gamma,B^\dagger}(G^\dagger)$; in particular, setting $B=1$ yields $|\Hom(\Gamma,G)/G|=|\Hom(\Gamma,G^\vee)/G^\vee|$.
\end{theorem}

\subsection{Proof outline}

In Section~\ref{sec:central-descent}, we establish residual central descent for elliptic modules under central quotients, construct the torus Fourier transform, and decompose $\wtP_{\Gamma,B}(G)$ into relative Weyl invariants of elliptic modules on Levi subgroups tensored with $\Sym(\fa_L^*)$. Combining these ingredients, Theorem~\ref{thm:transfer} reduces Theorem~\ref{thm:polyhedral-main} to full-center elliptic Fourier isomorphisms on simply connected compact simple groups.

In Section~\ref{sec:minima}, we combine Serre's minimum theorem with an equivariant cellular Euler characteristic on $S^2$ to show that any elliptic class on a simply connected simple group must attain the minimal fixed-subalgebra dimensions $d_2,d_3,d_k$ at the three branch points, with defect $\delta_k=0$ and central image of $-1$. This restricts the problem to finitely many simple root systems.

In Section~\ref{sec:realization}, we prove Theorem~\ref{thm:local-realization}, realizing any compatible triple of local conjugacy classes with $d_1+d_2+d_3=\dim\fg$ by a unique elliptic conjugacy class via minimizing $\Gamma$-equivariant holomorphic structures on $\PP^1$ and connectedness of quasi-Hamiltonian moment-map fibers (with an alternative Yang--Mills proof of existence and ellipticity in Appendix~\ref{app:yang-mills}). Its consequence, Corollary~\ref{cor:spherical-local}, expresses the global elliptic classes directly in terms of minimal local classes.

In Section~\ref{sec:odd-models}, we compare the ordinary elliptic factors for $2T$ and $2I$ using root-lattice bijections for odd orders and Cartan-matrix duality for the central basepoint set $H_G$.

In Section~\ref{sec:minimum-fourier}, we construct a canonical unitary Fourier transform between minimal $m$th roots of the center in $G$ and $G^\dagger$, extracted from the lowest nonzero homogeneous degree of $W$-alternating polynomials on finite torus subgroups.

In Sections~\ref{sec:tetrahedral}--\ref{sec:full-reciprocity}, we assemble the full-center elliptic Fourier isomorphisms for $2T$, $2O$, and $2I$ on simply connected simple groups, combining the local models with the cyclic Fourier transform and explicit low-rank calculations for $\SU(3)$ and $\SU(4)$.

In Section~\ref{sec:polyhedral-conclusion}, we feed the simple-factor transforms into Theorem~\ref{thm:transfer} to prove Theorem~\ref{thm:polyhedral-main}, and then compare the resulting graded packet multiplicities with the spherical Poincar\'e pairing on $S^3/\Gamma$.

\section{Elliptic descent and Levi reduction}\label{sec:central-descent}

\subsection{Residual central descent and torus Fourier transform}

\begin{lemma}\label{lem:residual-descent}
Let $A\subset C\subset D\subset Z(P)$ be finite central subgroups of a connected compact Lie group $P$, and set $B=C/A\subset Z(P/A)$. There is a canonical isomorphism, equivariant under all automorphisms of $P$ preserving $A,C,D$,
\begin{equation}\label{eq:residual-descent}
\wtE_B(P/A)\simeq\wtE_D(P)^{T_A,M_{C^\perp}},
\end{equation}
where $C^\perp=\{\chi\in\widehat D:\chi|_C=1\}$. Under~\eqref{eq:residual-descent}, the residual operator $T_{cA}$ ($c\in C$) is induced by $T_c$, and for $\chi\in\widehat{C/A}$, the residual operator $M_\chi$ is induced by $M_{\widetilde\chi}$ for any extension $\widetilde\chi\in\widehat D$ of the inflation of $\chi$ to $C$. Furthermore, if $H=(P/D)^\vee$, then Fourier exchange carries the invariant subspace in~\eqref{eq:residual-descent} to the corresponding invariant subspace for
\[
(P/C)^\vee=H/C^\perp,\qquad B^\dagger\simeq A^\perp/C^\perp.
\]
\end{lemma}
\begin{proof}
Write $\pi\colon P\to P/A$ for the quotient map. This quotient identifies the Lie algebras of $Z(P)^0$ and $Z(P/A)^0$, as well as the centralizer Lie algebras of $f$ and $\pi\circ f$, so ellipticity is preserved. The condition $\chi(d)=1$ for all $\chi\in C^\perp$ singles out the grades $d\in C$, on which projection gives
\[
(d,[f])\longmapsto(dA,[\pi\circ f]),\qquad d\in C.
\]

Consider a projective map $\bar f\colon\Gamma\to P/A$ of grade $\bar b\in C/A$. Choose $c\in C$ above $\bar b$ and a pointwise lift $\widetilde f$ with $\widetilde f(1)=1$. Its multiplication error has the form $c_c\alpha$, where $\alpha$ is a normalized $A$-valued cocycle. The identification $H^2(\Gamma,A)=A/A^\nu$ gives $\alpha=c_a\delta u$ for some $a\in A$ and an $A$-valued $1$-cochain $u$. Replacing $\widetilde f$ by $f=u^{-1}\widetilde f$ cancels $\delta u$, leaving the error $c_cc_a=c_{ca}$. Thus every quotient map has a normalized lift whose actual grade belongs to $C$.

If $(d,f)$ and $(d',f')$ project to the same map, write $f'=\lambda f$ with $\lambda$ valued in $A$. The equation $\delta\lambda=c_{d'/d}$ implies $d'/d=a^\nu$ for some $a\in A$, so the cochain $\lambda/a^\ell$ is a homomorphism $\Gamma\to A$ and hence equals $\eta^\ell$ for some $\eta\in A[\nu]$; thus $(d',f')=T_{a\eta}(d,f)$. Since conjugators in $P/A$ lift to $P$, the same argument shows that the fiber over a quotient conjugacy class is one $T_A$-orbit of lift classes.

Denote the orbit over $(\bar b,[\bar f])$ by $\cO_{\bar b,[\bar f]}$. It determines an invariant vector by summing the basis vectors of its distinct members, giving the identification
\[
e_{(\bar b,[\bar f])}\longmapsto
\sum_{(d,[f])\in\cO_{\bar b,[\bar f]}}e_{(d,[f])}.
\]
These vectors form the orbit-sum basis of the required invariant space. For $a\in A$ and $\chi\in C^\perp$, the equality $\chi(a^\nu)=1$ ensures that $T_A$ commutes with $M_{C^\perp}$. Translation by $c\in C$ carries each lift orbit bijectively to the orbit over the translated quotient class. An extension $\widetilde\chi\in\widehat D$ of a character of $C/A$ has the same value on all grades in a lift orbit, so its multiplier also agrees with the quotient action. Changing $c$ by an element of $A$, or $\widetilde\chi$ by an element of $C^\perp$, acts trivially on this invariant space. Since both the projection and the orbit sums commute with every automorphism of $P$ preserving $A,C,D$, the isomorphism is equivariant.

For the dual identification, let $X_D\subset X_C\subset X_A\subset X$ denote the character lattices of $P/D$, $P/C$, $P/A$, and $P$, respectively, and write $Y_J=X_J^*$ and $Y=X^*$ for their dual cocharacter lattices. The natural evaluation pairing identifies $D=Y_D/Y$ with the Pontryagin dual of $\widehat D=X/X_D$, under which $J^\perp=X_J/X_D$ for $J\in\{A,C\}$. Quotienting $H=(P/D)^\vee$ (with character and cocharacter lattices $(Y_D,X_D)$) by $C^\perp=X_C/X_D$ produces the lattices $(Y_C,X_C)$ of $(P/C)^\vee$. Under Fourier exchange, the invariance conditions $(T_A,M_{C^\perp})$ become $(M_A,T_{C^\perp})$, which match the right-hand side of~\eqref{eq:residual-descent} for the dual chain $C^\perp\subset A^\perp\subset\widehat D$ because $(A^\perp)^\perp=A$, with residual coefficient group $A^\perp/C^\perp\simeq\widehat{C/A}=B^\dagger$.
\end{proof}

\begin{lemma}\label{lem:central-descent}
Let $P=V/Y$ be a compact torus with character lattice $X=Y^*$, and let $B=Y_B/Y\subset P$ be a finite subgroup with dual lattice $X_B=Y_B^*$. Then $\wtE_B(P)=\CC[Y_B/\nu Y]$, and the full-lattice pairing
\[
Y_B/\nu Y\times X/\nu X_B\longrightarrow\UU(1),\qquad
(y,x)\longmapsto\eexp{\pair{x}{y}/\nu},
\]
is perfect. The associated positive-kernel Fourier transform intertwines $T_k\leftrightarrow M_k^\dagger$ and $M_\chi\leftrightarrow T_{\chi^{-1}}^\dagger$ between $(P,B)$ and $(P^\dagger,B^\dagger)$, where $P^\dagger=(P/B)^\vee$, and is equivariant with respect to contragredient lattice automorphisms.
\end{lemma}
\begin{proof}
For $b\in B$, surjectivity of the power map gives $t_0\in P$ with $t_0^\nu=b$. The identity $\delta(t_0^\ell)=c_b$ shows that $t_0^{-\ell}f$ is an ordinary homomorphism; factoring it through $C_\nu$ gives $t_0^{-\ell}f=\eta^\ell$ with $\eta^\nu=1$. Hence every map of grade $b$ has the form $f=t^\ell$ with $t=t_0\eta$ and $t^\nu=b$. Conversely, each such $t$ gives a map of that grade.

For $\nu>1$, the value at an element with $\ell(g)=1$ determines $t$ uniquely, whereas for $\nu=1$, the map is $f=1$ and the label determines $t=b$. Thus the parameter space for the pairs $(b,f)$ is $(\nu^{-1}Y_B)/Y$. Writing $t=y/\nu+Y$ identifies this space with $Y_B/\nu Y$, and the actual grade is $y+Y$.

An element $y\in Y_B$ pairs trivially with every $x\in X$ if and only if $y/\nu\in Y$, and an element $x\in X$ pairs trivially with every $y\in Y_B$ if and only if $x/\nu\in X_B$, so the displayed pairing is well-defined and perfect.

Let $e_y$ be the basis vector for $y\in Y_B/\nu Y$. Choose representatives $a=\kappa+Y\in B$ and $\chi=\lambda+X_B\in B^\dagger$, with $\kappa\in Y_B$ and $\lambda\in X$. The dual parameter $x\in X/\nu X_B$ has actual grade $x+X_B$. In these coordinates the operators take the form
\[
\begin{aligned}
T_a e_y&=e_{y+\nu\kappa},
& M_\chi e_y&=\eexp{\pair{\lambda}{y}}e_y,\\
T_\chi^\dagger e_x^\dagger&=e_{x+\nu\lambda}^\dagger,
& M_a^\dagger e_x^\dagger&=\eexp{\pair{x}{\kappa}}e_x^\dagger.
\end{aligned}
\]
With $N=|Y_B/\nu Y|=|X/\nu X_B|$, define
\[
\mathcal F e_y=N^{-1/2}\sum_{x\in X/\nu X_B}\eexp{\pair{x}{y}/\nu}e_x^\dagger.
\]
The kernel $K(x,y)=\eexp{\pair{x}{y}/\nu}$ satisfies
\[
\begin{aligned}
K(x,y+\nu\kappa)&=\eexp{\pair{x}{\kappa}}K(x,y),\\
K(x+\nu\lambda,y)&=\eexp{\pair{\lambda}{y}}K(x,y).
\end{aligned}
\]
The first equality gives $\mathcal F T_a=M_a^\dagger\mathcal F$. In the second, changing the summation variable by $\nu\lambda$ moves the dual basis index by $-\nu\lambda$, giving $\mathcal F M_\chi=T_{\chi^{-1}}^\dagger\mathcal F$. The paired contragredient lattice actions preserve $K$, so the transform has the stated equivariance.
\end{proof}

\subsection{Reconstructing representations via Levi subgroups}
We now pass from elliptic modules on Levi subgroups to the full graded projective module of a connected compact reductive group $G$.

\begin{lemma}\label{lem:levi}
Let $G$ be a connected compact reductive group and $B\subset Z(G)$ a finite central subgroup. For each Levi subgroup $L\subset G$, set $A_L=Z(L)^0$, $\fa_L=\Lie(A_L)_{\CC}$, and $W_G(L)=N_G(L)/L$. Then there is a canonical graded isomorphism
\begin{equation}\label{eq:levi-graded}
\wtP_{\Gamma,B}(G)\simeq\bigoplus_{[L]}\left(\wtE_B(L)\otimes\Sym(\fa_L^*)\right)^{W_G(L)},
\end{equation}
where $[L]$ runs over $G$-conjugacy classes of Levi subgroups of $G$. Taking degree zero gives $\wtR_{\Gamma,B}(G)\simeq\bigoplus_{[L]}\wtE_B(L)^{W_G(L)}$, and specializing to $B=1$ yields
\begin{equation}\label{eq:levi-classes}
\Hom(\Gamma,G)/G\simeq\coprod_{[L]}E_\Gamma(L)/W_G(L).
\end{equation}
\end{lemma}
\begin{proof}
For a projective map $f\colon\Gamma\to G$, put $H=C_G(f)^0$, choose a maximal torus $A\subset H$, and set $L=C_G(A)$. The centralizer of a torus is connected and is a Levi subgroup, and a maximal torus in a connected compact group is its own centralizer~\cite[15.14 and~23.13]{Milne}. Since $f(\Gamma)\subset L$, we have
\[
A\subseteq Z(L)^0\subseteq H.
\]
The maximality of $A$ gives $Z(L)^0=A$, while self-centralization in $H$ gives
\[
C_L(f)^0=C_H(A)=A.
\]
Thus $f$ is elliptic in $L$.

Conversely, let $f$ be elliptic in a Levi subgroup $L$ and put $A_L=Z(L)^0$. If a torus $R\subset C_G(f)^0$ contains $A_L$, then $R\subset C_G(A_L)=L$. Hence $R\subset C_L(f)^0=A_L$, so $A_L$ is a maximal torus of the ambient centralizer. An ambient conjugator between two elliptic maps can therefore be corrected inside the connected centralizer to align their chosen maximal tori, and hence conjugates the corresponding Levis as well. After fixing a representative $L$, ambient conjugacy is consequently described by the action of $W_G(L)$ on its elliptic classes.

Fix an elliptic $L$-class $[f]$. The map from $N_{C_G(f)}(A_L)$ to $W_G(L)$ has kernel $C_L(f)$ and image the stabilizer of $[f]$. Indeed, if $n\in N_G(L)$ carries $f$ to an $L$-conjugate, composing $n$ with an appropriate $L$-conjugator gives an element centralizing $f$ with the same action on $A_L$. Since $C_G(A_L)=L$, the kernel of the action on $\fa_L$ is also $C_L(f)$. Thus, as groups acting on $\fa_L$,
\[
W_f\simeq N_{C_G(f)}(A_L)/C_L(f)
\simeq\Stab_{W_G(L)}([f]).
\]

For a fixed actual grade, let $S$ be the set of elliptic $L$-classes, put $R=\Sym(\fa_L^*)$, and write $W=W_G(L)$. An invariant family $(p_s)_{s\in S}$ satisfies $p_{ws}=wp_s$. Its values on each orbit are therefore determined by a polynomial fixed by the stabilizer of a representative, and hence
\[
\left(\bigoplus_{s\in S}R\right)^W
\simeq\bigoplus_{[s]\in S/W}R^{W_s}.
\]
The stabilizer identification shows that these are the summands $\Sym(\mathfrak s_f^*)^{W_f}$ attached to the ambient classes. Summing over actual grades and Levi classes proves~\eqref{eq:levi-graded} in every degree.

Central twists preserve $C_G(f)$ and the associated Levi, and an $L$-conjugator acts trivially on $A_L$. They also commute with the relative Weyl action, because $B\subset Z(G)$ is fixed pointwise by $W_G(L)$. If a twist carries an elliptic class to another representative of its Weyl orbit, the relation $p_{ws}=wp_s$ gives the corresponding normalizer action on the polynomial. Grade multiplication acts by the same scalar on both descriptions. Thus~\eqref{eq:levi-graded} intertwines all $T_k$ and $M_\chi$.
\end{proof}

\begin{lemma}\label{lem:zero-packet}
Suppose that the elliptic packet multiplicities $n_P^{\mathrm{ell}}(e,\omega)$ of $(P,B)$ are supported only at $(e,\omega)=(0,0)$. Equivalently, elliptic classes occur only in actual grades $b\in B^\nu$, and twisting by $B[\nu]$ fixes every ordinary elliptic class. Then there is a canonical, automorphism-equivariant identification
\begin{equation}\label{eq:zero-packet}
\wtE_B(P)\simeq\CC[B^\nu]\otimes\CC[E_\Gamma(P)].
\end{equation}
The pairing
\begin{equation}\label{eq:zero-pairing}
B^\nu\times(B^\dagger)^\nu\longrightarrow\UU(1),\qquad
\pair{k^\nu}{\chi^\nu}=\chi(k)^\nu,
\end{equation}
is perfect. Its positive-kernel Fourier transform exchanges $T_k$ with $M_k^\dagger$ and $M_\chi$ with $T_{\chi^{-1}}^\dagger$ on the coefficient factor.

Let $P^\dagger=(P/B)^\vee$ and suppose that $(P^\dagger,B^\dagger)$ also satisfies the zero-packet condition. A given automorphism-equivariant isomorphism $\mathcal F_0\colon\CC[E_\Gamma(P)]\xrightarrow{\sim}\CC[E_\Gamma(P^\dagger)]$ then determines a full coefficient Fourier isomorphism
\[
\wtE_B(P)\xrightarrow{\sim}\wtE_{B^\dagger}(P^\dagger),
\]
obtained by tensoring $\mathcal F_0$ with the coefficient Fourier transform under the two zero-packet decompositions.
\end{lemma}
\begin{proof}
Zero packet support gives the asserted restriction on actual grades. At grade $1$, it also says that the finite abelian group $B[\nu]$ acts with only the trivial character, so it fixes the ordinary elliptic class basis. For $b=k^\nu$, transport by $T_k$ identifies the ordinary space with the grade-$b$ space. Two choices of $k$ differ by $B[\nu]$, and hence give the same identification, proving \eqref{eq:zero-packet} and its equivariance for automorphisms. Conversely, the two stated conditions imply zero packet support, since all allowed layers are obtained by transports commuting with $B[\nu]$.

The value in~\eqref{eq:zero-pairing} is unchanged when $k$ is multiplied by $B[\nu]$ or $\chi$ by $B^\dagger[\nu]$. If $k^\nu$ pairs trivially with every $\chi^\nu$, then $\chi(k^\nu)=1$ for every $\chi\in B^\dagger$, so $k^\nu=1$; since $|B^\nu|=|(B^\dagger)^\nu|$, the pairing is perfect. Translating by $k^\nu$ multiplies the positive Fourier kernel by evaluation at $k$, while multiplication by $\chi$ shifts the dual coefficient by $\chi^{-\nu}$. The claimed operator identities follow.

Applying~\eqref{eq:zero-packet} to both $(P,B)$ and $(P^\dagger,B^\dagger)$, where the central operators act on the coefficient factors, the tensor product of their Fourier transform with $\mathcal F_0$ gives the final assertion, equivariantly under the corresponding automorphisms.
\end{proof}

\begin{remark}\label{rem:coprime-packet}
Whenever $\gcd(\nu,|B|)=1$, both $B[\nu]$ and $B/B^\nu$ vanish, so Lemma~\ref{lem:zero-packet} applies automatically; in this case $B^\nu=B$, and the pairing~\eqref{eq:zero-pairing} is simply $(b,\chi)\mapsto\chi(b^{1/\nu})$.
\end{remark}

For a simply connected compact simple group $S$, we denote by
\[
S^\dagger=(S/Z(S))^\vee
\]
the simply connected compact simple group of the dual root system, so that $\widehat{Z(S)}\simeq Z(S^\dagger)$.

\begin{theorem}\label{thm:transfer}
Let $\Gamma\in\{2T,2O,2I\}$. Suppose that for every simply connected compact simple group $S$ there exists a full-center elliptic Fourier isomorphism
\[
\wtE_{Z(S)}(S)\xrightarrow{\sim}\wtE_{Z(S^\dagger)}(S^\dagger)
\]
intertwining $T_k\leftrightarrow M_k^\dagger$ and $M_\chi\leftrightarrow T_{\chi^{-1}}^\dagger$, equivariant under based root-datum automorphisms and compatible with permutations of isomorphic simple factors. Then, for every connected compact reductive group $G$ and every finite central subgroup $B\subset Z(G)$, there is a degree-preserving Fourier isomorphism
\[
\wtP_{\Gamma,B}(G)\xrightarrow{\sim}\wtP_{\Gamma,B^\dagger}(G^\dagger)
\]
intertwining $T_k\leftrightarrow M_k^\dagger$ and $M_\chi\leftrightarrow T_{\chi^{-1}}^\dagger$.
\end{theorem}
\begin{proof}
For a Levi subgroup $L\subset G$, use the product cover
\[
\pi_L\colon P=T\times\prod_i S_i\longrightarrow L,
\qquad T=Z(L)^0,
\]
where the $S_i$ are simply connected compact simple and $A=\ker\pi_L$ is finite and central. The preimage $C=\pi_L^{-1}(B)$ is finite and central as well; choose $n\ge1$ annihilating its projection to $T$ and put
\[
D=T[n]\times\prod_i Z(S_i),
\]
so that $A\subset C\subset D$.

Conjugation by $N_G(L)$ preserves $T$ and the derived group of $L$, and its action on the derived group lifts uniquely to the simply connected cover~\cite[16.10 and~16.21]{Milne}. Uniqueness makes these lifts compatible with composition, so they define an action on $P$ commuting with $\pi_L$. This action preserves $A$ and $C$, since $B$ is fixed pointwise, and it preserves $D$ because $T[n]$ and the simple centers are characteristic and isomorphic factors may be permuted. For $l\in L$, conjugation by any lift of $l$ to $P$ realizes the induced automorphism of $P$ and therefore acts trivially on conjugacy class spaces, so the action factors through $W_G(L)$.

Ellipticity in $P$ is factorwise. Tensor the torus Fourier transform of Lemma~\ref{lem:central-descent} with the assumed full-center transforms on the $S_i$. The torus transform respects contragredient lattice actions, while on the simple factors, the lifted action is a based root-datum automorphism up to inner conjugacy and permutation of isomorphic factors, so the hypotheses give equivariance of the tensor product. Taking $(T_A,M_{C^\perp})$-invariants by Lemma~\ref{lem:residual-descent} therefore yields a $W_G(L)$-equivariant elliptic Fourier isomorphism
\[
\mathcal F_L\colon\wtE_B(L)\xrightarrow{\sim}\wtE_{B^\dagger}(L^\dagger),
\qquad L^\dagger=(L/B)^\vee.
\]

Choose a maximal torus $T_0\subset L$ and a positive system making $L$ standard, with root subsystem $\Phi_I$. Write $(X,\Phi,Y,\Phi^\vee)$ for the root datum of $G$, and let $X_B\subset X$ be the character lattice of $G/B$, with $Y_B=X_B^*$. The root datum of $L^\dagger$ is
\[
(Y_B,\Phi_I^\vee,X_B,\Phi_I).
\]
Thus $L^\dagger$ is the corresponding Levi in $G^\dagger=(G/B)^\vee$. Aligning maximal tori inside $L$ identifies its relative Weyl group with the normalizer of the root subsystem modulo its Weyl group, and root duality gives the same quotient on the other side, namely
\[
W_G(L)\simeq N_W(\Phi_I)/W_I\simeq W_{G^\dagger}(L^\dagger),
\qquad W=W(G,T_0),\quad W_I=W(\Phi_I).
\]
The two actions are induced by the paired lattice actions of $N_W(\Phi_I)$ on the two root data, as in the construction of $\mathcal F_L$. The same root correspondence also matches the conjugacy classes of Levis on the two sides.

Put $U=\Lie(Z(L)^0)$ and $U^\dagger=\Lie(Z(L^\dagger)^0)$, as real vector spaces. In the real cocharacter spaces, they are the respective annihilators
\[
\begin{aligned}
U&=\{y\in Y\otimes_{\ZZ}\RR:\pair{\alpha}{y}=0\text{ for all }\alpha\in\Phi_I\},\\
U^\dagger&=\{x\in X\otimes_{\ZZ}\RR:\pair{x}{\alpha^\vee}=0\text{ for all }\alpha^\vee\in\Phi_I^\vee\}.
\end{aligned}
\]
The root and coroot spans pair nondegenerately, so the lattice pairing restricts to a nondegenerate real pairing between $U$ and $U^\dagger$. It is invariant under the common relative Weyl group, since the normalizer preserves the two root spans and $W_I$ acts trivially on both centers. Choose a positive definite inner product on $U$ and average it over this finite group; this identifies $U^*$ with $U$, while the pairing identifies $U$ with $(U^\dagger)^*$. Complexification and passage to symmetric algebras then give a degree-preserving, equivariant complex-linear isomorphism
\[
J_L\colon\Sym(\fa_L^*)\xrightarrow{\sim}\Sym(\fa_{L^\dagger}^*).
\]
 Tensoring $J_L$ with $\mathcal F_L$ and taking invariants under the common relative Weyl group gives the corresponding summand of Lemma~\ref{lem:levi} on each side. Summing over the matched Levi classes proves the theorem, with all polynomial degrees and central operator identities preserved.
\end{proof}

A related use of elliptic data on Levi subgroups appears in the nonabelian Fourier transform for unipotent representations of reductive $p$-adic groups: Aubert--Ciubotaru--Romano reconstruct the compact representation space by parabolic induction and, under compatibility hypotheses, reduce Fourier compatibility to the essentially elliptic spaces of Levi subgroups~\cite[Theorems~5.2 and~5.16]{ACR}. Here the transfer is obtained from residual central descent and the full polynomial reconstruction~\eqref{eq:levi-graded}.

\section{Minimality and centrality}\label{sec:minima}

The binary polyhedral groups $\Gamma\in\{2T,2O,2I\}$ have the presentation
\[
\Gamma=\langle x,y,z\mid x^2=y^3=z^k=xyz\rangle,
\qquad k=3,4,5,
\]
respectively.

\subsection{Minimality}

We first recall the minimum dimension of the fixed subalgebra of a periodic inner automorphism, and then apply it to the three local generators of a spherical group.

\begin{lemma}\label{lem:minima}
Let $L$ be the simply connected compact group with compact semisimple Lie algebra $\fl$. Choose a maximal torus and a positive root system $\Phi^+$, and put
\[
\rh=\frac12\sum_{\alpha\in\Phi^+}\alpha^\vee,\qquad
s_m=(2\rh)(e^{\pi i/m}),\qquad
d_m=\dim_{\RR}\fl^{\Ad(s_m)}
\quad(m\ge1).
\]
For an inner automorphism $\tau$ of $\fl$, write $\fl^\tau$ for its fixed subalgebra.
\begin{enumerate}[label=\textup{(\arabic*)},leftmargin=2em]
\item Every inner automorphism $\tau$ satisfying $\tau^m=1$ obeys
\[
\dim_{\RR}\fl^\tau\ge d_m.
\]
\item Let $\bar\Gamma=\langle \bar x,\bar y,\bar z\mid \bar x^2=\bar y^3=\bar z^k=\bar x\bar y\bar z=1\rangle$ be the spherical rotation group of type $(2,3,k)$, where $k\in\{3,4,5\}$. For a homomorphism $\sigma\colon\bar\Gamma\to\operatorname{Inn}(\fl)$, put $\sigma_{\bar x}=\sigma(\bar x)$, $\sigma_{\bar y}=\sigma(\bar y)$, and $\sigma_{\bar z}=\sigma(\bar z)$. If $\fl^{\sigma(\bar\Gamma)}=0$, then
\[
\dim_{\RR}\fl^{\sigma_{\bar x}}=d_2,\qquad
\dim_{\RR}\fl^{\sigma_{\bar y}}=d_3,\qquad
\dim_{\RR}\fl^{\sigma_{\bar z}}=d_k,\qquad
\dim_{\RR}\fl=d_2+d_3+d_k.
\]
\end{enumerate}
\end{lemma}

\begin{proof}
Part~(1) is Serre's minimum theorem applied to the complexification $\fl_{\CC}$~\cite[Theorem~1 and~\S3]{EKV}, where the hypothesis $\tau^m=1$ allows the order of $\tau$ to divide $m$. The assertion is immediate for $\fl=0$; in the remaining construction assume $\fl\ne0$.

For part~(2), realize $\bar\Gamma$ as the rotation group of the regular tetrahedron, octahedron, or icosahedron, according as $k=3,4,5$. The barycentric subdivision of its boundary, projected radially onto $S^2$, is a $\bar\Gamma$-equivariant triangulation. It has three vertex orbits, with stabilizers conjugate to $\langle \bar x\rangle$, $\langle \bar y\rangle$, and $\langle \bar z\rangle$, three edge orbits, and two face orbits, where edge and face stabilizers are trivial. Orienting edges by their endpoint types and faces by the orientation of $S^2$ makes the cellular actions permutation actions on the oriented cells.

Let $V$ be any finite-dimensional complex $\bar\Gamma$-module, and write $C_\bullet=C_\bullet(S^2,\CC)$ for this cellular chain complex. The invariant chain complex has the form
\[
(V\otimes_{\CC}C_\bullet)^{\bar\Gamma}:
\quad
0\longrightarrow V^{\oplus2}
\longrightarrow V^{\oplus3}
\longrightarrow V^{\bar x}\oplus V^{\bar y}\oplus V^{\bar z}
\longrightarrow0,
\]
with the three nonzero terms in degrees $2,1,0$. Averaging over $\bar\Gamma$ makes the invariant functor exact, and hence
\[
H_q\bigl((V\otimes_{\CC}C_\bullet)^{\bar\Gamma}\bigr)
\simeq
\bigl(V\otimes_{\CC}H_q(S^2,\CC)\bigr)^{\bar\Gamma}.
\]
The action on $H_0(S^2,\CC)$ and $H_2(S^2,\CC)$ is trivial, since the sphere is connected and rotations preserve orientation. Thus the zeroth and second homology of the invariant complex are both $V^{\bar\Gamma}$, while its first homology vanishes. Comparing Euler characteristics gives
\begin{equation}\label{eq:spherical-euler}
\dim_{\CC}V^{\bar x}+\dim_{\CC}V^{\bar y}+\dim_{\CC}V^{\bar z}
=\dim_{\CC}V+2\dim_{\CC}V^{\bar\Gamma}.
\end{equation}

For any group $K$ acting on the real Lie algebra $\fl$, fixed spaces commute with complexification:
\[
(\fl_{\CC})^K=\fl^K\otimes_{\RR}\CC,\qquad
\dim_{\CC}(\fl_{\CC})^K=\dim_{\RR}\fl^K.
\]
We use real dimensions for Lie subalgebras below. Apply~\eqref{eq:spherical-euler} to $\fl_{\CC}$ with the given action $\sigma$. This action may have a kernel, but its local images satisfy $\sigma_{\bar x}^2=\sigma_{\bar y}^3=\sigma_{\bar z}^k=1$. The hypothesis on common invariants and part~(1) therefore give
\[
\dim\fl
=\dim\fl^{\sigma_{\bar x}}+\dim\fl^{\sigma_{\bar y}}+\dim\fl^{\sigma_{\bar z}}
\ge d_2+d_3+d_k.
\]

To obtain the reverse inequality, we construct a reference action using the principal $\SU(2)$ of $L$. Choose simple coroots $H_i$ and root vectors $E_i,F_i$ with $[E_i,F_i]=H_i$, compatible with the compact form so that compact conjugation sends $E_i$ to $-F_i$. Write
\[
H_0=2\rh=\sum_i a_iH_i,\qquad
e=\sum_i\sqrt{a_i}E_i,\qquad
f=\sum_i\sqrt{a_i}F_i,
\qquad a_i>0.
\]
These elements satisfy $[H_0,e]=2e$, $[H_0,f]=-2f$, and $[e,f]=H_0$, giving the principal $\mathfrak{sl}_2$-triple~\cite[\S32.3]{EtingofExponents}. Its compact real form is $\operatorname{Span}_{\RR}\{iH_0,e-f,i(e+f)\}\subset\fl$. The resulting Lie algebra homomorphism integrates to
\[
\iota\colon\SU(2)\longrightarrow L,\qquad
\iota\bigl(\operatorname{diag}(t,t^{-1})\bigr)
=(2\rh)(t)\quad(t\in\UU(1)).
\]
For every root $\alpha$, the integer $\pair{\alpha}{2\rh}=2\rootht(\alpha)$ is even. Thus $\iota(-I)=(2\rh)(-1)$ is central in $L$, and the adjoint action factors through
\[
\varphi\colon\SO(3)\longrightarrow\operatorname{Inn}(\fl).
\]

Restrict $\varphi$ to the geometric subgroup $\bar\Gamma\subset\SO(3)$ and denote the resulting action by $\sigma_0$, with local images $\sigma_{0,\bar x},\sigma_{0,\bar y},\sigma_{0,\bar z}$. An order-$m$ rotation is conjugate to rotation through $2\pi a/m$ with $\gcd(a,m)=1$, so its image is conjugate to $(\Ad s_m)^a$. Since $(\Ad s_m)^m=1$, this power has the same fixed subalgebra as $\Ad s_m$. Consequently,
\[
\dim\fl^{\sigma_{0,\bar x}}=d_2,\qquad
\dim\fl^{\sigma_{0,\bar y}}=d_3,\qquad
\dim\fl^{\sigma_{0,\bar z}}=d_k.
\]
Applying~\eqref{eq:spherical-euler} to $\fl_{\CC}$ with this reference action yields
\[
d_2+d_3+d_k
=\dim\fl+2\dim\fl^{\sigma_0(\bar\Gamma)}
\ge\dim\fl.
\]
Combining the two comparisons gives
\[
\dim\fl\ge d_2+d_3+d_k\ge\dim\fl,
\]
so all inequalities are equalities; in particular,
\[
(\dim\fl^{\sigma_{\bar x}}-d_2)
+(\dim\fl^{\sigma_{\bar y}}-d_3)
+(\dim\fl^{\sigma_{\bar z}}-d_k)=0.
\]
Each summand is nonnegative by part~(1), so each vanishes, proving part~(2); moreover, the equality for the reference action gives $\fl^{\sigma_0(\bar\Gamma)}=0$.

The reference action depends only on $\fl$ and the geometric spherical group, so for every compact semisimple $\fl$, whether or not an elliptic action exists, the same construction gives
\[
\delta_k(\fl):=d_2+d_3+d_k-\dim\fl
=2\dim\fl^{\sigma_0(\bar\Gamma)}\ge0.
\]
\end{proof}

For a connected compact semisimple group $G$ with Lie algebra $\fg$, we also write $d_{m,\fg}$ or $d_m(G)$ for the minimum in part~(1), and $d_m$ when the Lie algebra is clear. An element $X\in G$ satisfying $(\Ad X)^m=1$ is called minimal if $\dim_{\RR}\fg^{\Ad(X)}=d_{m,\fg}$.

\subsection{Centrality}

\begin{lemma}\label{lem:involution}
Let $G$ be a connected compact semisimple group with Lie algebra $\fg$, and set $u=(2\rh)(i)$.
\begin{enumerate}[label=\textup{(\arabic*)},leftmargin=2em]
\item The inner involutions of $\fg$ whose fixed-point subalgebras have dimension $d_2$ form a single conjugacy class in $\operatorname{Inn}(\fg)$, represented by $\Ad(u)$.
\item If $U\in G$ induces such an involution, then for some $b\in Z(G)$,
\[
U^2=b^2z_G,\qquad z_G=u^2=(2\rh)(-1).
\]
\item Any such $U$ is $G$-conjugate to $bu$ for some $b\in Z(G)$.
\end{enumerate}
\end{lemma}

\begin{proof}
Assertions \textup{(1)} and \textup{(2)} are classical in~\cite[\S3.2]{SLinv}. For \textup{(3)}, part~\textup{(1)} gives $g\in G$ such that $\Ad(gUg^{-1})=\Ad(u)$. Thus $gUg^{-1}u^{-1}$ belongs to $\ker(\Ad)=Z(G)$, so $gUg^{-1}=bu$ for some $b\in Z(G)$, and hence $U=g^{-1}(bu)g$.
\end{proof}

\begin{corollary}\label{cor:squares}
Suppose $G$ is a simply connected compact semisimple group, and set $B=Z(G)$, $u=(2\rh)(i)$, and $z_G=u^2$. The minimal quadratic classes with central square form the set $\{[bu]_G:b\in B\}$. Every class in this set has stabilizer $B[2]$ under central translation, and squaring induces a canonical bijection
\[
\{[bu]_G:b\in B\}\xrightarrow{\sim}z_GB^2,
\qquad [U]_G\longmapsto U^2.
\]
\end{corollary}

\begin{proof}
Lemma~\ref{lem:involution} identifies these classes with the central translates of $[u]_G$. Their stabilizers coincide because $B$ is abelian and acts transitively. If $b$ fixes $[u]_G$, conjugacy of $bu$ and $u$ forces equality of their central squares, hence $b^2z_G=z_G$ and $b\in B[2]$.

It remains to show that every $b\in B[2]$ fixes $[u]_G$. A simply connected compact semisimple group is the product of its simply connected simple factors; both $b$ and $u$ decompose accordingly, so we work in one factor. Use the maximal torus $T$ and simple roots defining $u$. Simple coroots form a basis of $X_*(T)=Q^\vee$, and evaluation at $-1$ identifies $Q^\vee/2Q^\vee$ with $T[2]$. Consequently, there is a unique subset $I$ of the simple nodes such that
\[
b=\prod_{i\in I}\alpha_i^\vee(-1).
\]
Centrality of $b$ means that each simple root takes the value $1$ on $b$, or equivalently
\[
\sum_{i\in I}\pair{\alpha_j}{\alpha_i^\vee}
\equiv0\pmod2
\qquad\text{for every simple node }j.
\]

Consider the subgraph on the nodes $I$, and sum these congruences over $j\in I$. A diagonal entry contributes $2$, and the two Cartan entries belonging to a single, double, or triple edge sum to $-2$, $-3$, or $-4$, respectively. Only a double edge contributes an odd number. Each connected finite Dynkin diagram is a tree and has at most one double edge~\cite[\S23.8]{EtingofExponents}, so the subgraph on $I$ has no double edge. All its remaining edge entries are odd, and the individual congruences now say that each of its vertices has even degree. Since a nontrivial finite tree has a leaf, this forest has no edges, so the nodes in $I$ are pairwise nonadjacent.

For each $i\in I$, the equality $\pair{\alpha_i}{2\rh}=2$ gives $\alpha_i(u)=-1$. The simple reflection $s_i$ therefore satisfies
\[
s_i(u)=u\,\alpha_i^\vee\bigl(\alpha_i(u)^{-1}\bigr)
=\alpha_i^\vee(-1)u.
\]
The selected nodes are nonadjacent, so their simple reflections commute and each fixes the other selected coroots. It follows that
\[
\left(\prod_{i\in I}s_i\right)(u)
=\left(\prod_{i\in I}\alpha_i^\vee(-1)\right)u=bu.
\]
Thus $bu$ is conjugate to $u$. Applying this construction in each simple factor proves that the common stabilizer is $B[2]$. The class set is therefore $B/B[2]$, and its square map is
\[
bB[2]\longmapsto z_Gb^2.
\]
Its image is $z_GB^2$, and two classes have the same square if and only if their representatives differ by an element of $B[2]$. The map is therefore bijective.
\end{proof}

\begin{proposition}\label{prop:central-involution}
Let $G$ be a simply connected compact simple group, and let $k\in\{3,4,5\}$. Suppose that $X,Y,Z,h\in G$ satisfy
\[
h^2=1,\qquad [h,X]=[h,Y]=[h,Z]=1,\qquad X^2\in hZ(G),\qquad Y^3,Z^k,XYZ\in Z(G)\cup hZ(G),
\]
and that $\fg^{\langle X,Y,Z\rangle}=0$. Then $h\in Z(G)$.
\end{proposition}
\begin{proof}
Set $H=C_G(h)$ and $\fh=\Lie(H)$. By Steinberg's connectedness theorem~\cite[Theorem~8.1]{Steinberg}, the centralizer $C_{G_{\CC}}(h)$ is connected, and the $h$-equivariant polar decomposition identifies $C_{G_{\CC}}(h)$ with $H\times\fh$~\cite[\S2.4(3)]{Bergeron}, so $H$ is connected. Its Lie algebra has zero center, since
\[
\mathfrak z(\fh)\subseteq\fg^{\langle X,Y,Z\rangle}=0.
\]
Thus $H$ is semisimple. Both $h$ and $Z(G)$ lie in $Z(H)$, so the adjoint actions of $X,Y,Z$ on $\fh$ satisfy the spherical $(2,3,k)$ relations. Their common fixed subspace is zero, so Lemma~\ref{lem:minima} gives $\dim\fh^X=d_{2,\fh}$.

Suppose $h\notin Z(G)$. Choose a maximal torus $T$, simple roots $\Delta=\{\alpha_1,\ldots,\alpha_r\}$, and highest root $\theta=\sum_i a_i\alpha_i$. Conjugate $h$ to $\eexp{v}$ with $v$ in the closed fundamental alcove. The roots of $H$ are those $\beta$ with $\beta(v)\in\ZZ$. Their span is the span of the normals of the alcove walls through $v$, so its codimension is the dimension of the smallest face containing $v$. As $H$ has full rank and is semisimple, $v$ is a vertex, and it is nonzero since $h$ is noncentral. The alcove description~\cite[Propositions~3.3--3.4, Remark~3.6, and Theorem~3.7]{KostantBDS} gives
\[
v=\varpi_j^\vee/a_j,\qquad a_j=\operatorname{ord}(\Ad(h))=2,\qquad
\Delta_H=(\Delta\setminus\{\alpha_j\})\cup\{-\theta\},
\]
where $\varpi_j^\vee$ is the $j$th fundamental coweight and $\Delta_H$ is a simple system for $H$.

Every $b\in Z(H)$ belongs to $T$ and satisfies $\alpha_i(b)=1$ for $i\ne j$ and $\theta(b)=1$. The highest-root relation gives $\alpha_j(b)^2=1$ as well. Thus every simple root of $G$ is trivial on $b^2$, so $Z(H)^2\subseteq Z(G)$.

Let $\varrho_H^\vee$ be the half-sum of the positive coroots for $\Delta_H$, and write $h_\Phi=1+\sum_i a_i$ for the Coxeter number of $G$. Every root in $\Delta_H$ pairs to $2$ with $2\varrho_H^\vee$. Applying this to $-\theta$ gives
\[
\pair{\alpha_j}{2\varrho_H^\vee}=-1-\sum_{i\ne j}a_i=2-h_\Phi.
\]
The coefficient $a_j=2$ excludes type $A$, and every other irreducible type has even Coxeter number~\cite[p.~372]{EKV}. Hence every simple root of $G$ pairs evenly with $2\varrho_H^\vee$, and $z_H=(2\varrho_H^\vee)(-1)$ lies in $Z(G)$.

Lemma~\ref{lem:involution}, applied to the connected compact semisimple group $H$, now gives $X^2=b^2z_H$ for some $b\in Z(H)$. Both factors belong to $Z(G)$, whereas $X^2\in hZ(G)$ by hypothesis. This forces $h\in Z(G)$, a contradiction.
\end{proof}

The next proposition computes the necessary dimension condition for
the existence of elliptic projective representations, which will be used for $2T$, $2O$, and $2I$ in Sections~\ref{sec:odd-models}--\ref{sec:full-reciprocity}.

\begin{proposition}\label{prop:defect-screen}
Let $G$ be a simply connected compact simple group of rank $r$, with Lie algebra $\fg$. Choose homogeneous, algebraically independent generators
\[
\CC[\ft_{\CC}]^W=\CC[p_1,\ldots,p_r].
\]
The integers $e_j=\deg p_j-1$ are the exponents of $\fg$. The minimum $d_m$ in Lemma~\ref{lem:minima} is
\begin{equation}\label{eq:exponent-minima}
d_m=r+2\sum_{j=1}^r\left\lfloor\frac{e_j}{m}\right\rfloor
\qquad(m\ge1).
\end{equation}
For $k\in\{3,4,5\}$, set
$\delta_k=d_2+d_3+d_k-\dim\fg$. Then
\begin{equation}\label{eq:exponent-defect}
\delta_k=2\sum_{j=1}^r
\left(
1+\left\lfloor\frac{e_j}{2}\right\rfloor
+\left\lfloor\frac{e_j}{3}\right\rfloor
+\left\lfloor\frac{e_j}{k}\right\rfloor-e_j
\right).
\end{equation}
For $k=3,4,5$, the complete lists of simple types with $\delta_k=0$ are
\begin{center}
\begin{tabular}{ll}
\toprule
$\Gamma$ & Simple types with $\delta_k=0$\\
\midrule
$2T\ (k=3)$ & $A_1,\ A_2,\ G_2$\\
$2O\ (k=4)$ & $A_r\ (1\le r\le3),\ B_r,\ C_r\ (2\le r\le4),\ D_4,\ G_2,\ F_4$\\
$2I\ (k=5)$ & $A_r\ (1\le r\le5),\ B_r,\ C_r\ (2\le r\le7),$\\
 & $D_4,\ D_5,\ D_6,\ D_8,\ G_2,\ F_4,\ E_6,\ E_7,\ E_8$\\
\bottomrule
\end{tabular}
\end{center}
(Here $B_2=C_2$.)
For every type outside the appropriate row, all elliptic projective class
spaces vanish.

Among the tetrahedral types in the table, only $A_2$ has nontrivial
$3$-primary center. Among the octahedral types, only $A_2$ has nontrivial
odd-order center, and only $A_3$ has a $2$-primary center not annihilated
by $2$.
\end{proposition}

\begin{proof}
By Kostant's theorem, the highest weights for the principal triple are twice the exponents defined by the basic invariant degrees. Thus the principal $\mathfrak{sl}_2$ decomposition is
\[
\fg_{\CC}\simeq
\bigoplus_{j=1}^r\Sym^{2e_j}\CC^2;
\]
see~\cite[\S32.3]{EtingofExponents}. The principal element
$s_m=(2\rh)(e^{\pi i/m})$ acts on the weights
$-2e_j,-2e_j+2,\ldots,2e_j$ of the $j$th summand with eigenvalues
$e^{2\pi i a/m}$, $-e_j\le a\le e_j$. Exactly
$1+2\lfloor e_j/m\rfloor$ of these eigenvalues are $1$, proving \eqref{eq:exponent-minima}. Since
$\dim\fg=\sum_j(2e_j+1)$, equation~\eqref{eq:exponent-defect} follows.

Put
\[
a_k(e)=1+\lfloor e/2\rfloor+\lfloor e/3\rfloor+\lfloor e/k\rfloor-e.
\]
The recurrences
\[
a_3(e+6)=a_3(e)+1,\qquad
a_4(e+12)=a_4(e)+1
\]
reduce the calculation to the residue lists
\[
\begin{split}
(a_3(0),\ldots,a_3(5))&=(1,0,0,1,1,0),\\
(a_4(0),\ldots,a_4(11))&=(1,0,0,0,1,0,1,0,1,1,1,0).
\end{split}
\]
Consequently, for $k=3,4$, $a_k(e)\ge0$ for positive integers $e$, and its positive
integer zeros are
\[
\{e:a_3(e)=0\}=\{1,2,5\},\qquad
\{e:a_4(e)=0\}=\{1,2,3,5,7,11\}.
\]
For $k=5$ we have $a_5(e+30)=a_5(e)+1$. In the range $0\le e<30$, the value of $a_5(e)$ is zero precisely at
\[
1,\ 2,\ 3,\ 4,\ 5,\ 7,\ 8,\ 9,\ 11,\ 13,\ 14,\ 17,\ 19,\ 23,\ 29,
\]
and is one at the remaining residues. The recurrence gives $a_5(e)\ge0$ for every $e\ge0$ and shows that these are all its zeros.
Thus $\delta_k=0$ if and only if every exponent lies in the corresponding
zero set.

The exponent lists are standard; see~\cite[\S32.3]{EtingofExponents}
or~\cite[Ch.~VI, tables]{Bourbaki}. In the classical families they are
\begin{center}
\begin{tabular}{ll}
\toprule
$A_r$ & $1,2,\ldots,r$\\
$B_r,\ C_r$ & $1,3,\ldots,2r-1$\\
$D_r$ & $1,3,\ldots,2r-3,\ r-1$\\
\bottomrule
\end{tabular}
\end{center}
For $k=3$, the exponent $3$ excludes $A_r$ for $r\ge3$,
$B_r,C_r$ for $r\ge2$, and every $D_r$ for $r\ge4$.
For $k=4$, the exponent $4$ excludes $A_r$ for $r\ge4$ and $D_5$,
while the exponent $9$ excludes $B_r,C_r$ for $r\ge5$ and $D_r$
for $r\ge6$; $D_4$ has exponents $1,3,3,5$ and thus remains.
For $k=5$, $A_r$ first acquires the excluded exponent $6$ at $r=6$, and $B_r,C_r$ first acquire $15$ at $r=8$. In type $D_r$ the odd string $1,3,\ldots,2r-3$ contains $15$ for $r\ge9$. Among $4\le r\le8$, the additional exponent $r-1$ belongs to the zero set except at $r=7$, when $r-1=6$. This leaves $D_4,D_5,D_6,D_8$.
The exceptional exponent lists are
\begin{center}
\begin{tabular}{ll}
\toprule
$G_2$ & $1,5$\\
$F_4$ & $1,5,7,11$\\
$E_6$ & $1,4,5,7,8,11$\\
$E_7$ & $1,5,7,9,11,13,17$\\
$E_8$ & $1,7,11,13,17,19,23,29$\\
\bottomrule
\end{tabular}
\end{center}
which give the exceptional entries in the table. For $k=5$ every exponent in all five exceptional rows belongs to the zero set.

Let $f$ be a projective map with actual grade $b\in B$, and write $h=f(-1)$ and $(X,Y,Z)=(f(x),f(y),f(z))$. Since $-1$ has trivial image in the abelianization, $\ell(-1)=0$ and $\ell((-1)g)=\ell(g)$ for every $g$. Thus $c_b(-1,g)=c_b(g,-1)=1$, so $h^2=1$ and $h$ commutes with $X,Y,Z$. The binary spherical relations give $X^2,Y^3,Z^k,XYZ\in hB$. If $f$ is elliptic, then $\fg^{\langle X,Y,Z\rangle}=0$, since these three elements generate its adjoint image. Proposition~\ref{prop:central-involution} forces $h\in Z(G)$, so the adjoint action descends to the rotation group, and Lemma~\ref{lem:minima} gives $\delta_k=0$. Hence positive defect excludes elliptic classes in every actual grade.

Finally, the simply connected centers of the types in the tetrahedral and octahedral rows are
\begin{center}
\begin{tabular}{lccccccc}
\toprule
Type & $A_1$ & $A_2$ & $A_3$ & $B_r,\ C_r$ & $D_4$ & $G_2$ & $F_4$\\
\midrule
$Z(G)$ & $C_2$ & $C_3$ & $C_4$ & $C_2$ & $C_2\times C_2$ & $1$ & $1$\\
\bottomrule
\end{tabular}
\end{center}
These follow from the root and weight lattices in~\cite[Ch.~VI, tables]{Bourbaki};
see also~\cite[21.10]{Milne}. The claims about the relevant primary
centers follow immediately.
\end{proof}

\begin{remark}\label{rem:screening-use}
The first two rows of Proposition~\ref{prop:defect-screen} isolate the simple factors that require a central-coefficient calculation for $2T$ and $2O$. For $2I$ the abelianization is trivial, so the chosen carry cocycle is trivial in every actual grade. At positive defect the elliptic spaces vanish, and at zero defect ($\delta_5=0$), Sections~\ref{sec:odd-models} and~\ref{sec:full-reciprocity} compare the ordinary local class models uniformly. The third row identifies the simple groups with an ordinary elliptic $2I$-class. Indeed, the reference action in the proof of Lemma~\ref{lem:minima} gives $\delta_5(G)=2\dim\fg^{\Ad(\iota(2I))}$ for the restriction of the principal homomorphism $\iota\colon\SU(2)\to G$, so at zero defect this restriction already gives one ordinary elliptic class. Theorem~\ref{thm:local-realization} gives the finer realization and uniqueness for prescribed local classes.

For $2T$, the remaining types $A_1,A_2,G_2$ have centers $C_2,C_3,1$,
respectively. For $A_1$ and $G_2$, the center has order prime to $3$,
so they are covered by Remark~\ref{rem:coprime-packet}; thus only $A_2$, namely $\SU(3)$, requires a separate central-coefficient
calculation, carried out in Section~\ref{sec:tetrahedral}.

For $2O$, all types except $A_2$ and $A_3$ have centers of
exponent at most $2$. In Section~\ref{sec:octahedral}, these cases are
handled uniformly by the local quartic Fourier transform, while
$\SU(3)$ and $\SU(4)$, whose centers are $C_3$ and $C_4$, respectively,
will be treated separately.
\end{remark}

\section{Realization of local spherical classes via equivariant bundles on \texorpdfstring{$\PP^1$}{P1}}\label{sec:realization}

The main input from this section to the polyhedral Fourier comparison in Sections~\ref{sec:odd-models} and~\ref{sec:tetrahedral}--\ref{sec:full-reciprocity} is Corollary~\ref{cor:spherical-local}, which reduces global elliptic conjugacy classes for $2T,2O,2I$ to triples of minimal local conjugacy classes. We deduce it from a general local-to-global realization theorem for central extensions of spherical rotation groups.

\subsection{Main theorem}

Consider a central extension of finite groups
\[
1 \longrightarrow K \longrightarrow \Gamma
\longrightarrow H \longrightarrow 1,
\]
where $K$ is a finite group, $H=\langle \bar x,\bar y,\bar z\mid \bar x^2=\bar y^3=\bar z^k=\bar x\bar y\bar z=1\rangle \subset \SO(3)$ acts faithfully on
$\PP^1 \simeq S^2$, the quotient is again $\PP^1$,
and the quotient map has three branch points.
Choose representatives $p_i$ of the three exceptional orbits, and write
\[
\Gamma_i = \Stab_{\Gamma}(p_i),
\qquad
H_i = \Gamma_i/K \simeq C_{m_i}
\qquad (i=1,2,3).
\]
The quotient $\Gamma_i/K$ is cyclic even when the stabilizer $\Gamma_i$ itself is not cyclic, as can occur in the central extension in Lemma~\ref{lem:octahedral-local-model}.

Let $G$ be a simply connected compact semisimple Lie group,
let $\fg=\Lie(G)$, and let $G_{\CC}$
be its complexification. Fix
\[
\chi \colon K \longrightarrow Z(G),
\qquad
\phi_i \colon \Gamma_i \longrightarrow G,
\qquad
\phi_i|_K = \chi,
\]
where only the $G$-conjugacy class of each $\phi_i$ is prescribed.
Set
\[
d = \dim_{\RR}\fg
  = \dim_{\CC}\fg_{\CC},
\qquad
d_i = \dim_{\RR}
      \fg^{\phi_i(\Gamma_i)}.
\]
In this section, we prove the following theorem characterizing elliptic representations in terms of their local conjugacy classes.

\begin{theorem}\label{thm:local-realization}
If
\[
d_1+d_2+d_3=d,
\]
then there is a unique $G$-conjugacy class $[\rho]$ of homomorphisms
$\rho \colon \Gamma \to G$ satisfying
\[
\rho|_K = \chi,
\qquad
[\rho|_{\Gamma_i}] = [\phi_i].
\]
Moreover, $\fg^{\rho(\Gamma)}=0$; thus the resulting
homomorphism is elliptic.
\end{theorem}

Appendix~\ref{app:yang-mills} gives a second proof of existence and ellipticity using an invariant connection.

For $h\in Z(G)$ and $m\ge1$, set
\[
C_{h,m}(G)=\bigl\{[t]_G:t^m=h,\ \dim_{\RR}\fg^{\Ad(t)}=d_m(G)\bigr\}.
\]
Here $d_m(G)$ is the minimum over all $t\in G$ with $t^m\in Z(G)$, as in Lemma~\ref{lem:minima}. For a specified $h$ this set can therefore be empty. Write $E_\Gamma(G)$ for the elliptic conjugacy classes in $\Hom(\Gamma,G)/G$.

\begin{corollary}\label{cor:spherical-local}
Let $(\Gamma,k)$ be $(2T,3)$, $(2O,4)$ or $(2I,5)$. If $d_2+d_3+d_k=\dim\fg$, restriction to the standard binary spherical generators gives the bijection
\[
E_\Gamma(G)\xrightarrow{\sim}\coprod_{h\in Z(G)[2]} C_{h,2}(G)\times C_{h,3}(G)\times C_{h,k}(G),\qquad
[\rho]\longmapsto \bigl([\rho(x)],[\rho(y)],[\rho(z)]\bigr)_{h=\rho(-1)}.
\]
\end{corollary}

\begin{proof}
Each projection of an elliptic homomorphism to a simple factor of $G$ is elliptic, so Proposition~\ref{prop:central-involution} on those factors gives $\rho(-1)\in Z(G)[2]$. The adjoint action then descends to the rotation group, and Lemma~\ref{lem:minima} gives minimality at all three stabilizers. Conversely, in the summand indexed by $h$, the three local homomorphisms agree on $-1$, and their fixed dimensions sum to $\dim\fg$. Theorem~\ref{thm:local-realization} realizes them by a unique elliptic conjugacy class.
\end{proof}

The remainder of this section is devoted to the proof of Theorem~\ref{thm:local-realization}.

\subsection{Holomorphic \texorpdfstring{$\Gamma$}{Gamma}-equivariant principal \texorpdfstring{$G_{\CC}$}{GC}-bundle}

\begin{lemma}\label{lem:equivariant-bundle-construction} Let $X=\PP^1$ with the $\Gamma$-action above, of central kernel $K$, and put $H=\Gamma/K$. Let $p_1,p_2,p_3$ represent its three exceptional orbits, with stabilizers $\Gamma_i=\Stab_\Gamma(p_i)$. Suppose $G$ is a compact connected Lie group and we are given homomorphisms
\[
\chi\colon K\longrightarrow Z(G),
\qquad
\phi_i\colon\Gamma_i\longrightarrow G,
\qquad
\phi_i|_K=\chi.
\]
Then there exists a smooth $\Gamma$-equivariant principal
$G$-bundle $P\to X$ such that
\[
k\cdot e=e\,\chi(k)
\qquad(k\in K,\ e\in P),
\]
and for each $i$, there is an element $e_i\in P_{p_i}$ satisfying
\[
\gamma\cdot e_i=e_i\,\phi_i(\gamma)
\qquad(\gamma\in\Gamma_i),
\]
where the $\gamma$-action is induced by the $\Gamma$-equivariant structure and the $\phi_i(\gamma)$-action is induced by the principal $G$-bundle structure. Moreover, the complexification
\[
P_{\CC}=P\times_G G_{\CC}
\]
admits a holomorphic principal $G_{\CC}$-bundle structure for which the $\Gamma$-action is holomorphic.
\end{lemma}

\begin{proof}
Let $q\colon X\to Y=X/H$. Choose disjoint open disks $U_i$ about the three branch values, and a component $D_i$ of $q^{-1}(U_i)$ containing $p_i$. Over $q^{-1}(U_i)$ take the induced bundle $\Gamma\times_{\Gamma_i}(D_i\times G)$, with $\gamma(u,g)=(\gamma u,\phi_i(\gamma)g)$ for $\gamma\in\Gamma_i$. At $p_i$ the frame represented by $(1,p_i,1)$ realizes $\phi_i$.

Choose smaller concentric disks with $\overline{U_i'}\subset U_i$. Set $B=Y\setminus\bigcup_i\overline{U_i'}$ and $A_i=B\cap U_i$, an open annulus. The action gives a principal $H$-cover $q^{-1}(B)\to B$. Since $B$ retracts onto a pair of circles joined at one point, its fundamental group is free, so we may lift the monodromies of its free generators from $H$ to $\Gamma$. The lifted homomorphism determines a principal $\Gamma$-cover $Q\to B$, possibly disconnected, and an identification of $Q/K$ with $q^{-1}(B)$ as principal $H$-covers. Write the deck action on the left and form
\[
P_B=(Q\times G)/\bigl((ku,g)\sim(u,\chi(k)g)\bigr)\longrightarrow Q/K.
\]
Centrality of $K$ makes $\gamma[u,g]=[\gamma u,g]$ well defined. As $\chi(k)\in Z(G)$, the action of $k$ is right multiplication by $\chi(k)$, as it is on each disk model.

On $q^{-1}(A_i)$ form the smooth bundle of right $G$-equivariant isomorphisms from the disk-model fibers to the fibers of $P_B$, whose fiber is diffeomorphic to $G$. Each isomorphism commutes with the common $K$-action, so the action on this bundle factors through $H$, and the free $H$-cover descends it to a smooth fiber bundle on $A_i$. Smooth parallel transport along radial segments identifies the latter with the pullback of its restriction to a core circle (its fiber $G$ is compact). Over the circle, a cut and parallel transport give a mapping torus with identification $(1,g)\sim(0,\psi(g))$. A smooth path from $g_0$ to $\psi^{-1}(g_0)$, constant near both endpoints, defines a smooth section because $G$ is connected. Extend radially and pull back to obtain a $\Gamma$-equivariant smooth isomorphism on $q^{-1}(A_i)$. The three disjoint annuli now glue the models to a smooth equivariant principal bundle $P\to X$ with the required fiber actions.

Extend the action to the complexification $P_{\CC}=P\times_G G_{\CC}$ by $\gamma[p,g]=[\gamma p,g]$. The space of smooth connections on $P$ is affine, so averaging a connection over $\Gamma$ gives an invariant connection. Its complexification has vanishing $(0,2)$-curvature because $X$ is a complex curve, so its $(0,1)$-part defines a holomorphic structure on $P_{\CC}$~\cite[Lemmas~2.4--2.5]{Trautwein}; invariance makes the $\Gamma$-action holomorphic.

\end{proof}

\subsection{Associated bundle}

Given a holomorphic $\Gamma$-equivariant principal $G_{\CC}$-bundle $P_{\CC}\to X$, define its associated adjoint bundle via $\Ad\colon G_{\CC}\to\operatorname{GL}(\fg_{\CC})$ by
\[
F=\fg_{\CC}\overset{G_{\CC}}{\times} P_{\CC}.
\]

\begin{lemma}
\label{lem:invariant-riemann-roch}
Set
\[
d_i=\dim_{\CC}(F_{p_i})^{\Gamma_i}=\dim \fg^{\Gamma_i}.
\]
Then
\[
\dim_{\CC}H^0(X,F)^\Gamma-\dim_{\CC}H^1(X,F)^\Gamma=\frac{d_1+d_2+d_3-d}{2}.
\]
In particular, if $d_1+d_2+d_3=d$, then
\[
\dim_{\CC}H^0(X,F)^\Gamma=\dim_{\CC}H^1(X,F)^\Gamma.
\]
\end{lemma}

\begin{proof}
Since $\chi(K)$ is central, $K$ acts trivially on $F$, so the $\Gamma$-action on $F$ factors through $H$. Write
\[
\pi\colon X\longrightarrow Y=X/H\simeq\PP^1,\qquad N=|H|,\qquad \mathcal F=(\pi_*F)^H.
\]
Put $H_i=\Gamma_i/K$ and $m_i=|H_i|$. Choose a local coordinate $u$ centered at $p_i$ such that a generator $x$ of $H_i$ acts by $u\mapsto\zeta_i u$, where $\zeta_i=e^{2\pi i/m_i}$, so that the quotient coordinate is $v=u^{m_i}$. Shrinking the disk if necessary, we obtain an $H_i$-equivariant holomorphic trivialization of $F$. Write
\[
F_{p_i}=\bigoplus_{a=0}^{m_i-1}F_{i,a},
\]
where $x$ acts on $F_{i,a}$ by $\zeta_i^a\id$. Set $n_{i,a}=\dim_{\CC}F_{i,a}$.

For $e\in F_{i,a}$, an invariant section $f(u)e$ satisfies $f(\zeta_i u)=\zeta_i^a f(u)$, and hence has the form $f(u)=u^a h(u^{m_i})$ for a holomorphic function $h$. Thus the invariant sections form a free $\cO_Y$-module generated by $u^a e$. Together with the usual local description away from the branch points, this proves that $\mathcal F$ is locally free of rank $d$.

The natural evaluation map $\pi^*\mathcal F\longrightarrow F$ is an isomorphism away from the ramification points and is multiplication by $u^a$ on the above eigenspaces. There are $N/m_i$ points in the orbit of $p_i$. Comparing determinants gives
\[
\deg F=N\deg\mathcal F+\sum_{i=1}^3\frac{N}{m_i}\sum_{a=0}^{m_i-1}a\,n_{i,a}.
\]
The Killing form induces an equivariant isomorphism $F\simeq F^*$, so $\deg F=0$. Moreover, the pairing between $F_{i,a}$ and $F_{i,m_i-a}$ for $1\le a<m_i$ is nondegenerate, so $n_{i,a}=n_{i,m_i-a}$, while $n_{i,0}=d_i$. Therefore
\[
\sum_{a=1}^{m_i-1}\frac{a}{m_i}n_{i,a}=\frac{1}{2}\sum_{a=1}^{m_i-1}n_{i,a}=\frac{d-d_i}{2},
\]
and it follows that
\[
\deg\mathcal F=-\frac12\sum_{i=1}^3(d-d_i).
\]

Since $\pi$ is finite, its higher direct images on coherent sheaves vanish. As in the proof of Lemma~\ref{lem:minima}, the averaging projector shows that taking $H$-invariants is exact. Thus
\[
H^q(Y,\mathcal F)\simeq H^q(X,F)^H=H^q(X,F)^\Gamma\qquad(q=0,1).
\]
Applying Riemann--Roch on $Y\simeq\PP^1$, we obtain
\[
\begin{aligned}\dim_{\CC}H^0(X,F)^\Gamma-\dim_{\CC}H^1(X,F)^\Gamma&=d+\deg\mathcal F\\
&=d-\frac12\sum_{i=1}^3(d-d_i)\\
&=\frac{d_1+d_2+d_3-d}{2}.
\end{aligned}
\]
\end{proof}

\subsection{Equivariant Birkhoff--Grothendieck decomposition}

Let $E$ denote a $\Gamma$-invariant holomorphic structure on $P_{\CC}$, and put $F=\ad E$. In any Birkhoff--Grothendieck splitting $F\simeq\bigoplus_j\cO(a_j)$, the sum of the terms of degree at least $n$ defines a canonical subbundle $F_{\ge n}$. Indeed, a homomorphism $\cO(u)\to\cO(v)$ vanishes for $u>v$, so changes of splitting preserve this filtration. The same vanishing applied to the holomorphic bracket gives
\[
[F_{\ge n},F_{\ge m}]\subseteq F_{\ge n+m}.
\]
Fix a maximal torus $T_{\CC}\subset G_{\CC}$ and positive roots in its root system $\Phi$.

By Grothendieck's classification theorem~\cite[Th\'eor\`eme~1.1 and Corollaire~2]{Grothendieck}, $E$ admits a holomorphic reduction to $T_{\CC}$ and has a unique dominant cocharacter type $\mu$ with
\[
F\simeq(\cO_X\otimes\ft_{\CC})\oplus\bigoplus_{\alpha\in\Phi}\cO_X(\pair{\alpha}{\mu}).
\]
The bundle $E$ is holomorphically trivial if and only if $\mu(E)=0$.

\begin{lemma}\label{lem:nonzero-type-invariants}
If $\mu(E)\ne0$, there is an invariant holomorphic section $z\in H^0(X,F_{\ge0})^\Gamma$ whose bracket induces multiplication by $n$ on each graded quotient $F_n=F_{\ge n}/F_{\ge n+1}$. In particular $z\ne0$.
\end{lemma}

\begin{proof}
The canonical filtration is $\Gamma$-invariant, and the torus reduction identifies the fibers of $\operatorname{gr}F=\bigoplus_n F_n$ with $\fg_{\CC}$ as graded Lie algebras. The operator $D|_{F_n}=n\,\id_{F_n}$ is an invariant holomorphic derivation. Every derivation of a semisimple Lie algebra is inner~\cite[\S5.3]{HumphreysLie}. Since the center is zero, the fiberwise maps $\ad$ identify $\operatorname{gr}F$ holomorphically with its derivation bundle, so $D=\ad(z_{\mathrm{gr}})$ for a unique invariant holomorphic section. The identity $Dz_{\mathrm{gr}}=[z_{\mathrm{gr}},z_{\mathrm{gr}}]=0$ places it in $F_0$, and some root has nonzero pairing with $\mu(E)$, so $D\ne0$.

The bundle $F_{\ge1}$ is a sum of line bundles of positive degree, so $H^1(X,F_{\ge1})=0$. The sequence
\[
0\longrightarrow F_{\ge1}\longrightarrow F_{\ge0}\longrightarrow F_0\longrightarrow0
\]
therefore gives a holomorphic lift of $z_{\mathrm{gr}}$. Average this lift over $\Gamma$ to obtain $z$; its image is still $z_{\mathrm{gr}}$, and compatibility of the bracket with the filtration gives the required action on every $F_n$.
\end{proof}

\subsection{Triviality of a minimizing holomorphic structure}

\begin{lemma}\label{lem:minimal-type-zero}
Assume $d_1+d_2+d_3=d$. Among all $\Gamma$-invariant holomorphic structures $E$ on the fixed smooth bundle $P_{\CC}$, choose $E_0$ minimizing $\dim_{\CC}H^0(X,\ad E)$. Then $\mu(E_0)=0$, so $E_0$ is holomorphically trivial.
\end{lemma}

\begin{proof}
Lemma~\ref{lem:equivariant-bundle-construction} gives at least one such structure, and the dimensions are nonnegative integers, so a minimizer exists. Suppose $\mu(E_0)\ne0$ and put $F=\ad E_0$. Take the nonzero invariant section $z$ given by Lemma~\ref{lem:nonzero-type-invariants}; Lemma~\ref{lem:invariant-riemann-roch} and the dimension hypothesis give $H^1(X,F)^\Gamma\ne0$.

The endomorphism $\ad z$ preserves the degree filtration and induces $n\,\id$ on $F_n$. In a holomorphic splitting compatible with the filtration, the off-diagonal terms strictly raise degree, and the same triangular form persists on $H^1(X,F)$. Each $F_n$ is a sum of copies of $\cO(n)$, so a nonzero $H^1(X,F_n)$ requires $n\le-2$. Consequently all diagonal blocks on $H^1(X,F)$ are invertible, and $(\ad z)_*$ is an isomorphism there.

Choose a nonzero class $[a]\in H^1(X,F)^\Gamma$ and average a Dolbeault representative to obtain $a\in\Omega^{0,1}(X,F)^\Gamma$. On the pullback of $P_{\CC}$ to $X\times\Delta$, where $\Delta$ is a disk with coordinate $t$, consider the principal $(0,1)$-operator
\[
\bar\partial_{E_0}+\bar\partial_\Delta+t a.
\]
Its $(0,2)$-curvature along $X$ vanishes because $X$ is a curve, and the mixed term also vanishes because $t$ is holomorphic and $a$ is independent of $t$. The integrability theorem for principal bundles~\cite[Theorem~A.22]{TelemanHolomorphicBundles} yields a holomorphic family $\mathcal E$ with fibers $E_t$. The fixed smooth $\Gamma$-action preserves these structures since $a$ is invariant.

Put $\mathcal F=\ad\mathcal E$ and $F_t=\mathcal F|_{X\times\{t\}}$. Upper semicontinuity gives $h^0(X,F_t)\le h^0(X,F)$ near $0$, and minimality of $E_0$ gives the reverse inequality. For the proper projection $p\colon X\times\Delta\to\Delta$, Grauert's base-change theorem~\cite[Ch.~I, \S1]{OSS} therefore makes $p_*\mathcal F$ locally free with fibers $H^0(X,F_t)$. After shrinking the disk, the chosen section $z$ extends to a holomorphic section of this direct image, with fiber values $z_t$ and $z_0=z$. In the fixed smooth bundle the equation and its derivative at $0$ read
\[
\bar\partial_F z_t+t[a,z_t]=0,\qquad
\bar\partial_F\dot z_0+[a,z]=0.
\]
Thus $(\ad z)_*[a]=0$ in $H^1(X,F)$, since $[z,a]=-[a,z]$ is exact. Invertibility forces $[a]=0$, a contradiction, so $\mu(E_0)=0$.
\end{proof}

\subsection{Proof of the main theorem}

\begin{proof}[Proof of Theorem~\ref{thm:local-realization}]
Lemmas~\ref{lem:equivariant-bundle-construction} and~\ref{lem:minimal-type-zero} give a holomorphically trivial bundle with the prescribed equivariant fiber actions. In a holomorphic trivialization over $X=\PP^1$, each $\gamma\in\Gamma$ acts by $(x,g)\mapsto(\gamma x,f_\gamma(x)g)$. Since $G_{\CC}$ is affine and $X$ is compact and connected, the holomorphic map $f_\gamma\colon X\to G_{\CC}$ is constant, and its value $\rho(\gamma)$ satisfies the group law by equivariance. The fiber actions give $\rho|_K=\chi$ and $G_{\CC}$-conjugacy of $\rho|_{\Gamma_i}$ with $\phi_i$. By Lemma~\ref{lem:compact-reduction}(1) the finite-image homomorphism is conjugate into $G$, and its local restrictions are then $G$-conjugate to $\phi_i$. The adjoint action factors through $H$ since $\rho(K)$ is central. Formula~\eqref{eq:spherical-euler} for $V=\fg_{\CC}$ gives
\[
d_1+d_2+d_3=d+2\dim_{\RR}\fg^{\rho(\Gamma)},
\]
which by the dimension hypothesis proves ellipticity.

Choose oriented loops around the three branch values from one base point in the punctured quotient, in an order for which their product is trivial. Their monodromies give compatible generators $\bar x_i$ for the spherical presentation of $H$, after conjugating the chosen stabilizers and transporting their local data if necessary. Choose lifts $x_i\in\Gamma_i$ and write $x_i^{m_i}=k_i$ and $x_1x_2x_3=k_0$, with $k_i,k_0\in K$. Then $\Gamma_i=\langle K,x_i\rangle$, and a presentation of $\Gamma$ is obtained by adjoining $x_1,x_2,x_3$ to $K$ with relations
\[
[K,x_i]=1,\qquad x_i^{m_i}=k_i,\qquad x_1x_2x_3=k_0.
\]
Indeed the presented group surjects onto $\Gamma$ and induces the spherical presentation of $H$ upon quotienting by the image of $K$. Any kernel thus lies in that image, on which the map is injective by its prescribed restriction to $K$.

Let $C_i=[\phi_i(x_i)]_G\subset G$. Every $X_i\in C_i$ already satisfies $X_i^{m_i}=\chi(k_i)$ because $\chi(k_i)$ is central. Since $x_i$ and $K$ generate $\Gamma_i$, membership in $C_i$ fixes the entire local conjugacy class. The remaining relation identifies the required conjugacy classes with
\[
\Phi^{-1}(\chi(k_0))/G,\qquad
\Phi\colon C_1\times C_2\times C_3\longrightarrow G,\quad
(X_1,X_2,X_3)\longmapsto X_1X_2X_3.
\]
The product of conjugacy classes is a compact connected quasi-Hamiltonian $G$-space with fused moment map $\Phi$~\cite[Proposition~3.1 and Theorem~6.1]{AMM}. Since $G$ is simply connected, its moment-map fibers are connected~\cite[Theorem~7.2]{AMM}. The value $\chi(k_0)$ is central, so simultaneous conjugation preserves its fiber. The quotient is connected and Hausdorff, finite by Lemma~\ref{lem:compact-reduction}(2), and nonempty by existence, so it consists of a single point.
\end{proof}

\section{Ordinary elliptic classes for \texorpdfstring{$2T$ and $2I$}{2T and 2I}}\label{sec:odd-models}

Let $G$ be a simply connected compact semisimple group, and write $\cM_m(G)=C_{1,m}(G)$ and $B=Z(G)$. If $h\in B[2]$ and $m$ is odd, multiplication by $h$ gives the canonical bijection
\[
C_{h,m}(G)\xrightarrow{\sim}\cM_m(G),\qquad
[t]_G\longmapsto[ht]_G.
\]
Indeed, $(ht)^m=h^mt^m=h^2=1$, and the inverse is again multiplication by $h$; central translation preserves the adjoint fixed subalgebra. Corollary~\ref{cor:squares} shows that the quadratic factor $C_{h,2}(G)$ consists of one class if and only if $h$ lies in
\[
H_G=z_GB^2\cap B[2],
\]
and is empty otherwise.

For $\Gamma=2T$ take $k=3$, and for $\Gamma=2I$ take $k=5$. Suppose that $d_2+d_3+d_k=\dim\fg$. Put $h=\rho(-1)$, $Y=\rho(y)$ and $Z=\rho(z)$. In the disjoint union of Corollary~\ref{cor:spherical-local}, the square condition selects the indices $h\in H_G$ and determines the quadratic class uniquely for each such $h$. Applying the two odd-power bijections gives
\begin{equation}\label{eq:spherical-model}
E_\Gamma(G)\xrightarrow{\sim}H_G\times\cM_3(G)\times\cM_k(G),
\qquad [\rho]\longmapsto\bigl(h,[hY]_G,[hZ]_G\bigr).
\end{equation}
Writing $M_m(G)=\CC[\cM_m(G)]$, linear extension in the class bases gives
\[
\CC[E_\Gamma(G)]\simeq\CC[H_G]\otimes M_3(G)\otimes M_k(G).
\]

\subsection{Duality for the odd factors \texorpdfstring{$\cM_m(G)$}{Mm(G)}}

Let $G^\dagger=(G/Z(G))^\vee$ be the simply connected compact group with dual root system. Write $Q^\vee\subset X_*(T)\subset P^\vee$ for the coroot, cocharacter and coweight lattices of $G$, and $Q\subset X^*(T)\subset P$ for its root, character and weight lattices.

\begin{lemma}\label{lem:odd-local}
Let $G$ be a simply connected compact semisimple group and $m\in\{3,5\}$. There is a canonical bijection $\cM_m(G)\xrightarrow{\sim}\cM_m(G^\dagger)$, equivariant under based root-datum automorphisms and permutations of isomorphic simple factors.
\end{lemma}
\begin{proof}
Evaluation at $\exp(2\pi i/m)$ identifies the quotient of the cocharacter lattice by $m$ with the $m$-torsion in the maximal torus. For $G$ the cocharacter lattice is $Q^\vee$, and for the simply connected group $G^\dagger$ it is $Q$. Conjugacy classes of elements with $m$th power one are therefore parametrized by the Weyl orbits in $Q^\vee/mQ^\vee$ and $Q/mQ$, respectively.

On $E=Q^\vee\otimes\RR$ take the $W$-invariant real inner product for which the simple factors are orthogonal and the short coroots in each irreducible component have squared length two. Let $D(v)=(v,\,\cdot\,)\in E^*$. The root reflection formula gives
\[
D(\alpha^\vee)=d_\alpha\alpha,\qquad
d_\alpha=(\alpha^\vee,\alpha^\vee)/2\in\{1,2,3\},\qquad
\pair{Dv}{\alpha^\vee}=d_\alpha\pair{\alpha}{v}.
\]
Thus $D$ sends $Q^\vee$ to $Q$, and in the bases of simple coroots and simple roots its matrix is diagonal. On any factor where its diagonal entries are units modulo $m$, reduction gives a $W$-equivariant isomorphism $Q^\vee/mQ^\vee\to Q/mQ$. The $\alpha$-root space is fixed by the torus element represented by $v$ if and only if $\pair{\alpha}{v}\equiv0\pmod m$, which by the displayed identity is equivalent to the dual root space for $\alpha^\vee$ being fixed by the element represented by $Dv$. Adding the common rank to the number of fixed roots gives the same centralizer dimension on both sides.

The only exceptional pair is a factor of type $G_2$ with $m=3$. Label its simple roots $\alpha_1$ (short) and $\alpha_2$ (long). Interchanging the nodes identifies the Cartan matrix with its transpose. The resulting isomorphism to the dual based root datum sends $\alpha_i$ to $\alpha_{3-i}^\vee$ on root lattices and $\alpha_i^\vee$ to $\alpha_{3-i}$ on cocharacter lattices. Reducing the latter integral isomorphism modulo three gives the desired bijection on torus elements, which intertwines $s_i$ with $s_{3-i}^\vee$ and preserves root evaluations under the former correspondence, hence also the fixed-root count.

Taking products and passing to Weyl orbits gives a bijection preserving centralizer dimensions. Dual root systems have the same exponents, so Proposition~\ref{prop:defect-screen} gives the same $d_m$ on both sides. Restricting to the classes attaining this value yields the desired bijection. Since both the normalized metrics and the correspondence with the dual $G_2$ diagram are determined factor by factor from the based data, the map has the required equivariance.
\end{proof}

\subsection{Duality for the central factor \texorpdfstring{$H_G$}{HG}}

Next we compare $H_G$ and $H_{G^\dagger}$. Setting $S_G=Z(G)^2\cap Z(G)[2]$, we have $H_G=z_GS_G$, since $z_G\in Z(G)[2]$.

\begin{proposition}\label{prop:icosahedral-center}
For every simply connected compact simple group $G$, $S_G=Z(G)^2\cap Z(G)[2]$ has order at most $2$, and $|S_G|=|S_{G^\dagger}|$.

Translating the unique isomorphism $S_G\simeq S_{G^\dagger}$ by the canonical basepoints $z_G$ and $z_{G^\dagger}$ gives a canonical bijection
\[
H_G\xrightarrow{\sim}H_{G^\dagger}.
\]
\end{proposition}

\begin{proof}
Although this can be read off from the classification of centers~\cite[Ch.~24, \S c]{MilneAG2017}, we give a direct proof via the Cartan matrix $C$.

We first show that $|S_G|\le2$. The Smith normal form of $C$ gives the cyclic factors of $Z(G)$, and $S_G$ has one $C_2$ factor for each cyclic factor whose order is divisible by $4$. Recall that a connected Dynkin diagram is either a path with at most one multiple edge or a simply laced tree with a single trivalent vertex.

In the path case, order the vertices so that $C_{i,i+1}=-1$ for every $i$. Deleting the last row and the first column of $C$ gives a triangular minor with determinant $\pm1$. Hence $\coker(C)$ is cyclic by Smith normal form, and so is $Z(G)$. In particular, $|S_G|\le2$.

In the trivalent case, let $1\le a\le b\le c$ be the lengths of the three arms. A chain of length $a$ has Cartan determinant $a+1$ and contributes $a/(a+1)$ to the Schur complement at the trivalent vertex. Eliminating the three arms therefore gives
\[
\det C=(a+1)(b+1)(c+1)\left(\frac{1}{a+1}+\frac{1}{b+1}+\frac{1}{c+1}-1\right)>0.
\]
It follows that $a=1$ and $b\le2$. If $b=1$, then $\det C=4$. If $b=2$, then $\det C=5-c\le3$. Since $|Z(G)|=\det C$, the center is cyclic or has order four. In each case $|S_G|\le2$.

Finally, the Cartan matrix of $G^\dagger$ is $C^t$. Since $C$ and $C^t$ have the same Smith normal form, $Z(G)$ and $Z(G^\dagger)$ are isomorphic as abstract groups. Thus $S_G$ and $S_{G^\dagger}$ have the same order, at most two, so there is a unique isomorphism $u\colon S_G\xrightarrow{\sim}S_{G^\dagger}$. Using $H_G=z_GS_G$, we obtain the canonical bijection
\[
H_G\xrightarrow{\sim}H_{G^\dagger},\qquad z_Gs\longmapsto z_{G^\dagger}u(s).
\]
It is equivariant under based automorphisms, since they fix the distinguished basepoints and act trivially on a group of order at most two. For semisimple $G$, taking the product of these bijections over the simple factors gives the corresponding bijection compatibly with permutations of isomorphic factors.
\end{proof}

For $\Gamma=2T,2I$, the exponent formula in Proposition~\ref{prop:defect-screen} gives $\delta_k(G)=\delta_k(G^\dagger)$. At zero defect, take the product of the bijections for $H_G$ over the simple factors and the two bijections of Lemma~\ref{lem:odd-local} for $m=3,k$. Formula~\eqref{eq:spherical-model} then gives $E_\Gamma(G)\xrightarrow{\sim}E_\Gamma(G^\dagger)$, compatibly with based automorphisms and factor permutations. At positive defect both class sets are empty, by the factorwise vanishing in Proposition~\ref{prop:defect-screen}.

\section{Fourier transform of minimal cyclic classes}\label{sec:minimum-fourier}

Let $G$ be a simply connected compact semisimple group, and set $B=Z(G)=P^\vee/Q^\vee$, $G^\dagger=(G/Z(G))^\vee$, and $B^\dagger=Z(G^\dagger)$. Recall that $z_G=(2\rh)(-1)\in B$.

Define
\[
\cM_m^{\mathrm{cen}}(G):=\left\{[t]_G:t^m\in Z(G),\ \dim_{\RR}\fg^{\Ad(t)}=d_m(G)\right\}=\coprod_{h\in Z(G)} C_{h,m}(G).
\]
Equip $\CC[\cM_m^{\mathrm{cen}}(G)]$ with the Hermitian inner product for which the class basis $e_{[t]}$ is orthonormal. Define its shifted power grading by
\[
q_G([t])=t^m z_G^{-1}\in B.
\]
For $b\in B$ and $\chi\in B^\dagger=\Hom(B,\UU(1))$, define
\[
T_b e_{[t]}=e_{[bt]},\qquad
M_\chi e_{[t]}=\chi\bigl(q_G([t])\bigr)e_{[t]}.
\]
On the dual side, the corresponding operators are given by
\[
T_\chi^\dagger e_{[s]}=e_{[\chi s]},\qquad
M_b^\dagger e_{[s]}=q_{G^\dagger}([s])(b)e_{[s]},
\]
where $q_{G^\dagger}([s])=s^m z_{G^\dagger}^{-1}\in B^\dagger$.

\begin{theorem}\label{thm:minimum-fourier}
For every integer $m\ge 1$, there is a canonical unitary Fourier isomorphism
\[
\mathcal F\colon\CC[\cM_m^{\mathrm{cen}}(G)]\xrightarrow{\sim}\CC[\cM_m^{\mathrm{cen}}(G^\dagger)].
\]
It intertwines the operators defined above by
\[
\mathcal F T_b=M_b^\dagger\mathcal F,\qquad \mathcal F M_\chi=T_{\chi^{-1}}^\dagger\mathcal F.
\]
Moreover,
\[
d_m(G)=d_m(G^\dagger).
\]
The construction is equivariant under based root-datum automorphisms and permutations of isomorphic simple factors.
\end{theorem}

\begin{proof}
Recall the notation $Q^\vee\subset X_*(T)\subset P^\vee$ for the coroot, cocharacter, and coweight lattices of $G$, and $Q\subset X^*(T)\subset P$ for its root, character, and weight lattices. Set $E=Q^\vee\otimes_{\ZZ}\RR$ and $E^*=Q\otimes_{\ZZ}\RR$. Since $G$ is simply connected, $X_*(T)=Q^\vee$ and $X^*(T)=P$, so that $B=Z(G)=P^\vee/Q^\vee$ and $B^\dagger=Z(G^\dagger)=P/Q$. Because $(Q^\vee)^*=P$ and $(P^\vee)^*=Q$, the finite abelian groups
\[
A_m=\left(\frac{1}{m}P^\vee\right)/Q^\vee,\qquad A_m^\dagger=\left(\frac{1}{m}P\right)/Q
\]
are in perfect duality under $(x,y)\mapsto\eexp{m\pair{y}{x}}$, and the $W$-orbits in $A_m$ parametrize the conjugacy classes of elements $t\in G$ satisfying $t^m\in B$.

Choose a $W$-invariant positive definite real inner product $(\, ,\,)$ on $E$, invariant also under based automorphisms and permutations of isomorphic factors, and equip $E^*$ with its dual metric. The real metric isomorphism $R\colon E\to E^*$, $v\mapsto(v,\,\cdot\,)$, induces
\[
J=\Sym\bigl((R^{-1})_{\CC}\bigr)\colon
\Sym(E_{\CC}^*)\xrightarrow{\sim}\Sym(E_{\CC}).
\]
On each polynomial algebra use the Fischer Hermitian inner product, given in orthonormal coordinates by $\langle X^a,X^b\rangle=\delta_{ab}\prod_i a_i!$. The map $R^{-1}$ takes a dual orthonormal basis to the corresponding orthonormal basis of $E$, so $J$ is a complex-linear, degree-preserving isometry and commutes with $W$. Consider the graded space
\[
\mathscr A_G=\{f\colon A_m\to\Sym(E_{\CC}^*): f(wx)=\det(w)\,w f(x)\text{ for all }w\in W\},
\]
and define $\mathscr A_{G^\dagger}$ with $A_m^\dagger$ and $\Sym(E_{\CC})$ in place of $A_m$ and $\Sym(E_{\CC}^*)$. With the sum of the Fischer inner products over the finite group, the transform
\[
(\mathbf F_m f)(y)=|A_m|^{-1/2}\sum_{x\in A_m}\eexp{m\pair{y}{x}}J(f(x))
\]
is unitary by character orthogonality and the isometry of $J$. The lattice pairing is $W$-invariant and $J$ is equivariant, so the transform and its inverse preserve the alternating condition, yielding a degree-preserving isomorphism $\mathscr A_G\xrightarrow{\sim}\mathscr A_{G^\dagger}$.

We next identify the lowest nonzero homogeneous part of $\mathscr A_G$. For each $W$-orbit in $A_m$, take its unique point $x$ in the closed fundamental alcove (and the product of the fundamental alcoves for the semisimple case), and set
\[
\Phi_x=\{\alpha\in\Phi:\alpha(x)\in\ZZ\},\qquad N_x=\frac{|\Phi_x|}{2}.
\]
By the description of point stabilizers in the affine Weyl group $Q^\vee\rtimes W$ as the groups generated by reflections in the walls through the point, the stabilizer of the class of $x$ modulo $Q^\vee$ is $W_x=W(\Phi_x)$. The inward gradients of the alcove walls containing $x$ determine a positive system $\Phi_x^+$. Define
\[
\Delta_x(v)=\prod_{\alpha\in\Phi_x^+}\alpha(v)\in \Sym(E_{\CC}^*).
\]
On the orbit of $x$, the defining condition for $\mathscr A_G$ says that $f(x)$ is $W_x$-alternating. Every such polynomial vanishes on each reflection hyperplane and is therefore divisible by $\Delta_x$, with $W_x$-invariant quotient. Thus
\[
\Sym(E_{\CC}^*)^{W_x,\det}=\Delta_x\,\Sym(E_{\CC}^*)^{W_x}.
\]
Consequently, this orbit contributes first in degree $N_x$, and its contribution in that degree is one-dimensional.

For $t_x=\eexp{x}$, as in the proof of Lemma~\ref{lem:odd-local},
\[
\dim_{\RR}\fg^{\Ad(t_x)}=r+|\Phi_x|=r+2N_x,
\]
where $r=\operatorname{rank}G$. Hence the lowest nonzero degree of $\mathscr A_G$ is
\[
N_m=\min_{x\in A_m}N_x,\qquad d_m(G)=r+2N_m.
\]
Since $\mathbf F_m$ preserves degrees, the corresponding lowest degree on the dual side is also $N_m$; as the two groups have the same rank, $d_m(G)=d_m(G^\dagger)$. For each minimal orbit, define a unit vector $v_x$, supported on $W\cdot x$, by
\[
v_x(wx)=\frac{\det(w)\,w\Delta_x}{\sqrt{|W:W_x|}\,\|\Delta_x\|}.
\]
These vectors form an orthonormal basis of $(\mathscr A_G)_{N_m}$. The natural bijection $\cM_m^{\mathrm{cen}}(G)\leftrightarrow \{W\cdot x\in A_m/W:N_x=N_m\}$ identifies $(\mathscr A_G)_{N_m}$ with $\CC[\cM_m^{\mathrm{cen}}(G)]$. Using the analogous identification on the dual side, the restriction of $\mathbf F_m$ defines the transformation $\mathcal F\colon\CC[\cM_m^{\mathrm{cen}}(G)]\xrightarrow{\sim}\CC[\cM_m^{\mathrm{cen}}(G^\dagger)]$.

We now determine the central operators. For $b\in B=P^\vee/Q^\vee$, let $a_b(v)=w_bv+\mu_b$ ($w_b\in W$, $\mu_b\in P^\vee$) be the unique element $(\mu_b,w_b)$ in $\widetilde{W}=P^\vee\rtimes W$ that preserves the fundamental alcove and represents $b$, and define translation on $\mathscr A_G$ by $(\tau_{\mu_b}f)(x)=f(x-\mu_b)$. Since $a_b$ preserves the inward wall gradients, $w_b\Delta_x=\Delta_{a_b(x)}$; it also preserves the stabilizer order and the norm of the discriminant, so $\tau_{\mu_b}v_x=\det(w_b)\,v_{a_b(x)}$. Under the class-basis identification, this reads $T_b=\varepsilon_b\tau_{\mu_b}$ with $\varepsilon_b=\det(w_b)$. To compute this sign, consider the Weyl denominator
\[
D(v)=\prod_{\alpha\in\Phi^+}2\sin(\pi\alpha(v)),
\]
which is positive in the interior of the fundamental alcove and satisfies
\[
D(wv)=\det(w)D(v),\qquad D(v+\mu)=(-1)^{\pair{2\varrho}{\mu}}D(v)\quad(\mu\in P^\vee),
\]
where $\varrho=\tfrac12\sum_{\alpha\in\Phi^+}\alpha$ is the half-sum of the positive roots. Since $a_b$ preserves the alcove, we obtain
\[
\varepsilon_b=(-1)^{\pair{2\varrho}{\mu_b}}=z_{G^\dagger}(b),
\]
where $z_{G^\dagger}=(2\varrho)(-1)$ is viewed as a character of $B$ through $B^\dagger=\widehat B$.

Let $\chi\in B^\dagger=P/Q$ be represented by $\lambda\in P$, and define $(V_\lambda f)(x)=\eexp{m\pair{\lambda}{x}}f(x)$ and, on the dual side, $(V_\mu^\dagger h)(y)=\eexp{m\pair{y}{\mu}}h(y)$ for $\mu\in P^\vee$. Direct calculation gives
\[
\mathbf F_m\tau_\mu=V_\mu^\dagger\mathbf F_m,\qquad\mathbf F_m V_\lambda=\tau_{-\lambda}^\dagger\mathbf F_m.
\]
In the normalized class bases, the shifted operators satisfy $M_\chi=\chi(z_G)^{-1}V_\lambda$ and $M_b^\dagger=z_{G^\dagger}(b)V_{\mu_b}^\dagger$, while the dual version of the translation calculation gives $\tau_{-\lambda}^\dagger=\chi(z_G)T_{\chi^{-1}}^\dagger$. Combining these identities yields $\mathcal F T_b=M_b^\dagger\mathcal F$ and $\mathcal F M_\chi=T_{\chi^{-1}}^\dagger\mathcal F$.

It remains to verify canonicity and equivariance. The construction depends on two auxiliary choices: the $W$-invariant inner product and the positive root system. On a simple factor, rescaling the inner product by $c>0$ multiplies normalized vectors of degree $N$ in $\mathscr A_G$ by $c^{N/2}$, the map $J$ by $c^{-N}$, and normalized vectors in $\mathscr A_{G^\dagger}$ by $c^{-N/2}$; these factors cancel in every matrix coefficient and hence leave $\mathcal F$ unchanged. A simultaneous Weyl group action on the two positive systems changes the basis identifications on both sides by the same determinant sign, which also cancels. Finally, based root-datum automorphisms and permutations of isomorphic factors preserve the pairing, the alcove construction, and the normalized local discriminants, so $\mathcal F$ is canonical and equivariant as claimed.
\end{proof}

\section{The binary tetrahedral group \texorpdfstring{$2T$}{2T}}\label{sec:tetrahedral}

In this section we take $\Gamma=2T$, with spherical type $(2,3,3)$ and abelianization $C_3$. Building on the ordinary elliptic comparison of Section~\ref{sec:odd-models}, we incorporate the actual carry grades and central translation operators required by Theorem~\ref{thm:transfer}, and deduce elliptic lifting under central isogenies.

Throughout Sections~\ref{sec:tetrahedral}--\ref{sec:full-reciprocity}, we recall the local minimal-class notation
\begin{equation}\label{eq:local-minimal-classes}
\begin{gathered}
C_{h,m}(G)=\left\{[t]_G:t^m=h,\ \dim_{\RR}\fg^{\Ad(t)}=d_m(G)\right\},\\
M_m(G)=\CC[\cM_m(G)],\qquad \cM_m(G)=C_{1,m}(G).
\end{gathered}
\end{equation}
The class spaces carry their orthonormal class bases, and we use the positive Fourier convention for global projective modules,
\begin{equation}\label{eq:global-fourier-convention}
\mathcal F T_k=M_k^\dagger\mathcal F,
\qquad
\mathcal F M_\chi=T_{\chi^{-1}}^\dagger\mathcal F.
\end{equation}

\subsection{Zero-packet support}

\begin{lemma}\label{lem:tetrahedral-zero-packet}
Let $G$ be simply connected compact simple and put $B=Z(G)$. An elliptic projective class of $2T$ can occur only at an actual grade $b\in B^3$, and every twist by $B[3]$ fixes each such conjugacy class. Thus the elliptic multiplicities of Definition~\ref{def:packets} are supported only at $(e,\omega)=(0,0)$. If
\[
\delta_3(G)=d_2(G)+2d_3(G)-\dim\fg
\]
is positive, all actual-grade elliptic spaces vanish. If $\delta_3(G)=0$, then $G$ has type $A_1$, $A_2$, or $G_2$.
\end{lemma}
\begin{proof}
Proposition~\ref{prop:defect-screen} gives the vanishing and the three zero-defect types. In types $A_1$ and $G_2$ the center has order prime to $3$, so $B^3=B$ and $B[3]=1$.

In type $A_2$, write $G=\SU(3)$ and identify $B$ with the scalar matrices $\mu_3$. Let $f$ be elliptic of actual grade $b\in\mu_3$. The centrality argument for carry cocycles in Proposition~\ref{prop:defect-screen} gives $f(-1)\in B[2]=\{I\}$. Set $Z=f(z)$. Since $\ell(z)=1$, $\ell(-1)=0$, and $z^3=-1$, multiplication with the prescribed carry gives
\begin{equation}\label{eq:tetrahedral-su3-relations}
Z^3=bI.
\end{equation}
The adjoint action of $f$ factors through $A_4$, and Lemma~\ref{lem:minima} gives $\dim\fg^{\Ad Z}=d_3(G)=2$. Thus $Z$ is regular and its eigenvalues are the three distinct roots of $t^3-b$. Their product is $b$, so $\det Z=1$ forces $b=1$.

At grade $1$, every regular cubic class has spectrum $\{1,\omega,\omega^2\}$, where $\omega=e^{2\pi i/3}$; this spectrum gives an element of $\SU(3)$. Hence $\cM_3(G)$ is a singleton, and $H_G=\{I\}$ since $B[2]=\{I\}$ and $z_G\in B[2]$. The ordinary model~\eqref{eq:spherical-model} identifies $E_{2T}(G)$ with $H_G\times\cM_3(G)^2$, whose unique class is fixed by every twist in $B[3]$.
\end{proof}

\subsection{The full-center Fourier transform}

\begin{theorem}\label{thm:tetrahedral-fourier}
Let $G$ be simply connected compact simple, put $B=Z(G)$, and let $G^\dagger=(G/B)^\vee$ with $B^\dagger=Z(G^\dagger)=\widehat B$. There is a canonical unitary Fourier isomorphism
\[
\wtE_{2T,B}(G)\xrightarrow{\sim}\wtE_{2T,B^\dagger}(G^\dagger)
\]
satisfying~\eqref{eq:global-fourier-convention}. It is equivariant under based root-datum automorphisms, and its tensor products commute with permutations of isomorphic simple factors.
\end{theorem}
\begin{proof}
Theorem~\ref{thm:minimum-fourier} gives the same defect on the two sides. At positive defect both spaces vanish by Lemma~\ref{lem:tetrahedral-zero-packet}. At zero defect the possible centers are $C_2,C_3,$ and $1$, each of which has $B^2\cap B[2]=\{1\}$, so $H_G=\{z_G\}$. The ordinary model~\eqref{eq:spherical-model} gives $E_{2T}(G)\simeq\cM_3(G)^2$ and the corresponding identification for $G^\dagger$.

Lemma~\ref{lem:tetrahedral-zero-packet} applies to both $(G,B)$ and $(G^\dagger,B^\dagger)$. Applying Lemma~\ref{lem:zero-packet} to each pair gives the decomposition
\begin{equation}\label{eq:tetrahedral-zero-packet}
\wtE_{2T,B}(G)
\simeq\CC[B^3]\otimes\CC[E_{2T}(G)]
\simeq\CC[B^3]\otimes M_3(G)^{\otimes2}
\end{equation}
and its dual analogue. Use the positive Fourier transform for the perfect pairing
\[
B^3\times(B^\dagger)^3\longrightarrow\UU(1),
\qquad (k^3,\chi^3)\longmapsto\chi(k)^3,
\]
from Lemma~\ref{lem:zero-packet}, and tensor it with the two linear extensions of the cubic class bijection in Lemma~\ref{lem:odd-local}. The central operators act on the first factor of~\eqref{eq:tetrahedral-zero-packet}, so the stated intertwining relations follow from Lemma~\ref{lem:zero-packet}. The pairing and the odd local bijections are equivariant under based automorphisms, and their tensor products give the asserted compatibility with factor permutations.
\end{proof}

\subsection{Elliptic lifting under central isogenies}

\begin{corollary}\label{cor:tetrahedral-lifting}
Every central isogeny $\pi\colon\widetilde G\to G$ of connected compact semisimple groups induces a canonical bijection
\[
\pi_*\colon E_{2T}(\widetilde G)\xrightarrow{\sim}E_{2T}(G),
\]
equivariant under automorphisms of the isogeny. For every connected compact semisimple $P$ and finite $B\subset Z(P)$, the tetrahedral elliptic packet is consequently supported only at $(0,0)$.
\end{corollary}
\begin{proof}
First write $U=\prod_iG_i$ with each $G_i$ simply connected compact simple, put $D=Z(U)$, and take $A\subset D$. If a factor has positive tetrahedral defect, its projective elliptic spaces vanish, so factorwise ellipticity and Lemma~\ref{lem:residual-descent} imply $E_{2T}(U/A)=\varnothing$.

Otherwise all factors have type $A_1$, $A_2$, or $G_2$. Hence $D$ has exponent dividing $6$, and
\[
D^3=D_{(2)},\qquad A\cap D^3=A^3.
\]
Indeed, cubing is the identity on the $2$-primary subgroup of $D$ and kills its $3$-primary subgroup. In particular, if $a\in A\cap D^3$, then $a^2=1$ and $a=a^3\in A^3$. The reverse inclusion follows from $A\subset D$. By Lemma~\ref{lem:tetrahedral-zero-packet}, only grades in $D^3$ occur and $D[3]$ fixes every ordinary elliptic class. The normalized transports of Lemma~\ref{lem:zero-packet} thus identify the set of all actual-grade elliptic classes with
\[
D^3\times E_{2T}(U),
\qquad
T_a(b,[\rho])=(a^3b,[\rho]).
\]
The class-level descent in Lemma~\ref{lem:residual-descent}, with $C=A$, identifies $E_{2T}(U/A)$ with the $T_A$-orbits whose grades lie in $A$. The allowed grades are $A^3$, on which $T_A$ is transitive, and the ordinary class coordinate is unchanged. Hence projection induces a bijection $E_{2T}(U)\to E_{2T}(U/A)$.

Every central isogeny of connected compact semisimple groups has the form $U/A\to U/C$ for a common simply connected cover and $A\subset C\subset Z(U)$. Applying the preceding bijection twice proves the first assertion and its equivariance.

Let $q\colon P\to P/B$ be the quotient. An elliptic projective map $f$ of grade $b$ projects to an ordinary elliptic homomorphism, since $q$ identifies the adjoint Lie algebras. Surjectivity of $q_*$ gives an ordinary elliptic $\rho\colon 2T\to P$ whose projection is conjugate to $q\circ f$. Lift a conjugating element to $P$ and conjugate $\rho$ so that $q\circ\rho=q\circ f$. Then $f(g)=\lambda(g)\rho(g)$ for a $B$-valued cochain $\lambda$, and comparing multiplication laws gives $c_b(g,g')=\lambda(g)\lambda(g')\lambda(gg')^{-1}$; hence $[c_b]=0$, and~\eqref{eq:polyhedral-cohomology} gives $b\in B^3$.

If $\eta\in B[3]$, the ordinary elliptic homomorphisms $\rho$ and $\eta^\ell\rho$ have identical projections. Injectivity of $q_*$ gives $[\eta^\ell\rho]=[\rho]$. For any $b=k^3$, transport $T_k$ identifies the ordinary layer with the layer at $b$ and commutes with $T_\eta$. Thus all twists by $B[3]$ fix every class in every supported layer, so the elliptic packet is supported only at $(0,0)$.
\end{proof}

Corollary~\ref{cor:tetrahedral-lifting} gives uniqueness up to conjugacy. For a fixed elliptic homomorphism $\bar\rho\colon 2T\to G$, a conjugating element downstairs lifts to $\widetilde G$, so the bijection gives a homomorphism lift with projection $\bar\rho$. Its homomorphism lifts form a torsor under $\Hom(2T,\ker\pi)=(\ker\pi)[3]$, acting by pointwise multiplication.

\section{The binary octahedral group \texorpdfstring{$2O$}{2O}}\label{sec:octahedral}

In this section $\Gamma=2O$, with spherical type $(2,3,4)$ and abelianization $C_2$. We first parametrize the elliptic projective classes locally and then treat the center cases from Proposition~\ref{prop:defect-screen} to establish the full-center elliptic hypothesis of Theorem~\ref{thm:transfer} for $2O$.

\subsection{Local realization}

\begin{lemma}\label{lem:octahedral-local-model}
Let $G$ be simply connected compact semisimple, put $B=Z(G)$, and assume $d_2+d_3+d_4=\dim\fg$. The class basis of $\wtE_{2O,B}(G)$ is canonically parametrized by
\begin{equation}\label{eq:octahedral-local-model}
\coprod_{b\in B}\ \coprod_{h\in B[2]}
C_{bh,2}(G)\times C_{h,3}(G)\times C_{b^2h,4}(G).
\end{equation}
The first index is the actual carry grade, and $h$ is the image of $-1\in2O$.
\end{lemma}
\begin{proof}
Fix $b\in B$ of order $d$ and let $\epsilon(g)\in\{0,1\}$ represent the parity of $g\in2O$. The fiber product for the parity map and reduction $C_{2d}\to C_2$ is
\[
\Gamma_b=2O\times_{C_2}C_{2d}
=\{(g,\bar n):n\equiv\epsilon(g)\pmod 2\}.
\]
The projection to $2O$ has cyclic kernel generated by $\beta=(1,\bar 2)$ of order $d$. Its section $s(g)=(g,\overline{\epsilon(g)})$ satisfies $s(g)s(g')=\beta^{\epsilon(g)\epsilon(g')}s(gg')$. Therefore evaluation $f(g)=\rho(s(g))$ takes each homomorphism $\rho\colon\Gamma_b\to G$ with $\rho(\beta)=b$ to a projective map of grade $b$. Every element of $\Gamma_b$ is uniquely $\beta^j s(g)$ with $j\in\ZZ/d\ZZ$, and the inverse assignment
\[
\rho(\beta^j s(g))=b^j f(g)
\]
is multiplicative because the carry in $f(g)f(g')$ is $b^{\epsilon(g)\epsilon(g')}$. These assignments commute with conjugation and give the same centralizer, since $\rho(\Gamma_b)$ is generated by $b$ and $f(2O)$.

The action on the sphere factors through $2O/\{\pm1\}\simeq S_4$ and has central kernel $K=\langle\beta\rangle\times\langle h_0\rangle$, where $h_0=s(-1)=(-1,\bar0)$. For the lifts $\widetilde x=s(x)$, $\widetilde y=s(y)$, and $\widetilde z=s(z)$ we have
\[
\widetilde x^2=\beta h_0,\qquad \widetilde y^3=h_0,\qquad
\widetilde z^4=\beta^2h_0,\qquad \widetilde x\widetilde y\widetilde z=\beta h_0.
\]
For elliptic $f$, apply Proposition~\ref{prop:central-involution} in each simple factor to the images $X=f(x)$, $Y=f(y)$, $Z=f(z)$, and $h=f(-1)=\rho(h_0)$. The displayed relations give its power and product hypotheses, and $h_0$ is a central involution, so $h\in B[2]$ and the adjoint action factors through $S_4$. By Lemma~\ref{lem:minima}, the local images are minimal at all three exceptional stabilizers, so $([X],[Y],[Z])$ lies in the $(b,h)$-summand.

Conversely, fix $h\in B[2]$ and a triple in that summand. Define $\chi\colon K\to B$ by $\chi(\beta)=b$, $\chi(h_0)=h$. The three stabilizers are generated by $K$ and respectively $\widetilde x,\widetilde y,\widetilde z$; their quotients by $K$ are cyclic of orders $2,3,4$. The displayed power relations show that sending these lifts to representatives $X,Y,Z$ gives local homomorphisms extending $\chi$. Their fixed dimensions sum to $d_2+d_3+d_4=\dim\fg$ because $\chi(K)$ is central and each local class is minimal. Theorem~\ref{thm:local-realization} gives a unique elliptic conjugacy class of homomorphisms $\rho\colon\Gamma_b\to G$ with the prescribed local classes. Restriction along $s$ gives the required projective class of grade $b$ with $f(-1)=h$, and taking the disjoint union over $b$ and $h$ gives~\eqref{eq:octahedral-local-model}.
\end{proof}

\subsection{Case analysis}

By Proposition~\ref{prop:defect-screen}, the simple types with $\delta_4=0$ are $A_r$ ($1\le r\le 3$), $B_r$ and $C_r$ ($2\le r\le 4$), $D_4$, $G_2$, and $F_4$; we treat these according to the structure of $B=Z(G)$.

\begin{lemma}\label{lem:octahedral-exponent-two}
Under the hypotheses of Lemma~\ref{lem:octahedral-local-model}, suppose that $B^2=1$. There is a canonical unitary identification
\begin{equation}\label{eq:octahedral-exponent-two}
\wtE_{2O,B}(G)
\simeq M_3(G)\otimes\CC[\cM_4^{\mathrm{cen}}(G)],
\end{equation}
under which the actual carry grade is the shifted quartic grade and both central operators act on the second factor. Put $G^\dagger=(G/B)^\vee$ and $B^\dagger=Z(G^\dagger)\simeq\widehat B$. The cubic class bijection and the quartic Fourier transform induce a canonical unitary isomorphism
\[
\wtE_{2O,B}(G)\xrightarrow{\sim}\wtE_{2O,B^\dagger}(G^\dagger)
\]
satisfying~\eqref{eq:global-fourier-convention} and equivariant under based root-datum automorphisms.
\end{lemma}
\begin{proof}
Corollary~\ref{cor:squares} gives a unique minimal quadratic class, with square $z_G$. The local equations in Lemma~\ref{lem:octahedral-local-model} consequently give
\[
bh=z_G,\qquad Y^3=h,\qquad Z^4=h.
\]
Thus $b=Z^4z_G^{-1}$ and $(hY)^3=1$. Conversely, $[Y']\in\cM_3(G)$ and $[Z]\in\cM_4^{\mathrm{cen}}(G)$ determine
\[
h=Z^4,\qquad b=hz_G^{-1},\qquad [Y]=[hY'],
\]
and the unique quadratic class completes the local triple. This proves the identification by $[f]\mapsto([hY],[Z])$.

A central twist by $k\in B$ sends $(X,Y,Z)$ to $(kX,Y,kZ)$ and leaves $h$ unchanged; it fixes the unique quadratic class and therefore acts as $1\otimes T_k$. Since $b=Z^4z_G^{-1}$, multiplication by a grade character acts as $1\otimes M_\chi$.

Theorem~\ref{thm:minimum-fourier} identifies the three minimal dimensions on the two sides; since $\dim G^\dagger=\dim G$, the dual defect is also zero. Since $B^\dagger\simeq\widehat B$ has exponent dividing two, the decomposition~\eqref{eq:octahedral-exponent-two} applies on the dual side as well. In these decompositions tensor the linear extension of the cubic bijection of Lemma~\ref{lem:odd-local} with the quartic Fourier transform of Theorem~\ref{thm:minimum-fourier}. Each factor is unitary and equivariant under based automorphisms. The quartic intertwining identities match~\eqref{eq:global-fourier-convention} directly.
\end{proof}

\begin{lemma}\label{lem:octahedral-su3}
For $G=\SU(3)$ and $B=\mu_3$, there is a unique elliptic projective class of $2O$ in every actual grade. The resulting full-center space is canonically $\CC[B]$, with its positive Fourier transform for the pairing of Lemma~\ref{lem:zero-packet}.
\end{lemma}
\begin{proof}
Here $(d_2,d_3,d_4)=(4,2,2)$, their sum is $\dim G=8$, and $B[2]=1$. At ordinary grade $b=1$, the required minimal local classes have unique spectra
\[
\{1,-1,-1\},\qquad
\{1,\omega,\omega^2\},\qquad
\{1,i,-i\},\qquad \omega=e^{2\pi i/3}.
\]
The quadratic multiplicities are $(1,2)$ and determinant one forces the signs shown. The cubic and quartic classes are regular, and for the latter, determinant one forces the omitted fourth root to be $-1$. Lemma~\ref{lem:octahedral-local-model} therefore gives one ordinary elliptic class.

The centers of $G$ and its simply connected dual-root group $G^\dagger$ of type $A_2$ both have order three. In each center squaring is bijective, so Remark~\ref{rem:coprime-packet} and Lemma~\ref{lem:zero-packet} identify the respective class spaces with $\CC[B]$ and $\CC[B^\dagger]$, with one class per actual grade. Choose $a=\omega I$ and the character $\chi$ satisfying $\chi(a)=\omega$. The unique square root of $a^r$ is $a^{2r}$, so the pairing~\eqref{eq:zero-pairing} takes the value $\chi^s(a^{2r})=\omega^{2rs}$ on $(a^r,\chi^s)$. In the class bases indexed by these grades the normalized positive transform is
\begin{equation}\label{eq:octahedral-su3-fourier}
\mathcal F(e_r)=\frac1{\sqrt3}\sum_{s=0}^2\omega^{2rs}e_s^\dagger.
\end{equation}
Lemma~\ref{lem:zero-packet} gives the intertwining identities~\eqref{eq:global-fourier-convention}, and naturality of the pairing yields equivariance under based automorphisms.
\end{proof}

\begin{lemma}\label{lem:octahedral-su4}
Let $G=\SU(4)$, let $B=\langle\zeta\rangle$ with $\zeta=iI$, and choose $\chi(\zeta)=i$. There is one elliptic projective class of $2O$ in each actual grade $I,-I$, and none in grades $\pm iI$. Both supported classes have $f(-1)=-I$. Order the class basis $(e_+,e_-)$ by the grades $I,\zeta^2$. Then
\begin{equation}\label{eq:octahedral-su4-operators}
T_\zeta=\begin{pmatrix}0&1\\1&0\end{pmatrix},
\qquad
M_\chi=\begin{pmatrix}1&0\\0&-1\end{pmatrix}.
\end{equation}
With the dual class basis ordered by $1,\chi^2$, the canonical positive Fourier transform is
\begin{equation}\label{eq:octahedral-su4-fourier}
\mathcal F=\frac1{\sqrt2}\begin{pmatrix}1&1\\1&-1\end{pmatrix}.
\end{equation}
It satisfies~\eqref{eq:global-fourier-convention} and is equivariant under based root-datum automorphisms.
\end{lemma}
\begin{proof}
For $\SU(4)$, we have $(d_2,d_3,d_4)=(7,5,3)$ and $7+5+3=15=\dim G$. Lemma~\ref{lem:octahedral-local-model} applies, with
\[
h\in\{I,-I\},\qquad X^2=bh,\qquad Y^3=h,\qquad Z^4=b^2h.
\]
Here $B^2=B[2]=\{I,-I\}$ contains $z_G$. By Corollary~\ref{cor:squares}, the possible squares of minimal quadratic classes are $z_G B^2=B^2$, with a unique class for each square. Since $X^2=bh$ and $h\in B^2$, we must have $b\in B^2$, so $b^2=I$ and $Z^4=h$. The quartic class is regular, so for $h=cI$ its four eigenvalues are all the roots of $t^4-c$, whose product is $-c=1$; this forces $h=-I$ and determines a unique quartic class.

A minimal cubic class with $Y^3=-I$ has multiplicities $(2,1,1)$. If $\lambda$ is the repeated eigenvalue, the product of the three distinct roots is $-1$, so $\det Y=-\lambda=1$. Hence its unique spectrum is
\[
\{-1,-1,e^{\pi i/3},e^{-\pi i/3}\}.
\]
For each $b=I,-I$, the prescribed square $X^2=-b$ determines one quadratic class, so Lemma~\ref{lem:octahedral-local-model} realizes one global class in each of these grades.

Translation by $\zeta$ exchanges the two grades, and $M_\chi$ takes the values $1,-1$ on them, giving~\eqref{eq:octahedral-su4-operators}. Each element of $B[2]$ preserves both actual grades and hence fixes their unique classes, and the support is $B^2$, so the zero-packet condition holds. The same computation applies to $G^\dagger$ of simply connected type $A_3$, so Lemma~\ref{lem:zero-packet} identifies the two coefficient spaces with $\CC[B^2]$ and $\CC[(B^\dagger)^2]$. In the stated bases the only pairing value different from one is
\[
\pair{\zeta^2}{\chi^2}=\chi(\zeta)^2=-1.
\]
Formula~\eqref{eq:zero-pairing} now gives the normalized matrix~\eqref{eq:octahedral-su4-fourier}, and Lemma~\ref{lem:zero-packet} gives its operator identities and equivariance.
\end{proof}

\subsection{The full-center Fourier transform}

\begin{theorem}\label{thm:octahedral-fourier}
For every simply connected compact simple group $G$, with $B=Z(G)$, $G^\dagger=(G/B)^\vee$, and $B^\dagger=Z(G^\dagger)=\widehat B$, there is a canonical unitary Fourier isomorphism
\[
\wtE_{2O,B}(G)\xrightarrow{\sim}\wtE_{2O,B^\dagger}(G^\dagger)
\]
satisfying~\eqref{eq:global-fourier-convention}. It is equivariant under based root-datum automorphisms and compatible with permutations of isomorphic simple factors.
\end{theorem}
\begin{proof}
The defect $\delta_4(G)=d_2+d_3+d_4-\dim\fg$ is nonnegative by Lemma~\ref{lem:minima}, and is unchanged by root duality by Theorem~\ref{thm:minimum-fourier}. At positive defect both spaces vanish by Proposition~\ref{prop:defect-screen}. Among its zero-defect simple types, only $A_2$ and $A_3$ have centers of exponent greater than two, namely $C_3$ and $C_4$. The other centers are trivial, $C_2$, or $C_2\times C_2$ (in type $D_4$), so Lemma~\ref{lem:octahedral-exponent-two} applies. Lemmas~\ref{lem:octahedral-su3} and~\ref{lem:octahedral-su4} give the transforms for $A_2$ and $A_3$. Each transform satisfies~\eqref{eq:global-fourier-convention} and is equivariant under based automorphisms, and their tensor products commute with permutations of isomorphic factors.
\end{proof}

\section{The binary icosahedral group \texorpdfstring{$2I$}{2I}}\label{sec:full-reciprocity}

In this section $\Gamma=2I$, with spherical type $(2,3,5)$ and trivial abelianization. Because the carry cocycle is trivial in every actual grade, the full-center elliptic Fourier transform follows directly from the ordinary model of Section~\ref{sec:odd-models}.

\begin{theorem}\label{thm:icosahedral-fourier}
Let $G$ be simply connected compact simple, put $B=Z(G)$, and let $G^\dagger=(G/B)^\vee$ with $B^\dagger=Z(G^\dagger)=\widehat B$. There is a canonical unitary Fourier isomorphism
\[
\wtE_{2I,B}(G)\xrightarrow{\sim}\wtE_{2I,B^\dagger}(G^\dagger)
\]
satisfying~\eqref{eq:global-fourier-convention}. It is equivariant under based root-datum automorphisms and compatible with permutations of isomorphic simple factors.
\end{theorem}
\begin{proof}
Put $\delta_5(G)=d_2+d_3+d_5-\dim\fg$. By Proposition~\ref{prop:central-involution} and Lemma~\ref{lem:minima}, the existence of an elliptic class forces $\delta_5(G)=0$, and Theorem~\ref{thm:minimum-fourier} gives $\delta_5(G^\dagger)=\delta_5(G)$; thus, if $\delta_5(G)$ is nonzero, both elliptic spaces vanish.

Suppose that $\delta_5(G)=0$. For an ordinary elliptic homomorphism $\rho\colon 2I\to G$, put $h=\rho(-1)$. Its local powers satisfy
\[
X^2=Y^3=Z^5=h,\qquad h\in B[2].
\]
Corollary~\ref{cor:squares} says that a minimal quadratic class with square $h$ exists if and only if
\[
h\in H_G=z_GB^2\cap B[2],
\]
and is then unique. Since $h^2=1$, multiplication by $h$ identifies $C_{h,m}(G)$ with $\cM_m(G)$ for $m=3,5$. Conversely, these local classes meet the dimension hypothesis of Theorem~\ref{thm:local-realization}, applied to $2I\to A_5$ with $-1\mapsto h$. We therefore obtain the parametrization
\begin{equation}\label{eq:icosahedral-ordinary-model}
E_{2I}(G)\xrightarrow{\sim}H_G\times\cM_3(G)\times\cM_5(G),\qquad [\rho]\longmapsto\bigl(h,[h\rho(y)],[h\rho(z)]\bigr).
\end{equation}
Proposition~\ref{prop:icosahedral-center} and Lemma~\ref{lem:odd-local} identify the three factors with their dual counterparts, equivariantly under based automorphisms.

Here $\nu=1$, so the chosen carry cocycle is trivial in every actual grade and the central twist changes only that grade. Hence, equivalently by Lemma~\ref{lem:zero-packet},
\begin{equation}\label{eq:icosahedral-full-center}
\wtE_{2I,B}(G)\simeq\CC[B]\otimes\CC[E_{2I}(G)]\simeq\CC[B]\otimes\CC[H_G]\otimes M_3(G)\otimes M_5(G).
\end{equation}
On the first factor use the ordinary positive Fourier transform
\[
e_b\longmapsto |B|^{-1/2}\sum_{\chi\in B^\dagger}\chi(b)e_\chi^\dagger.
\]
Tensoring it with the linear extension of the ordinary class bijection proves the assertion, while the central pairing and the three factor bijections give the stated equivariance and compatibility with products.
\end{proof}

\section{Conclusion}\label{sec:polyhedral-conclusion}

We now prove Theorem~\ref{thm:polyhedral-main}, which completes Theorem~\ref{thm:unified}\textup{(1)} and~\textup{(2)} for $\Gamma\in\{2T,2O,2I\}$.

\begin{proof}[Proof of Theorem~\ref{thm:polyhedral-main}]
For simply connected compact simple $S$, the required full-center elliptic Fourier isomorphisms
\[
\wtE_{\Gamma,Z(S)}(S)\xrightarrow{\sim}\wtE_{\Gamma,Z(S^\dagger)}(S^\dagger)
\]
are given by Theorem~\ref{thm:tetrahedral-fourier} for $\Gamma=2T$, Theorem~\ref{thm:octahedral-fourier} for $\Gamma=2O$, and Theorem~\ref{thm:icosahedral-fourier} for $\Gamma=2I$. Each of these isomorphisms satisfies \eqref{eq:global-fourier-convention}, is equivariant under based root-datum automorphisms, and is compatible with permutations of isomorphic simple factors. Theorem~\ref{thm:transfer} therefore gives \eqref{eq:final-graded-fourier} for arbitrary connected compact reductive $G$ and finite $B\subset Z(G)$, together with
\[
\mathcal F T_k=M_k^\dagger\mathcal F,\qquad \mathcal F M_\chi=T_{\chi^{-1}}^\dagger\mathcal F.
\]
The polynomial factors in this construction are those of Lemma~\ref{lem:levi}; hence they preserve the invariants of the full centralizers. Taking degree zero gives the asserted isomorphism of projective class modules by Definition~\ref{def:projective-modules}.
\end{proof}

For the binary polyhedral groups, let $\nu=|\Gamma_{\mathrm{ab}}|$, so $\nu=3,2,1$ for $2T,2O,2I$, respectively. For a connected compact reductive group $G$ and finite $B\subset Z(G)$, write
\[
n_G(e,\omega;u)\in\ZZ[[u]],\qquad e\in\widehat{B[\nu]},\quad \omega\in B/B^\nu,
\]
for the Hilbert series of the $T_\eta$-eigenspace with character $e$ in the actual grade $b$ summand of $\wtP_{\Gamma,B}(G)$, where $b$ represents $\omega$, that is,
\[
n_G(e,\omega;u)=\sum_{q\ge0}u^q\dim_{\CC}\left\{v\in\bigoplus_{[f]\in\cR_b(G)}\left(\Sym^q(\mathfrak s_f^*)\right)^{W_f}:T_\eta v=e(\eta)v\text{ for all }\eta\in B[\nu]\right\}.
\]
This is the graded version of Definition~\ref{def:packets}, and the transports $T_k$ show that it is independent of the representative. Below we write the packet labels additively.

\begin{proposition}\label{prop:packet-reciprocity}
Fix $\Gamma\in\{2T,2O,2I\}$. Let $G$ be connected compact reductive, let $B\subset Z(G)$ be finite, and put $G^\dagger=(G/B)^\vee$.

Under the canonical evaluation identifications
\[
\widehat{B[\nu]}\simeq B^\dagger/(B^\dagger)^\nu,\qquad B/B^\nu\simeq\widehat{B^\dagger[\nu]}
\]
the packet Hilbert series satisfy
\begin{equation}\label{eq:hilbert-reciprocity}
n_G(e,\omega;u)=n_{G^\dagger}(-\omega,e;u)=n_{G^\dagger}(\omega,e;u).
\end{equation}

For each of these three groups, the corresponding degree-zero cohomological modules satisfy the reciprocity of~\cite[Conjecture~1.9]{kojima-tachikawa} under the Poincar\'e duality swap on $S^3/\Gamma$.
\end{proposition}

\begin{proof}
Denote the actual-grade $b$ summand by
\[
V_b=\bigoplus_{[f]\in\cR_b(G)}\Sym(\mathfrak s_f^*)^{W_f},
\]
and, for $e\in\widehat{B[\nu]}$, set
\[
V_{b,e}=\{v\in V_b:T_\eta v=e(\eta)v\text{ for every }\eta\in B[\nu]\}.
\]
Thus
\[
\operatorname{Hilb}(V_{b,e};u)=n_G(e,\omega;u).
\]
For $\omega\in B/B^\nu$, consider the graded subspace
\[
V_{e,\omega}=\bigoplus_{b'\in bB^\nu}V_{b',e}.
\]
The transports $T_k$ identify the summands in this direct sum and preserve polynomial degree, since $T_kT_\eta=T_\eta T_k$ and $T_k$ sends grade $b'$ to grade $k^\nu b'$; hence
\[
\operatorname{Hilb}(V_{e,\omega};u)=|B^\nu|\,n_G(e,\omega;u).
\]

Define the corresponding spaces on the dual side in the same way, using a superscript $\dagger$. For $\eta\in B[\nu]$, the identity $\mathcal F T_\eta=M_\eta^\dagger\mathcal F$ identifies the character $e$ with the dual grade coset $e\in B^\dagger/(B^\dagger)^\nu$. For $\psi\in B^\dagger[\nu]$, $M_\psi$ has eigenvalue $\psi(b)$ throughout $V_{e,\omega}$, since $\psi$ is trivial on $B^\nu$. The identity
\[
T_\psi^\dagger\mathcal F=\mathcal F M_{\psi^{-1}}
\]
shows that the dual translation character is $\psi\mapsto\psi(b)^{-1}$, corresponding to $-\omega$. It follows that
\[
\mathcal F(V_{e,\omega})=V^\dagger_{-\omega,e}.
\]
Comparing Hilbert series and using $|B^\nu|=|(B^\dagger)^\nu|$ (the $\nu$th-power maps have the same image size on each cyclic factor and its dual) gives
\[
n_G(e,\omega;u)=n_{G^\dagger}(-\omega,e;u).
\]

To remove the minus sign, note that each polynomial summand $\Sym(\mathfrak s_f^*)^{W_f}$ is the complexification of the corresponding invariant polynomial algebra on the real Lie algebra $\Lie(S_f)$. The central twists act through permutations of the classes and real linear maps on these Lie algebras. Thus every fixed actual-grade summand admits a degree-preserving, conjugate-linear involution $\tau$ commuting with $T_\eta$ for $\eta\in B[\nu]$, which exchanges eigenspaces with characters $e$ and $-e$, so
\[
n_G(e,\omega;u)=n_G(-e,\omega;u).
\]
Applying the same argument on the dual side proves the second equality in \eqref{eq:hilbert-reciprocity}.

It remains to compare this equality with the topological Poincar\'e duality pairing. Put $M=S^3/\Gamma$. Since the classifying map $M\to B\Gamma$ has homotopy fiber $S^3$, the low-degree edge maps identify $H^i(\Gamma,C)$ with $H^i(M,C)$ for $i=1,2$ and every finite trivial coefficient group $C$. Moreover, $H_1(M,\ZZ)\simeq C_\nu$, and for $\nu>1$, the nonsingular linking pairing on a generator $a$ has the form
\[
\lambda(a,a)=\frac{q}{\nu}\pmod{\ZZ},\qquad \gcd(q,\nu)=1.
\]
In the coordinates determined by the integral carry class and the abelianization generator, integral Poincar\'e duality $H^2(M,\ZZ)\to H_1(M,\ZZ)$ is multiplication by a unit modulo $\nu$. Since $H^3(M,\ZZ)$ and $H_0(M,\ZZ)$ are free, the coefficient maps identify $H^2(M,C)$ and $H_1(M,C)$ with their integral counterparts tensored with $C$. Naturality of cap product gives the same unit on $C/\nu C$ for every $C$, and the reverse identification uses its inverse. Since $\nu=2$ or $3$, every such unit equals its inverse and is $1$ or $-1$ modulo $\nu$. Hence the Poincar\'e duality swap differs from the coordinate swap by at most a simultaneous sign on the two packet labels. For $\nu=1$, both label groups are trivial. The positive swap is already covered by \eqref{eq:hilbert-reciprocity}; for the simultaneous negative swap, the preceding equalities give
\[
n_G(e,\omega;u)=n_G(-e,\omega;u)=n_{G^\dagger}(-\omega,-e;u).
\]
Thus the packet multiplicities also agree under the Poincar\'e duality swap.

Finally, choose one actual-grade representative for each cohomology class and take polynomial degree zero, as in Remark~\ref{rem:kt-module}. Since these cohomological modules are representations of finite abelian groups, matching all joint character multiplicities implies their equivalence under the Poincar\'e duality swap. This proves \cite[Conjecture~1.9]{kojima-tachikawa}.
\end{proof}

\clearpage
\appendix
\part*{Appendices}

\section{Poincar\'e duality for \texorpdfstring{$S^3/\Dic_m$}{S3/Dic(m)}}\label{sec:poincare}

For cap products, Poincar\'e duality, and transfer under finite coverings, we use~\cite[\S\S3.3 and~3.G]{hatcher}.

We use the generator classes $a,b$ and the cohomology classes $\rho_1,\rho_2$ defined in Section~\ref{sec:packets}. As before, $M=S^3/\Gamma$ denotes the quotient by the left action $\gamma\cdot x=\gamma x$. We mark a based loop by the endpoint of its lift from $1$. Set
\[
P\colon H^2(M,\ZZ)\longrightarrow H_1(M,\ZZ),\qquad
P(\eta)=\eta\cap[M],\qquad
S_0=\begin{pmatrix}1&1\\0&1\end{pmatrix}.
\]
The domain and codomain of $P$ are the two copies of $\coker D$ with the coordinates fixed in Section~\ref{sec:packets}.

\begin{proposition}\label{prop:poincare}
For an overall orientation sign $\varepsilon\in\{1,-1\}$, the map $P$ is
\begin{equation}\label{eq:PD-integral}
P([v])=[\varepsilon S_0^{-1}v]\qquad\text{in }\coker D.
\end{equation}
Consequently, for every finite abelian coefficient group $B$ with $B^\dagger=\widehat B$, the spherical Poincar\'e pairing is
\begin{equation}\label{eq:PD}
\iota_B(k)(\bar\eta)=\pair{k}{\varepsilon S_0^{-1}\eta},
\qquad k\in\ker D_B,\quad
\bar\eta\in\coker D_{B^\dagger},
\end{equation}
where the evaluation pairing takes values in $\UU(1)$.
\end{proposition}
\begin{proof}
Realize $\Gamma$ in the unit quaternions by
\[
a=\exp(\pi i/m),\qquad b=j,
\]
and take $q=\exp(\pi i/(2m))$, so that $q^2=a$. The identity $j\exp(it)=\exp(-it)j$ gives
\[
qaq^{-1}=a,\qquad qbq^{-1}=ab.
\]
Hence conjugation by $q$ normalizes $\Gamma$ and induces the automorphism
\[
\sigma(a)=a,\qquad \sigma(b)=ab.
\]
The map $F(x)=qxq^{-1}$ satisfies $F(\gamma x)=\sigma(\gamma)F(x)$, so it descends to a based, orientation-preserving diffeomorphism
\[
f\colon M\longrightarrow M,\qquad [x]\longmapsto[qxq^{-1}],
\]
whose action on the fundamental group is $\sigma$.
On $H_1(M,\ZZ)$ the matrix of $f_*$ is therefore $S_0$.

The pullback on $H^2$ is computed using central extensions. Replacing generator lifts $(t,u)$ by $(t,tu)$ changes the relator labels by
\[
(\eta_1,\eta_2)\longmapsto(\eta_1+\eta_2,\eta_2).
\]
Indeed, if $utu^{-1}=\eta_2^{-1}t^{-1}$ in multiplicative notation for the central kernel, then $(tu)^2=\eta_2^{-1}u^2$. The first label is therefore multiplied by $\eta_2$, while the second is unchanged. Thus $f^*$ on $H^2(M,\ZZ)$ also has matrix $S_0$.

Naturality of cap product, together with $f_*[M]=[M]$, gives
\[
f_*\bigl(f^*\eta\cap[M]\bigr)
=\eta\cap f_*[M]
=\eta\cap[M].
\]
In our two coordinate systems this becomes
\[
S_0PS_0=P.
\]

If $m$ is odd, then $\coker D\simeq\ZZ/4$, with $a=2b$ and $\rho_1=2\rho_2$. Any isomorphism between these cyclic groups sends $\rho_2$ to $b$ or $-b$. Since $S_0^{-1}$ acts as $-1$ on $\coker D$, this already proves \eqref{eq:PD-integral} up to the overall orientation sign.

Suppose now that $m$ is even. Then $\coker D\simeq(\ZZ/2)^2$, and $S_0^{-1}=S_0$ on $\coker D$. The displayed naturality identity says that $P$ commutes with $S_0$. An invertible matrix over $\FF_2$ commuting with $S_0$ has the form
\[
P=\begin{pmatrix}1&c\\0&1\end{pmatrix},\qquad c\in\{0,1\}.
\]
In particular,
\[
P(\rho_1)=a,\qquad P(\rho_2)=b+ca.
\]
We determine $c$ by lifting cycles to a double cover.

Write $m=2n$ and put
\[
H=\langle\alpha=a^2,\ \beta=b\rangle\simeq\Dic_n,
\qquad
p\colon N=S^3/H\longrightarrow M.
\]
For $n=1$, the same presentation defines $C_4$. Let $\rho'_1,\rho'_2$ be the classes for the relators $\alpha^n\beta^{-2}$ and $(\beta\alpha\beta^{-1}\alpha)^{-1}$, using the cellular coordinates of Section~\ref{sec:packets}. Restricting a central extension to $H$ leaves the first relator label unchanged and doubles the second, so
\[
p^*(\eta_1,\eta_2)=(\eta_1,2\eta_2),\qquad
p^*\rho_2=2\rho'_2.
\]

Let $\tau\colon H_1(M,\ZZ)\to H_1(N,\ZZ)$ send a cycle to the sum of all its lifts. For $a$, the two lifted paths concatenate to $a^2=\alpha$. For $b$, the two lifted loops, starting in the sheets represented by $1,a$, represent $\beta$ and $aba^{-1}=\alpha\beta$. Hence
\[
\tau(a)=\alpha,\qquad \tau(b)=\alpha+2\beta.
\]
The same sum-of-lifts construction on chains commutes with cap product after pulling back the cohomology class. With the covering orientation on $N$, it sends $[M]$ to $[N]$, and therefore
\begin{equation}\label{eq:PD-cover}
\tau P=P_Np^*,\qquad P_N(\xi)=\xi\cap[N].
\end{equation}

We only need the value of \eqref{eq:PD-cover} on $\rho_2$. The group $H_1(N,\ZZ)$ has order four. If it is $(\ZZ/2)^2$, then $2P_N(\rho'_2)=2\beta=0$. If it is cyclic, $\rho'_2$ and $\beta$ generate the respective cyclic groups of order four, so any isomorphism again gives $2P_N(\rho'_2)=2\beta$. Thus
\[
2\beta+(c+1)\alpha
=\tau(b+ca)
=\tau P(\rho_2)
=P_N(2\rho'_2)
=2\beta.
\]
Since $\alpha$ is a nonzero element of order two, $c=1$. This gives
\[
P(\rho_1)=a,\qquad P(\rho_2)=a+b=b-a,
\]
which is \eqref{eq:PD-integral}. For even $m$ the overall sign acts trivially; reversing orientation changes the sign in the odd case. The matrix formula is well defined on the cokernel, since
\[
S_0^{-1}D=D(S_0^t)^{-1}.
\]

For any finite coefficient group $C$, universal coefficients give natural identifications
\[
H^2(M,C)=H^2(M,\ZZ)\otimes C,
\qquad
H_1(M,C)=H_1(M,\ZZ)\otimes C.
\]
Here $H^3(M,\ZZ)=H_0(M,\ZZ)=\ZZ$, so no additional torsion terms occur. Naturality of cap product therefore extends \eqref{eq:PD-integral} to these coefficients. Evaluating $k\in H^1(M,B)=\Hom(H_1(M,\ZZ),B)$ on the dual cycle of $\bar\eta\in H^2(M,B^\dagger)$ proves \eqref{eq:PD}.
\end{proof}

\begin{corollary}\label{cor:poincare-packets}
In the evaluation coordinates of Section~\ref{sec:packets}, the graded packet multiplicities satisfy
\begin{equation}\label{eq:PDpackets}
n_G(e,\omega;u)
=n_{G^\dagger}(\varepsilon S_0^{-1}\omega,\varepsilon S_0e;u).
\end{equation}
If electric characters are labeled by $E=\varepsilon S_0e$ using \eqref{eq:PD}, this is the transposition $(E,\omega)\mapsto(\omega,E)$ of electric and magnetic cohomology labels.
\end{corollary}
\begin{proof}
The automorphism $\sigma$ of $\Gamma$ sends a pair $(t,u)$ to $(t,tu)$, sending actual grades to $S_0b$ and lift translations to $S_0^tk$, so its action on packet labels is
\[
(e,\omega)\longmapsto(S_0^{-1}e,S_0\omega).
\]
The full pair centralizer is unchanged, since $C_G(t,tu)=C_G(t,u)$. Thus the automorphism preserves polynomial degree and gives this symmetry of packet Hilbert series for either group.

Combining it with coordinate Fourier reciprocity \eqref{eq:packet} and electric inversion symmetry \eqref{eq:realform} yields
\[
\begin{aligned}
n_G(e,\omega;u)
&=n_G(\varepsilon e,\omega;u)\\
&=n_{G^\dagger}(-\omega,\varepsilon e;u)\\
&=n_{G^\dagger}(-S_0^{-1}\omega,\varepsilon S_0e;u)\\
&=n_{G^\dagger}(\varepsilon S_0^{-1}\omega,\varepsilon S_0e;u).
\end{aligned}
\]
The first and last equalities use the real forms, and the middle equalities use Fourier reciprocity and the automorphism $\sigma$. Formula~\eqref{eq:PD} identifies an evaluation character labeled by $e$ with the Poincar\'e label $E=\varepsilon S_0e$, giving the stated transposition.

The symmetry groups of $\cH_B(G)$ and $\cH_{B^\dagger}(G^\dagger)$, defined in Section~\ref{sec:packets}, are identified by
\[
\begin{aligned}
\phi_B\colon A_B\times A_{B^\dagger}&\longrightarrow A_{B^\dagger}\times A_B,\\
(k,\chi)&\longmapsto
\bigl(\varepsilon S_0^t\chi,
      \varepsilon(S_0^t)^{-1}k\bigr).
\end{aligned}
\]
The identity $D S_0^t=S_0D$ shows that $\phi_B$ is well defined and invertible. The character on $A_{B^\dagger}\times A_B$ with labels $(\varepsilon S_0^{-1}\omega,\varepsilon S_0e)$ pulls back under $\phi_B$ to the character on $A_B\times A_{B^\dagger}$ with labels $(e,\omega)$. Thus \eqref{eq:PDpackets} and semisimplicity in each polynomial degree give a degree-preserving $\phi_B$-equivariant isomorphism
\[
\cH_B(G)\simeq\cH_{B^\dagger}(G^\dagger),
\]
whose degree-zero restriction gives the cohomological class-module reciprocity.
\end{proof}

\section{An invariant Yang--Mills proof of local realization}\label{app:yang-mills}

We keep the setup and notation of Section~\ref{sec:realization} ($1\to K\to\Gamma\to H\to1$, $X=\PP^1$, $\Gamma_i=\Stab_\Gamma(p_i)$, $G$ simply connected compact semisimple with Lie algebra $\fg$, $\chi\colon K\to Z(G)$, and $\phi_i\colon\Gamma_i\to G$ with $\phi_i|_K=\chi$ and $d_1+d_2+d_3=d$) and give an alternative proof of existence and ellipticity in Theorem~\ref{thm:local-realization} via the Yang--Mills functional.

Lemma~\ref{lem:equivariant-bundle-construction} gives a smooth $\Gamma$-equivariant principal $G$-bundle $P\to X$ with these fiber actions, which we keep fixed. If $P$ admits an invariant flat connection, parallel transport on the simply connected sphere turns the fiber actions into a single homomorphism.

Give $X$ the unit round metric and $\fg$ an $\Ad(G)$-invariant inner product. Set $V=P\times_{\Ad}\fg$ and let $\mathcal A^\Gamma(P)$ be the affine space of smooth connections preserved by the lift. All norms and integrated pairings below are taken in $L^2(X)$.

\begin{proposition}\label{prop:app-ym-minimum}
The functional $\mathrm{YM}(A)=\tfrac12\|F_A\|^2$ attains its minimum on $\mathcal A^\Gamma(P)$. Every minimizer $A$ satisfies $d_A^*F_A=0$, and its second variation $Q_A(a)$ is nonnegative for every $a\in\Omega^1(X,V)^\Gamma$.
\end{proposition}

\begin{proof}
Averaging a connection over $\Gamma$ gives an element of $\mathcal A^\Gamma(P)$. Every base orbit is finite and hence has codimension two. R{\aa}de's existence theorem for a fixed lift~\cite[Corollary~1]{RadeInvariant} therefore gives a smooth invariant minimizer $A$. For each invariant $a$, the line $A+ta$ lies in $\mathcal A^\Gamma(P)$, so
\[
0=\left.\frac{d}{dt}\right|_{0}\mathrm{YM}(A+ta)   =\langle d_A^*F_A,a\rangle.
\]
The lift preserves $A$ and both metrics, so $d_A^*F_A$ is itself an admissible variation, and taking $a=d_A^*F_A$ gives the equation. The second derivative along the same lines is nonnegative by minimality.
\end{proof}

The dimension hypothesis enters through the invariant operator
\[
D_A=d_A^*\oplus d_A\colon\Omega^1(X,V)^\Gamma \longrightarrow \bigl(\Omega^0(X,V)\oplus\Omega^2(X,V)\bigr)^\Gamma.
\]
Its index can be computed for any invariant connection.

\begin{proposition}\label{prop:app-ym-index}
For every smooth $A\in\mathcal A^\Gamma(P)$, the operator $D_A$, with domain completed in $W^{1,2}$ and target in $L^2$, is Fredholm and
\[
\operatorname{ind}_{\RR}D_A=d-d_1-d_2-d_3.
\]
\end{proposition}

\begin{proof}
The unrestricted operator is elliptic and commutes with $\Gamma$. Averaging a parametrix over $\Gamma$ and restricting to invariant sections gives a parametrix for $D_A$. Elliptic regularity identifies its kernel and its orthogonal cokernel with smooth invariant solutions of $D_Aa=0$ and $D_A^*(u,\eta)=0$, respectively.

Let $F=(V_{\CC},\bar\partial_A)$ be the holomorphic adjoint bundle. The $(0,2)$-curvature vanishes on the curve, and the invariant connection makes $F$ equivariant. Taking the $(0,1)$-part identifies real $V$-valued one-forms with $F$-valued $(0,1)$-forms. In an oriented orthonormal coframe $e^1,e^2$ normal at a point, write $a=a_1e^1+a_2e^2$. Then $a^{0,1}=\tfrac12(a_1+i a_2)(e^1-i e^2)$, so at that point
\[
2\bar\partial_A^*a^{0,1}=-(\nabla_1-i\nabla_2)(a_1+i a_2)=d_A^*a-i* d_Aa.
\]
Here $*$ is the Hodge star and $\nabla_j$ is covariant differentiation in the dual frame. Thus the equivariant real-linear isomorphism $(u,\eta)\mapsto u-i*\eta$ identifies $D_A$ with $2\bar\partial_A^*$. Dolbeault Hodge theory on a curve identifies the complex kernel and cokernel of this last operator on invariants with $H^1(X,F)^\Gamma$ and $H^0(X,F)^\Gamma$, respectively. Passing to real dimensions and applying Lemma~\ref{lem:invariant-riemann-roch} gives
\[
\operatorname{ind}_{\RR}D_A=2\bigl(\dim_{\CC}H^1(X,F)^\Gamma-\dim_{\CC}H^0(X,F)^\Gamma\bigr)=d-d_1-d_2-d_3.\qedhere
\]
\end{proof}

Now fix a minimizer as in Proposition~\ref{prop:app-ym-minimum}. When the index is zero, vanishing of $\ker D_A$ implies vanishing of $\ker D_A^*$, which contains the curvature by the Yang--Mills equation.

\begin{proposition}\label{prop:app-ym-flat}
If $d_1+d_2+d_3=d$, then every invariant minimizer $A$ satisfies $\ker D_A=\ker D_A^*=0$ on invariant sections. In particular $F_A=0$.
\end{proposition}

\begin{proof}
For a real invariant one-form $a$, set $[a\wedge a](v,w)=2[a(v),a(w)]$. Expanding $F_{A+ta}=F_A+t\,d_Aa+\tfrac12t^2[a\wedge a]$ gives
\[
Q_A(a)=\|d_Aa\|^2+\langle F_A,[a\wedge a]\rangle \ge 0.
\]
The same curvature term occurs in the integrated Weitzenb\"ock formula for $V$-valued one-forms~\cite{BourguignonLawson}. Indeed, if $F_{ji}=F_A(e_j,e_i)$ in an oriented orthonormal frame, its pointwise contraction with $a$ is
\[
\sum_{i,j}\langle[F_{ji},a_j],a_i\rangle=2\langle F_{12},[a_1,a_2]\rangle.
\]
Since the Ricci tensor of the unit sphere equals its metric, the formula is
\[
\|d_Aa\|^2+\|d_A^*a\|^2=\|\nabla_Aa\|^2+\|a\|^2+\langle F_A,[a\wedge a]\rangle.
\]
For $a\in\ker D_A$, both $d_Aa$ and $d_A^*a$ vanish. Substituting into the second variation yields
\[
0\le Q_A(a)=-\|\nabla_Aa\|^2-\|a\|^2,
\]
hence $a=0$. Proposition~\ref{prop:app-ym-index} and the dimension hypothesis now give $\ker D_A^*=0$. On invariant pairs the formal adjoint is $D_A^*(u,\eta)=d_Au+d_A^*\eta$, and since $F_A$ is invariant, Proposition~\ref{prop:app-ym-minimum} gives $D_A^*(0,F_A)=0$, so $F_A=0$.
\end{proof}

It remains to read off the representation in a parallel frame.

\begin{proof}[Proof of existence and ellipticity in Theorem~\ref{thm:local-realization}]
Choose an invariant minimizer $A$, which is flat by Proposition~\ref{prop:app-ym-flat}. Since $X$ is simply connected, a point of one fiber extends by parallel transport to a global parallel section $s$ of $P$. For each $\gamma\in\Gamma$, the lift takes parallel frames to parallel frames. Hence there is a single $\rho(\gamma)\in G$ such that
\[
\gamma s(x)=s(\gamma x)\rho(\gamma)\qquad(x\in X).
\]
Applying this identity successively to two elements shows that $\rho$ is a homomorphism. For $k\in K$, the action on $P$ is right multiplication by $\chi(k)$, so $\rho(k)=\chi(k)$. At $p_i$, write the frame realizing $\phi_i$ as $s(p_i)g_i$. Computing its stabilizer action in the two frames gives $\rho(\gamma)=g_i\phi_i(\gamma)g_i^{-1}$ for $\gamma\in\Gamma_i$. This proves existence with the prescribed local data.

In the frame $s$, each $v\in\fg^{\rho(\Gamma)}$ gives an invariant parallel section $u$ of $V$. Then $D_A^*(u,0)=d_Au=0$. By Proposition~\ref{prop:app-ym-flat}, $u=0$ and hence $v=0$, so $\rho$ is elliptic. Uniqueness up to $G$-conjugacy then follows from the quasi-Hamiltonian argument in Section~\ref{sec:realization}.

\end{proof}

\section{The cyclic case}\label{sec:cyclic}

The cyclic counting equality is due to Djokovi\'c~\cite[Theorem~1]{djokovic}. Kojima and Tachikawa~\cite[Theorems~2.1 and~2.2]{kojima-tachikawa} give the finite-torus argument and prove its central-coefficient refinement. We extend this argument to full-centralizer polynomial modules for connected compact (and hence, by Lemma~\ref{lem:compact-reduction}, connected complex) reductive groups.

Let $C_n=\langle c\mid c^n=1\rangle$, with $n\ge1$, and let $B\subset Z(G)$ be finite. As in Part~\ref{part:dihedral}, we write tori and finite coefficient groups additively, and set $G^\dagger=(G/B)^\vee$ and $B^\dagger=\widehat B$. Put
\begin{equation}\label{eq:cyclic-classes}
\cR_B(C_n,G)=\{t\in G:t^n\in B\}/G,
\qquad \gamma(t)=t^n.
\end{equation}
For $[t]\in\cR_B(C_n,G)$, choose a maximal torus $S_t\subset C_G(t)^0$ and let $W_t$ be the image of $N_{C_G(t)}(S_t)$ on $\mathfrak s_t=\Lie(S_t)_{\CC}$. Define
\begin{equation}\label{eq:cyclic-classes-mods}
\cP_B(C_n,G)=\bigoplus_{[t]\in\cR_B(C_n,G)}
\Sym(\mathfrak s_t^*)^{W_t},
\end{equation}
with linear forms of degree one. As in \eqref{eq:gradeddef}, conjugation and changes of maximal torus canonically identify the polynomial fibers. Central multiplication preserves $C_G(t)$, so $\tau_b(t)=bt$, for $b\in B$, defines translations on this module. Set
\begin{equation}\label{eq:cyclic-operators}
(T_bf)(t)=f(\tau_{-b}t),\qquad
(M_\chi f)(t)=\chi(t^n)f(t),\qquad \chi\in B^\dagger.
\end{equation}

\begin{proposition}\label{prop:cyclic-fourier} For every $n\ge1$, every connected compact reductive group $G$, and every finite $B\subset Z(G)$, there is a complex-linear isomorphism preserving polynomial degree,
\[
\mathcal F^{\mathrm{cyc}}_B\colon\cP_B(C_n,G)\xrightarrow{\ \sim\ }
\cP_{B^\dagger}(C_n,G^\dagger),
\]
such that
\begin{equation}\label{eq:cyclic-fourier-rules}
\mathcal F^{\mathrm{cyc}}_B T_b=M_b^\dagger\mathcal F^{\mathrm{cyc}}_B,
\qquad
\mathcal F^{\mathrm{cyc}}_B M_\chi=T_{-\chi}^\dagger\mathcal F^{\mathrm{cyc}}_B.
\end{equation}
Its degree-zero restriction compares the function spaces on the actual class sets. For $B=1$, it gives $|\Hom(C_n,G)/G|=|\Hom(C_n,G^\vee)/G^\vee|$.
\end{proposition}
\begin{proof}
Fix a maximal torus $T=V/\Lambda$ of $G$, with Weyl group $W=N_G(T)/T$. Every $t$ in \eqref{eq:cyclic-classes} is conjugate into $T$, with conjugacy in $T$ given by $W$. Using \eqref{eq:coefficientlattice}, the torus parameters form the finite group
\[
E_{B,n}=\frac{n^{-1}\Lambda_B}{\Lambda},\qquad
\gamma(x)=nx\pmod\Lambda.
\]
The torus $T$ is maximal in $C_G(t)^0$, and
\[
N_{C_G(t)}(T)/T=\Stab_W(t).
\]
Indeed, a representative in $N_G(T)$ centralizes $t$ if and only if its Weyl element fixes $t$. Evaluation at orbit representatives now gives
\begin{equation}\label{eq:cyclic-model}
\cP_B(C_n,G)\simeq
\bigl[\Fun(E_{B,n},\CC)\otimes\Sym(V_{\CC}^*)\bigr]^W.
\end{equation}
The polynomial fiber over an orbit represented by $x$ is $\Sym(V_{\CC}^*)^{\Stab_W(x)}$.

By \eqref{eq:duallattices}, the dual torus is $V^*/\Lambda_B^*$ and its coefficient lattice is $\Lambda^*$. Thus its finite parameter group is
\[
E_{B,n}^\dagger=\frac{n^{-1}\Lambda^*}{\Lambda_B^*},\qquad
\gamma^\dagger(y)=ny\pmod{\Lambda_B^*}.
\]
The pairing
\begin{equation}\label{eq:cyclic-pairing}
P_n(x,y)=\eexp{n\pair{\widetilde y}{\widetilde x}}
\end{equation}
is well defined and perfect, since multiplication by $n$ identifies the second group with $\Lambda^*/n\Lambda_B^*$, the character group of $n^{-1}\Lambda_B/\Lambda$. The pairing is $W$-invariant for the contragredient actions. Moreover,
\[
\frac{P_n(x+b,y)}{P_n(x,y)}=\gamma^\dagger(y)(b),\qquad
\frac{P_n(x,y+\chi)}{P_n(x,y)}=\chi(\gamma(x)).
\]

A $W$-invariant positive definite form on $V$ induces a degree-preserving $W$-equivariant isomorphism $J\colon\Sym(V_{\CC}^*)\xrightarrow{\sim}\Sym(V_{\CC})$. For polynomial-valued functions define
\[
(\mathcal F^{\mathrm{cyc}}_Bf)(y)=
\frac{1}{\sqrt{|E_{B,n}|}}\sum_{x\in E_{B,n}}P_n(x,y)J(f(x)).
\]
Finite Fourier inversion shows that $\mathcal F^{\mathrm{cyc}}_B$ is an isomorphism; it commutes with $W$ and preserves degree, so \eqref{eq:cyclic-model} transports it to the stated modules. The two pairing identities give \eqref{eq:cyclic-fourier-rules} with the conventions \eqref{eq:cyclic-operators}, and restricting to constants gives the degree-zero assertion.
\end{proof}

For the cohomological formulation, put
\begin{equation}\label{eq:cyclic-cohomology}
A_B=B[n]=H^1(C_n,B),\qquad Q_B=B/nB=H^2(C_n,B).
\end{equation}
An actual grade $b=t^n$ represents $\omega=[b]\in Q_B$, and changing the lift $t$ by $k\in B$ changes $b$ by $nk$. Write $\cP_b(C_n,G)$ for the grade-$b$ summand, choose $b_\omega$ for each $\omega\in Q_B$, and define
\begin{equation}\label{eq:cyclic-H}
\cH_B^{\mathrm{cyc}}(G)=\bigoplus_{\omega\in Q_B}\cP_{b_\omega}(C_n,G).
\end{equation}
The operators $T_k$, $k\in A_B$, and $M_\chi$, $\chi\in A_{B^\dagger}$, commute. Evaluation identifies $\widehat{A_B}=Q_{B^\dagger}$ and $\widehat{Q_B}=A_{B^\dagger}$. Let $n_G^{\mathrm{cyc}}(e,\omega;u)$ be the Hilbert series of the $e$-eigenspace of $A_B$ in $\cP_{b_\omega}(C_n,G)$. Translations identify the spaces for different representatives of $\omega$ and commute with $A_B$, so this series is well defined.

To specify the spherical coordinates, let $M_n=S^3/C_n$ with its quotient orientation. The classifying map identifies cohomology in degrees one and two. Let $\rho\in H^2(M_n,\ZZ)$ be the generator given by the relator $c^n$, and also write $c$ for the generator of $H_1(M_n,\ZZ)$. Both groups are cyclic of order $n$, so
\[
\rho\cap[M_n]=r c
\]
for an integer $r$ prime to $n$. Naturality with respect to coefficients gives $\pair{\omega}{\chi}_{\mathrm{PD}}=\pair{r\omega}{\chi}$. Here multiplication by $r$ is invertible on every $Q_B$ and $A_B$; its inverse is denoted by $r^{-1}$.

\begin{corollary}\label{cor:cyclic-packets} There is a degree-preserving isomorphism
\[
\cH_B^{\mathrm{cyc}}(G)\simeq
\cH_{B^\dagger}^{\mathrm{cyc}}(G^\dagger)
\]
equivariant under
\[
\phi_B^{\mathrm{cyc}}\colon A_B\times A_{B^\dagger}\xrightarrow{\sim}
A_{B^\dagger}\times A_B,\qquad
(k,\chi)\longmapsto(r^{-1}\chi,rk).
\]
The corresponding packet identity is
\begin{equation}\label{eq:cyclic-PDpackets}
n_G^{\mathrm{cyc}}(e,\omega;u)=
n_{G^\dagger}^{\mathrm{cyc}}(r\omega,r^{-1}e;u).
\end{equation}
With electric characters labeled by $E=r^{-1}e$ using the spherical pairing, this is the transposition $(E,\omega)\mapsto(\omega,E)$. Its degree-zero restriction gives the cyclic central-coefficient refinement.
\end{corollary}
\begin{proof}
Applying Proposition~\ref{prop:cyclic-fourier} to the commuting coefficient actions gives
\begin{equation}\label{eq:cyclic-packet}
n_G^{\mathrm{cyc}}(e,\omega;u)=n_{G^\dagger}^{\mathrm{cyc}}(-\omega,e;u).
\end{equation}
As in the proof of Corollary~\ref{cor:packets}, Fourier first compares all actual grades in a fixed cohomology class. They contribute $|nB|$ equal copies, while the dual side contributes $|nB^\dagger|=|nB|$ copies. Cancellation in each degree proves \eqref{eq:cyclic-packet}. Each fixed-grade module has a real form, so $n_G^{\mathrm{cyc}}(e,\omega;u)=n_G^{\mathrm{cyc}}(-e,\omega;u)$.

Choose an integer $a\equiv r^{-1}\pmod n$ prime to $n|B|$; such a lift exists by the Chinese remainder theorem. The power map $t\mapsto t^a$ is invertible on the elements under consideration and satisfies $C_G(t^a)=C_G(t)$. It sends grades to $ab$ and coefficient twists to $ak$, hence gives the degree-preserving packet symmetry
\[
(e,\omega)\longmapsto(a^{-1}e,a\omega).
\]
Apply this symmetry on the dual side of \eqref{eq:cyclic-packet}, then use its real form, to obtain
\[
n_G^{\mathrm{cyc}}(e,\omega;u)
=n_{G^\dagger}^{\mathrm{cyc}}(-r\omega,r^{-1}e;u)
=n_{G^\dagger}^{\mathrm{cyc}}(r\omega,r^{-1}e;u).
\]
The character on $A_{B^\dagger}\times A_B$ with labels $(r\omega,r^{-1}e)$ pulls back under $\phi_B^{\mathrm{cyc}}$ to the character on $A_B\times A_{B^\dagger}$ with labels $(e,\omega)$, so semisimplicity in each polynomial degree gives the module isomorphism. Finally, the formula for the spherical pairing identifies $\phi_B^{\mathrm{cyc}}$ with the Poincar\'e swap of Kojima and Tachikawa~\cite[Definition~1.8]{kojima-tachikawa}.
\end{proof}

\clearpage
\bookmarksetup{startatroot}
\addtocontents{toc}{\protect\vspace{12pt plus 2pt}}
\nocite{Husemoller}
\bibliographystyle{amsplain}
\bibliography{references}
\end{document}